\documentclass[10pt]{article}
\usepackage{bm}
\usepackage[top=0.9in,bottom=1.1in,left=1.05in,right=1.05in]{geometry}

\usepackage{textgreek,mathtools}
\usepackage{amsmath}
\usepackage{graphicx}
\usepackage[colorinlistoftodos]{todonotes}
\usepackage[colorlinks=true, allcolors=blue]{hyperref}
\usepackage{amsmath}
\usepackage{upgreek}
\usepackage{amssymb}
\usepackage{amsfonts}
\usepackage{amsthm}
\usepackage{mathrsfs}
\usepackage{fontenc}
\usepackage{empheq}
\usepackage{enumerate}
\usepackage{comment}
\usepackage{appendix}
\usepackage{mathtools}
\usepackage{bbm}
\usepackage{subfig}
\usepackage{float}
\usepackage{dsfont}
\usepackage{longtable}
\mathtoolsset{showonlyrefs}

\numberwithin{equation}{section}

\DeclareMathOperator*{\argmin}{arg\,min} 
\DeclareMathOperator*{\argmax}{arg\,max}

\newcommand{\F}{\mathcal{F}}

\newcommand{\R}{\mathbb{R}}

\newcommand{\bX}{\bm{X}}
\newcommand{\bY}{\bm{Y}}
\newcommand{\bZ}{\bm{Z}}

\newcommand{\eps}{\epsilon}
\newcommand{\ud}{\,\mathrm{d}}

\newcommand{\EE}{\mathbb{E}}

\newcommand{\PP}{\mathbb{P}}

\newcommand{\RR}{\mathbb{R}}
\newcommand{\mcF}{\mathcal{F}}

\newcommand{\mc}[1]{\mathcal{#1}}
\newcommand{\half}{\frac{1}{2}}
\newcommand{\mt}[1]{\mathtt{#1}}

\theoremstyle{plain}

\newtheorem{theo}{Theorem}[section]

\newtheorem{remark}[theo]{Remark}

\newtheorem{rem}[theo]{Remark}

\newtheorem{defn}[theo]{Definition}
\newtheorem{defi}[theo]{Definition}

\newcommand{\bx}{\bm{x}}
\newcommand{\bz}{\bm{z}}
\newcommand{\bW}{\bm{W}}
\newcommand{\balpha}{\bm{\alpha}}

\newcommand{\MCF}{\mathcal{F}}
\newcommand{\cN}{\mathcal{N}}

\newcommand{\abs}[1]{\left|#1\right|}

\newcommand{\transpose}{^{\operatorname{T}}}

\makeatletter
\let\@fnsymbol\@arabic
\makeatother
\title{Recent Developments in Machine Learning Methods 
\\ for Stochastic Control and Games}
\author{Ruimeng Hu \thanks{Department of Mathematics, and Department of Statistics and Applied Probability, University of California, Santa Barbara, CA 93106-3080, USA, {\em rhu@ucsb.edu}.} \and Mathieu Lauri\`{e}re \thanks{Shanghai Frontiers Science Center of Artificial Intelligence and Deep Learning; NYU-ECNU Institute of Mathematical Sciences, NYU Shanghai, 567 West Yangsi Road, Shanghai, 200126, People’s Republic of China, {\em mathieu.lauriere@nyu.edu}. Corresponding author.}}
\date{}%

\newcommand{\mbA}{\mathbb{A}}
\newcommand{\NN}{\mathbb{N}}

\newcommand{\ltwonorm}[1]{\left\lVert#1\right\rVert_2}

\newcommand{\figwidth}{{0.95\textwidth}}

\usepackage{mathtools,empheq,algorithm,algorithmic}
\usepackage{etoolbox,booktabs}
\newcommand{\algorithmicdoinparallel}{\textbf{do in parallel}}
\makeatletter
\AtBeginEnvironment{algorithmic}{%
  \newcommand{\FORALLP}[2][default]{\ALC@it\algorithmicforall\ #2\ %
    \algorithmicdoinparallel\ALC@com{#1}\begin{ALC@for}}%
}
\makeatother

\newcommand{\Tr}{\mathrm{Tr}}
\newcommand{\Hess}{\mathrm{Hess}}

\newcommand{\ctrl}{\alpha}
\newcommand{\ctrldim}{k}
\newcommand{\dom}{\mathcal{Q}}

\newcommand{\cA}{\mathcal{A}}
\newcommand{\cE}{\mathcal{E}}

\newcommand{\cP}{\mathcal{P}}

\newcommand{\cU}{\mathcal{U}}

\newcommand{\bN}{\mathbf{N}}

\newcommand{\grad}{\nabla}
\newcommand{\indic}{\mathbf{1}}

\newcommand{\diver}{\mathrm{div}}
\newcommand{\Law}{\mathcal{L}}

\usepackage{titlesec}

\titleformat{\paragraph}
{\normalfont\normalsize\bfseries}{\theparagraph}{1em}{}
\titlespacing*{\paragraph}
{0pt}{3.25ex plus 1ex minus .2ex}{1.5ex plus .2ex}

\newcommand{\checkX}{\check{X}}

\begin{document}

\maketitle
\begin{abstract}
Stochastic optimal control and games have a wide range of applications, from finance and economics to social sciences, robotics, and energy management. Many real-world applications involve complex models that have driven the development of sophisticated numerical methods. Recently, computational methods based on machine learning have been developed for solving stochastic control problems and games. In this review, we focus on deep learning methods that have unlocked the possibility of solving such problems, even in high dimensions or when the structure is very complex, beyond what traditional numerical methods can achieve. We consider mostly the continuous time and continuous space setting. Many of the new approaches build on recent neural-network-based methods for solving high-dimensional partial differential equations or backward stochastic differential equations, or on model-free reinforcement learning for Markov decision processes that have led to breakthrough results. This paper provides an introduction to these methods and summarizes the state-of-the-art works at the crossroad of machine learning and stochastic control and games.

\end{abstract}

\tableofcontents

\section{Introduction}

Stochastic optimal control and games have been extensively studied throughout the twentieth century and have found a wide range of applications in various areas such as finance, social sciences, operations research, and epidemic management problems, among others.
In recent years, computational methods for stochastic control and games have made great progress with the help of machine learning tools. A striking example of recent breakthroughs in applied mathematics using such tools is the numerical resolution of general nonlinear parabolic partial differential equations and backward stochastic differential equations in high dimensions~\cite{hutzenthaler2019multilevel, MR3736669,HaJeE:18,SiSp:18}. In short, stochastic control problems study how an agent optimally controls a stochastic dynamical system. The agent perceives some observations of the system's state and can decide to influence the evolution of the state based on these observations. The goal is to optimize an objective function that typically incorporates the cost of controlling the system and the reward for reaching some state. One of the most popular methods to solve such problems is dynamic programming, developed by Richard Bellman in the 1950s~\cite{bellman1957markovian}. However, this method suffers from what Bellman called the \emph{curse of dimensionality}, meaning that its complexity increases drastically with the number of possible states. This is a significant issue for systems that  evolve in continuous and high-dimensional spaces, since they can not be approximated by a small number of states. In such cases, using exact dynamic programming becomes computationally infeasible. 
Additionally, complexity can also arise from the structure of the system's evolution.
For example, in some cases, the system's evolution or its observation may be subject to delay, which appears in many real-world applications, {\it e.g.}, in economics, mechanics, or biology.  
To model the delay feature, the dynamics of the controlled system will depend not only on the current state but also on the history prior to the current time. This makes the problem path-dependent and, therefore, infinite-dimensional.

On the other hand, stochastic differential game theory, which was initiated by \cite{Isaacs1965}, combines theory and optimal control, and provides a framework for modeling and analyzing the behavior of strategic agents in the context of a dynamical system. The theory has been extensively employed in many disciplines, including management science,  economics, social science, and biology. One of the core objectives in differential games is to compute Nash equilibria, \textit{i.e.}, strategy profiles according to which no player has an incentive to deviate unilaterally~\cite{nash1951noncoop}. However, computing Nash equilibria in $N$-agent games is a notoriously hard problem, and direct computation of Nash equilibria is extremely demanding in terms of time and memory \cite{DaGoPa:2009} even for moderately large $N$. The mean-field game paradigm has been introduced independently by Lasry and Lions in \cite{LaLi:2007} and by Huang, Malham\'{e} and Caines in \cite{HuMaCa:06} to provide a tractable approximation for games with very large populations. A mean-field game is a game with a continuum of infinitesimal agents, where any single agent does not influence the rest of the population, which forms the mean field with which each agent interacts. This framework provides an efficient way to compute approximate Nash equilibria for symmetric $N$-agent games when $N$ is large. However, challenges remain in terms of computational complexity for games with high dimensional or complex environments, or when common noise affects the population dynamics. Furthermore, in the intermediate regime when $N$ is moderately large, the mean-field theory does not provide a good approximation of the $N$-player game. In such cases, one may still face a high-dimensional problem. To make these challenges more concrete, we discuss some illustrative examples.

\subsection{Some high-dimensional examples in applications}
We first discuss an example in financial markets. Many problems in economics or finance involve multiple interacting agents. For instance, we may consider a group of traders who buy and sell stocks in a financial market such as the S\&P 500, a free-float weighted measurement stock market index of 500 of the largest companies listed on stock exchanges in the United States. Each trader's portfolio describes the investment in stocks available on the market. To describe the investments of the whole group of traders, we need to incorporate all the traders' portfolios, and hence this description can be very high dimensional. However, if the traders have similar risk preferences, it is sufficient to study how one representative trader optimizes their payoff to understand the whole group's behavior. Any single agent has only a negligible impact on the stocks' prices. However, the impact of the group might be significant because if everyone wants to buy or sell the same stock, then the price will probably be shifted up or down. The portfolio optimization problem then becomes to find a Nash equilibrium. Each agent can anticipate that every other agent will behave like themselves and can thus predict the impact of the group on the stocks' prices. The mean-field game paradigm provides a rigorous framework in which each agent, taken individually, has no impact at all on the group's dynamics, and the problem is to find a fixed point at the population level. This simplifies the analysis. A first approach to describe the solution is through a forward-backward system of partial differential equations (PDEs) composed of a Fokker-Plank (FP) equation for the population distribution and a Hamilton-Jacobi-Bellman (HJB) equation for the value function of an infinitesimal player. However, several difficulties arise. Firstly, even solving the optimal control problem for the representative agent can be challenging: Given the number of stocks, the problem is in high dimension, and using an exact dynamic programming algorithm can be computationally too expensive. Secondly, even if the agent's portfolio state is in a low dimension, the optimal strategy may depend on the average state or control of the group. If the agents are trading the same stocks, then their trading strategies are all subject to the same source of randomness, which implies that the group's average strategy itself is stochastic. In the context of mean-field games, this is formalized through the notion of common noise. In such cases, the PDEs of the forward-backward system characterizing the solutions are no longer deterministic but stochastic {\cite{peng1992stochastichjb,carmonadelarue2014mastereqlarge}, making the system considerably harder to solve. Another approach \cite[Chapter~2]{CaDe2:17} characterizes the solution by a forward-backward system of stochastic differential equations (FBSDEs) which involves the conditional distribution of the forward and the backward processes, given the common noise. A third approach consists of describing the solution by a PDE on the space of probability distributions \cite{cardaliaguetdelaruelasrylions2019master}, but here again, this PDE is difficult to solve since it is posed on a high or infinite-dimensional space. In any case, in this application to optimal trading and as in many others, realistic models lead to the need to approximate functions of high or infinite-dimensional inputs. Altogether, this makes it crucial to develop efficient and accurate deep learning algorithms and theories for computing optimal controls and Nash equilibria in high dimensions.

A second example arises in infectious disease control of multiple regions. In a classic compartmental epidemiological model, each individual in a geographical region is assigned a label, for instance, \textbf{S}usceptible, \textbf{E}xposed, \textbf{I}nfectious, \textbf{R}emoved, \textbf{V}accinated. The transmission of a virus, being infected or recovered, moves individuals from one compartment to another, and this transition is usually described by stochastic dynamical equations. When a disease outbreak is reported, the region planner needs to take measurements to control its spread. The ongoing COVID-19 includes issuing lockdown or work from home policies, developing vaccines and later expanding equitable vaccine distribution, providing telehealth programs, distributing free testing kits, educating the public on how the virus transits, and focusing on surface disinfection. As a region planner, the decisions are usually made by weighting different costs, including the economic loss due to less productivity during lockdown policy or work from home policy, the economic value of life due to the death of infected individuals, and various social-welfare costs due to measurements mentioned above, and many more. Moreover, as the world is more interconnected than ever before, one region's decision will inevitably influence its neighboring regions. For instance, in the US, the decision made by the New York governor will affect the situation in New Jersey as so many people travel daily between the two states. Imagining that both state governors make decisions representing their own benefits but also take into account others' rational decisions, and they may even compete for scarce resources ({\it e.g.}, frontline workers and personal protective equipment), these are precisely the features of a non-cooperative game. A Nash equilibrium computed from such a game will definitely provide some qualitative guidance and insights for policymakers on the impact of certain policies. However, even with only three states (New York, New Jersey, and Pennsylvania) and a simple stochastic SEIR model as in \cite{xuan2020optimal,xuan2020ams}, this problem's state space is already twelve dimensions. Figure~\ref{fig:pandemic} below showcases the equilibrium lockdown policy corresponding to the multi-region SEIR model solved by a deep learning algorithm proposed in \cite{HaHu:19} (see Section~\ref{sec:MarkovianNE}) between the three states. The model parameters are estimated from real data posted by the Centers for Disease Control and Prevention (CDC). In general, the problem dimension is proportional to the number of compartments in the epidemiological model multiplied by the number of regions considered. For the most basic SIR model, the dimension of the problem for US governors will be $3 \times 50 = 150$.

\begin{figure}[!htb]
    \centering
    \includegraphics[width = 0.75\textwidth, trim = {2em 2em 5em 6em}, clip, keepaspectratio=True]{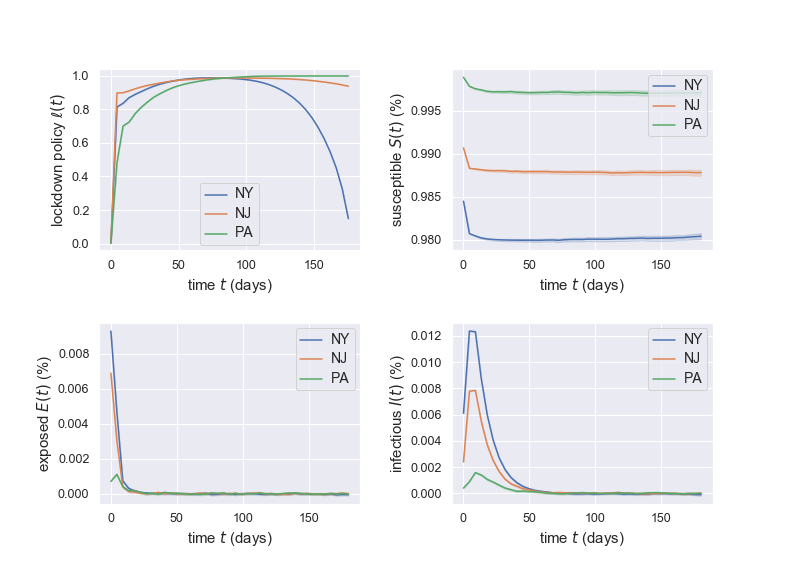}
    \caption{A case study of the COVID-19 pandemic in three states: New York (NY), New Jersey (NJ), and Pennsylvania (PA) in \cite{xuan2020optimal}. Plots of optimal policies (top-left), Susceptibles (top-right), Exposed (bottom-left), and Infectious (bottom-right) for three states: New York (blue), New Jersey (orange), and Pennsylvania (green). Large $\ell$ indicates high intensity of lockdown policy. Choices of parameters are referred to \cite[Section 4.2]{xuan2020optimal}. %
    }
    \label{fig:pandemic}
\end{figure}

\subsection{An illustrative linear quadratic model}
\label{sec:intro-LQsysrisk}
To illustrate numerical methods and show that they can correctly compute the problem's solution, it is convenient to have examples with analytical or semi-explicit solutions. We present here an example introduced in \cite{CaFoSu:15} to model the interactions in a system of banks. This model and other similar models with linear-quadratic structures admit a closed-form solution and have found applications in various fields.

We consider a stochastic differential game with $N$ players, and we denote by $\mathcal{I} =  \{1, 2, \ldots, N\}$ the set of players. Each player is interpreted as a bank and the state is its log-reserve.  Let $T$ be a finite time horizon. At each time $t \in [0,T]$,  player $i \in \mathcal{I}$ has a state $X_t^i \in \RR$ and takes an action $\alpha_t^i \in \RR$. The information structure will be discussed later and for now, we proceed informally but we can think of $\alpha^i_t$ as a stochastic process adapted to a filtration that represents the information available to player $i$. The dynamics of the controlled state process on $[0,T]$ are given by
$$
    \ud X_t^i = [a(\overline{X}_t - X_t^i)   + \alpha_t^i ] \ud t  + \sigma \left(\rho \ud W_t^0 + \sqrt{1-\rho^2} \ud W_t^i\right), \quad X_0^i \sim \mu_0, \quad i \in \mc{I}, \quad \overline{X}_t = \frac{1}{N}\sum_{i=1}^N X_t^i,
$$
where $\mu_0$ is a given initial distribution and $\bm{W} =[W^0, W^1,\ldots, W^N]$ are $(N+1)$ $m$-dimensional independent Brownian motions. %
We shall call $W^i$ the idiosyncratic (i.e., individual) noises and $W^0$ the common noise. The parameter $\rho \in[0,1]$ characterizes the noise correlation between agents. Here $a(\overline{X}_t - X_t^i)$ represents the rate at which bank $i$ borrows from or lends to other banks in the lending market, while $\alpha_t^i$  denotes its control rate of cash flows to a central bank. %
Furthermore, $\overline{X}_t = \frac{1}{N} \sum_{i=1}^N X^i_t$ denotes the average state. The $N$ dynamics are thus coupled since all the states $\bm{X}_t = [X_t^1, \ldots, X_t^N]$ affect the drift of every agent. Given a set of strategies $(\bm{\alpha}_t)_{t \in [0,T]} = ([\alpha_t^1, \ldots, \alpha_t^N])_{t \in [0,T]}$, the cost incurred to player $i$ is
\begin{equation}%
    J^i(\bm{\alpha}) = \EE\left[\int_0^T f^i(t, \bm X_t, \bm\alpha_t) \ud t + g^i(\bm X_T)\right],
\end{equation}
where the running cost $f^i: [0,T] \times \RR^{N}\times \RR^N \to \RR$ and the terminal cost $g^i: \RR^{N} \to \RR$ are given by
\begin{equation}
    f^i(t, \bm{x},\balpha) = \half (\alpha^i)^2 - q \alpha^i(\overline{x} - x^i) + \frac{\eps}{2}(\overline{x} - x^i)^2,  \quad g^i(\bm{x}) = \frac{c}{2}(\overline{x} - x^i)^2, \quad \overline{x} = \frac{1}{N} \sum_{i=1}^N x^i. 
\end{equation}
where $\bm{x} = [x^1, \ldots, x^N]$ and $\bm{\alpha} = [\alpha^1, \ldots, \alpha^N]$. 
All the parameters are non-negative. Here $\half(\ctrl^i)^2$ denotes the quadratic cost of the control, and $-q\ctrl^i(\overline{x} - x^i)$ models the incentive to borrowing or lending: Bank $i$ will want to borrow if $X_t^i$ is smaller than $\overline{X}_t$ and lend if $X_t^i$ is larger than $\overline{X}_t$. The quadratic term $(\overline{x} - x^i)^2$ in $f^i$ and $g^i$ penalizes the deviation from the average, given the other players' states. Player $i$ chooses $(\alpha_t^i)_{t \in [0,T]}$ to minimize her cost $J^i(\balpha)$ within some set of admissible strategies. We assume $q\le \epsilon^2$ so that the Hamiltonian is jointly convex in state and control variables, ensuring that there is at most one best response and then at most one Nash equilibrium. In the original work \cite[Section~3.1]{CaFoSu:15}, open-loop and closed-loop equilibria are characterized by semi-explicit formulas using ordinary differential equations.

As the number of agents $N$ grows to infinity, the idiosyncratic noises have a smaller and smaller influence on $\overline{X}$, which, in the limit, depends only on the common noise $W^0$. This is formalized in the following mean-field game (MFG). Let $(W_t)_{0 \leq t \leq T}$ and $(W_t^0)_{0 \leq t \leq T}$  be independent $m$-dimensional Brownian motions. We shall refer to $W$ as the {idiosyncratic noise} of the representative player and to $W^0$ as the {common noise} of the system. We consider the stochastic control problem
\begin{equation*} \inf_{\alpha} \EE\biggl\{\int_0^T \left[
    \frac{\alpha_t^2}{2}-q\alpha_t(m_t-X_t)+\frac{\epsilon}{2}(m_t-X_t)^2
    \right]\ud t +\frac{c}{2}(m_T-X_T)^2 \biggr\}, \end{equation*}
 \begin{equation}\label{def:LQ_SDE} \text{where }\quad\displaystyle \ud X_t = [a(m_t-X_t)+\alpha_t]\ud t  + \sigma(\rho \ud W_t^0 + \sqrt{1-\rho^2}\ud W_t), \quad X_0 \sim \mu_0,\end{equation} and the representative agent controls her state $X$ through a control process $\alpha$. Here $m_t = \EE[X_t|\mcF^{W^0}_t]$ is the conditional population mean given the common noise. As in the $N$-player case, one advantage of LQ models lies in the existence of an analytical solution for the mean-field equilibrium, which can provide a benchmark to test numerical algorithms. In this model, at equilibrium, we have
\begin{align}
    & m_t = \EE[X_0] + \rho\sigma W_t^0, \quad t\in[0,T], \label{def:LQ_m} \\
    & \alpha_t = (q+\eta_t) (m_t -X_t), 
    \quad t\in[0,T], \label{def:LQ_alpha}
\end{align}
where $\eta$ is a deterministic function of time solving the Riccati equation,
\begin{equation*}
    \dot \eta_t = 2(a+q)\eta_t + \eta_t^2 - (\eps - q^2), \quad \eta_T = c.
\end{equation*}
The solution is given by
\begin{equation*}
    \eta_t = \frac{-(\epsilon-q^2)(e^{(\delta^+-\delta^-)(T-t)}-1) -c(\delta^+e^{(\delta^+-\delta^-)(T-t)}-\delta^-)}{(\delta^-e^{(\delta^+-\delta^-)(T-t)} - \delta^+) -c(e^{(\delta^+-\delta^-)(T-t)} -1)},
\end{equation*}
where $\delta^\pm = -(a+q)\pm \sqrt{R}$, $R = (a+q)^2 + (\epsilon-q^2)>0$. At equilibrium, \textit{i.e.}, when all the payers use the equilibrium control, and the minimal expected cost for a representative player is 
\begin{equation*}
u(0, x_0 - \EE[x_0]),\quad\text{with}\quad    u(t,x) = \frac{\eta_t}{2}x^2 + \half \sigma^2(1-\rho^2)\int_t^T \eta_s \ud s.
\end{equation*}

Even though the state is in dimension one only, the presence of the common noise means that the optimal control is a function of the common noise. In this example, the equilibrium control actually depends on the common noise only through the first conditional moment of the distribution. In this case, the equilibrium can be found using neural networks which take as inputs not only the individual player's state but also an estimate of this first conditional moment. This idea can be extended to scenarios where the dependence on the common noise occurs only through a finite dimensional vector of information; see~\cite[Test cases 5 and 6]{carmona2019convergence2} for more details. However, this approach requires estimating aggregate quantities, for example by simulating a finite but large population of particles for many realizations of the common noise. When the interactions are through moments, another approach is to use only one realization of the idiosyncratic noise for each realization of the common noise: 
In \cite{MinHu:21}, based on the rough path theory, a single-loop algorithm called signatured deep fictitious play has been proposed. The proposed algorithm can accurately capture the effect of common uncertainty changes on mean-field equilibria without further training of neural networks, as previously needed in the existing machine learning algorithms. We will provide more details on this method in Section~\ref{sec4_MFG_with_CN} below. Figure~\ref{fig:LQ} showcases the performance for this LQ MFG with common noise, where the benchmark trajectories are simulated according to \eqref{def:LQ_SDE}  with $m_t$ and $\alpha_t$ in \eqref{def:LQ_m} and \eqref{def:LQ_alpha}. 

\begin{figure}[!htb]
    \centering
    \subfloat[$X_t$]{
         \includegraphics[width=0.3\columnwidth]{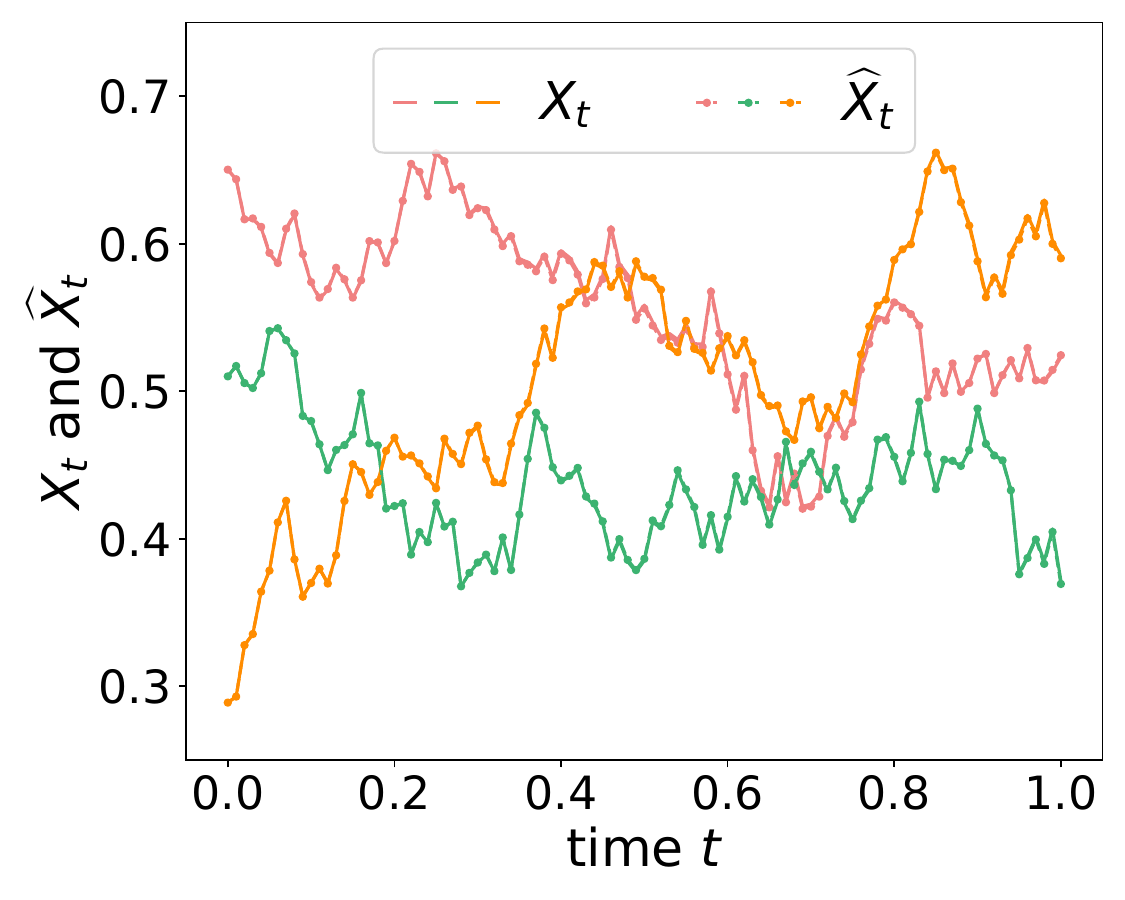}}
    \subfloat[$m_t = \EE(X_t \vert \mcF_t^{W^0})$]{
         \includegraphics[width=0.3\columnwidth]{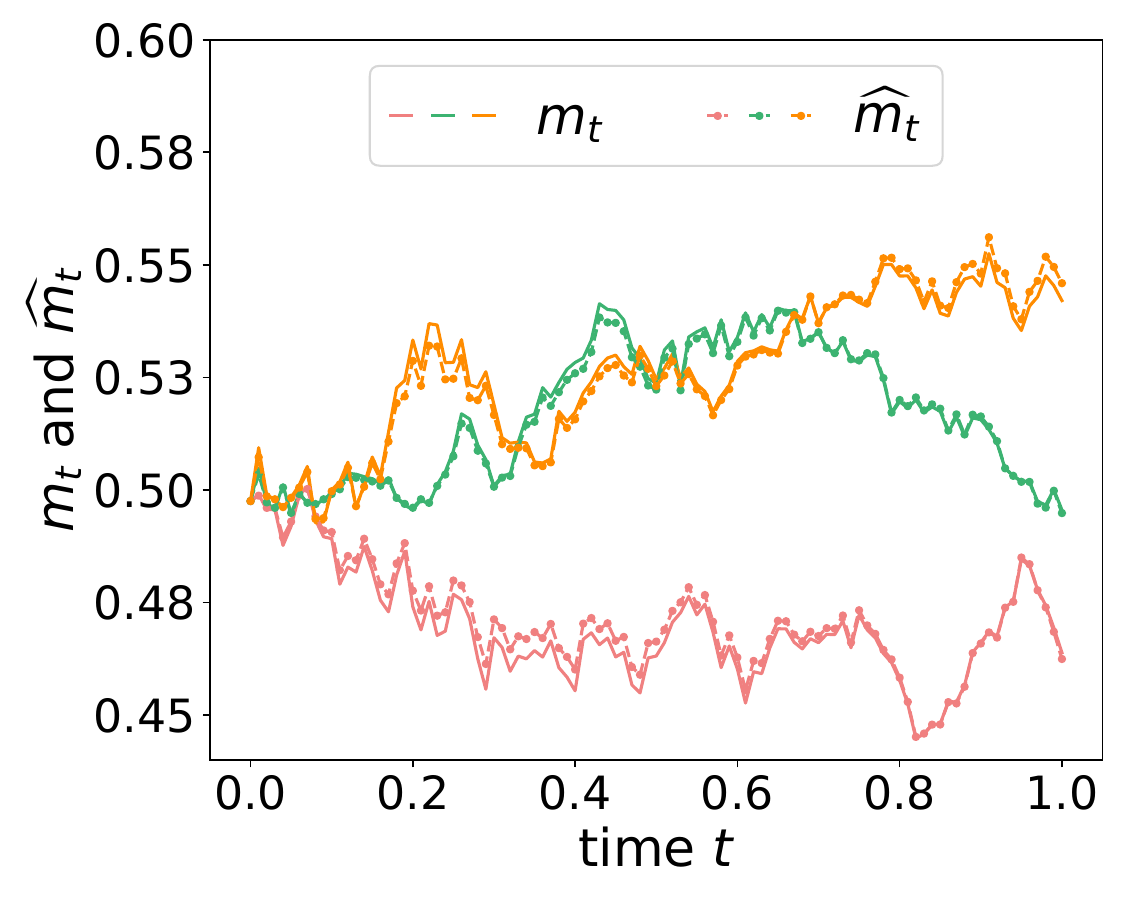}}
    \subfloat[Minimized Cost]{
         \includegraphics[width=0.3\columnwidth]{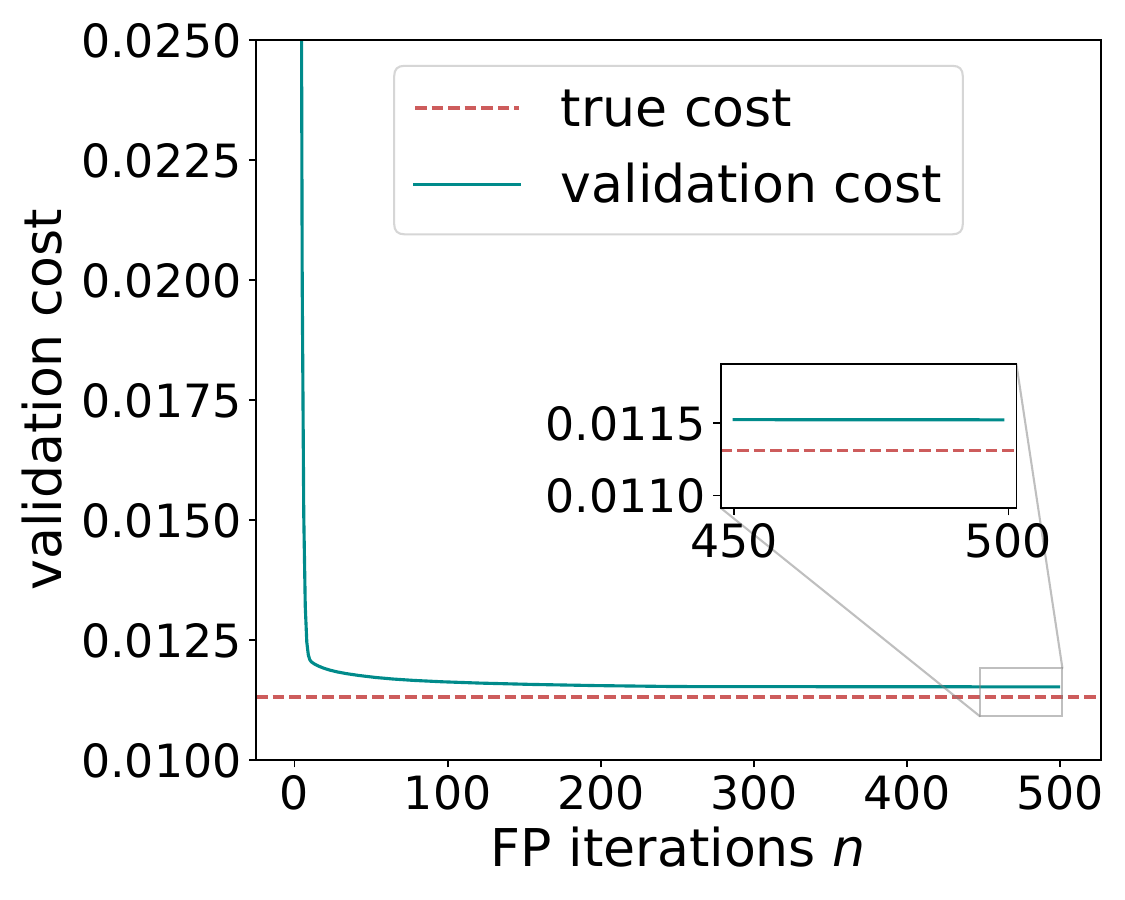}}
    \caption{The illustrative linear quadratic model in Section~\ref{sec:intro-LQsysrisk}. Panels (a) and (b) give three trajectories of $X_t$, $m_t = \EE[X_t \vert \mcF_t^{W^0}]$ (solid lines) and their approximations $\widehat X_t$ (dashed lines) using different realizations of $(X_0, W, W^0)$ from validation data. Panel (c) shows the minimized cost computed using validation data over fictitious play iterations. 
    Parameter choices are given in \cite[Section 5]{MinHu:21}. %
    }
    \label{fig:LQ}
\end{figure}

\subsection{Organization of the survey}
In the rest of this survey, we review recent developments in machine learning methods and theory for stochastic control and differential games. We shall also identify unsolved challenges, and make connections to real applications. The topics are organized based on the number of players involved in the problem, first in the model-based setting, and then in the model-free setting.
We first review deep learning algorithms for stochastic control problems in Section~\ref{sec:SCP}, including mean field control problems, viewed as a type of control problem where the individual state’s dynamics are influenced by its own law. In Section~\ref{sec:SDG}, we focus on deep learning for stochastic differential games, including (moderately large) $N$-player games and mean-field games. In Section~\ref{sec:RL}, we review the basic principles underpinning model-free reinforcement learning methods for stochastic control and games. We make conclusive remarks and discuss unsolved challenges in Section~\ref{sec:conclusion}. In the appendix, we provide more background on two main tools of modern machine learning, namely, neural network architectures and stochastic gradient descent. We also summarize all the acronyms and the notations frequently used in this paper.

\medskip
\noindent{\bf The scope of the survey. }
This paper focuses on the recently developed neural network-based algorithms aiming at high-dimensional problems. The main reason why we focus on the stochastic setting is that in the deterministic setting, the problems can usually be tackled using open loop controls by reducing the problem to a two-point boundary value problem \cite{kierzenka2001bvp,zang2022machine}, which can be solved in high dimension without machine learning. 

Besides the methods reviewed in this paper, there is a rich literature on methods, some of them related to machine learning but without neural networks. For example, to cite just a few examples, Markov chain based methods \cite{budhiraja2007convergent}, regression based methods \cite{barrera2006numerical,bouchard2004discrete}, and approaches based on PDEs and BSDEs \cite{chen2008semi,forsyth2007numerical}; see, {\it e.g.}, %
the survey papers \cite{kushner1990numerical,jin2022survey} and the references therein. Furthermore, we focus mostly on standard classes of stochastic control problems and games but many other problems are considered in the literature. For instance, we do not discuss in this paper optimal switching and optimal stopping problems, for which numerical algorithms have been extensively developed, \emph{e.g.}, in \cite{kohler2010pricing,HuLudkovski17,becker2019deep,becker2020pricing,hu2020deep,hurephamwarin2020deep,lapeyre2021neural,reppen2022deep,reppen2022neural,gao2022convergence,bayraktar2022deep}.

\section{Stochastic Control Problems}\label{sec:SCP}

We start by reviewing the development of machine learning methods for stochastic optimal control problems, also called simply stochastic control (SC) problems. SC is a long-standing topic that studies stochastic dynamic systems that are controlled so as to achieve optimal performance. It has  applications in many areas, including but not limited to engineering, economics, and mathematical finance. Some existing methods have been reviewed in, {\it e.g.}, works by Pham and co-authors \cite{pham2005some,germain2021neuralfinance}. In this section, we shall focus on reviewing methods based on deep learning for solving standard SC problems, or problems with delay and mean-field type SC problems, BSDE-based deep learning algorithms, primal-dual approaches, and PDE-based algorithms. We focus here on the algorithms and we refer to  Appendix~\ref{sec:DLtools} for a basic introduction to neural networks and deep learning.

\subsection{Formulation of stochastic control}\label{sec:control_intro}

Let $T<\infty$ be a time horizon. Let $(\Omega, \MCF, \PP)$ be a probability space supporting an $m$-dimensional Brownian motion $W = (W_t)_{t \in [0,T]}$, and $\mathbb{F} = (\mc{F}_t)_{t\in[0,T]}$ be the natural filtration generated by $W$. In the most common case, a stochastic control problem is formulated as follows. Let $d$ and $k$ be integers for the dimensions of the state and the action. Let $\mc{A} \subset \RR^k$ denote the set of admissible actions. Let $b, \sigma, f, g$ be Borel-measurable functions,
\begin{equation}
\label{def:oc-bsigmaf}
    (b, \sigma, f) : [0,T] \times \RR^d \times \mc{A} \to (\RR^d, \RR^{d \times m}, \RR), \quad g : \RR^d \to \RR.
\end{equation}
We denote by $\mathbb{A}$ the set of so-called admissible controls. This set describes the integrability and measurability conditions required on $\alpha$. We usually require that $\alpha$ is square-integrable. For measurability, two popular choices are that either $\alpha$ should be a $\mc{F}_t$-progressively measurable process, or $\alpha$ should be expressed as a measurable function of $(t, X_t)$. The former case is called \emph{adaptive controls} (also called \emph{open-loop}) while the latter one is referred to as \emph{Markovian controls} or \emph{closed-loop controls in feedback form} (we use the two terms interchangeably in the sequel).  We will sometimes consider \emph{closed-loop controls}, which means controls that are adapted to the filtration generated by $X$. 
We shall see that different choices of $\mathbb{A}$ might lead to different algorithm designs in Section~\ref{sec:control-direct}. 

\begin{defi}[Stochastic control problem]
An agent controls her state process $X$ through an action process $\alpha$ taking values respectively in $\RR^d$ and $\mc{A}$, where the dynamics of $X$ are given by the stochastic differential equation (SDE), 
\begin{equation}\label{def:control-Xt}
    \ud X_t = b(t, X_t, \alpha_t) \ud t + \sigma(t, X_t, \alpha_t) \ud W_t, \quad X_0 = x_0.
\end{equation}
The agent aims to minimize the expected cost
\begin{equation}\label{def:control-cost}
    J(\alpha): \ctrl \mapsto \EE \left[\int_0^T f(t, X_t, \alpha_t) \ud t + g(X_T) \right],
\end{equation}
over the set of admissible action processes, denoted by $\mathbb{A}$ and to be discussed below. 
\end{defi}

We briefly describe how the above problem can be tackled by PDEs, or (F)BSDEs.

\noindent
\textbf{PDE approach.} When considering Markovian controls, one can define the value function $u:[0,T] \times \RR^d \to \RR$,
\begin{equation}
    u(t, x) = \inf_\alpha \EE\left[\int_t^T f(s, X_s, \alpha_s) \ud s + g(X_T) \vert X_t = x\right], 
\end{equation}
and employ the dynamic programming principle (DPP) \cite[Section~3]{PH:09}: $u(T, x) = g(x)$ and for any stopping time $\tau \in [t, T]$, 
\begin{equation}
    u(t,x) = \inf_\alpha \EE\left[\int_t^\tau f(s, X_s, \alpha_s) \ud s + u(\tau, X_\tau) \vert X_t = x\right].
\end{equation}
Then, one can derive the Hamilton-Jacobi-Bellman (HJB) equation, which describes the evolution of the value function. Under suitable conditions, $u$ solves
\begin{equation}\label{def:control-HJB}
    \begin{cases}
   	\partial_t u (t,x)  + \min_{\alpha \in \mc{A}} H(t, x, \grad_x u(t,x), \Hess_x u(t, x), \alpha) = 0, \qquad (t,x) \in [0,T) \times \RR^d 
   	\\
   	u(T,x) = g(x), \qquad x \in \RR^d,
   	\end{cases}
\end{equation}
where the Hamiltonian $H$ is defined as
\begin{equation}\label{def:control-H}
    H(t, x, p, q, \alpha) = b(t, x, \alpha) \cdot p + \half \Tr(\sigma(t, x, \alpha)\sigma(t,x, \alpha)\transpose q) + f(t, x, \alpha),
\end{equation} 
and $\Hess_x u(t, x)$ as the Hessian matrix of $u$ with respect to $x$. 
If \eqref{def:control-HJB} has a classical solution, then the optimal control is given by
\begin{equation*}
	\hat{\ctrl}(t,x) = \hat{\ctrl}(t, x, \grad_x u(t,x), \Hess_x u(t,x)) = 
	\argmin_{\alpha \in \mc{A}} H(t, x, \grad_x u(t, x), \Hess_x u(t, x), \alpha).
\end{equation*}

\noindent
\textbf{BSDE approach.} The connection between SC and BSDEs can be established in two different ways: by representing the value function or its derivative as the solution of a BSDE. When the volatility is uncontrolled,  that is, $\sigma(t, x, \alpha)$ is free of $\alpha$, then $\hat \ctrl$ does not depend on $\Hess_x u(t,x)$ and the PDE \eqref{def:control-HJB} becomes semi-linear
\begin{multline}
	\partial_t u (t,x)  + \half \Tr(\sigma(t, x)\sigma(t,x)\transpose \Hess_x u(t, x)) + b(t, x, 	\hat{\ctrl}(t, x, \grad_x u(t,x))) \cdot  \grad_x u(t,x) \\
	+ f(t, x, \hat{\ctrl}(t, x, \grad_x u(t,x))) = 0.
\end{multline}
In this case, suppose that there exist functions $\tilde b(t, x)$ and $h(t, x, z)$ such that 
$$
    \tilde b(t, x) \cdot \grad_x u(t,x) + h(t, x, \sigma(t,x) \transpose \grad_x u(t,x)) = b(t, x, \hat{\ctrl}(t, x, \grad_x u(t,x))) \cdot  \grad_x u(t,x) + f(t, x, \hat{\ctrl}(t, x, \grad_x u(t,x))).
$$
Then the nonlinear Feynman-Kac formula \cite{PaPe:90} gives the following BSDE interpretation of $u$,
\begin{equation}\label{def_control_BSDEreform}
    \begin{dcases}
    \ud \mc{X}_t = \tilde b(t, \mc{X}_t) \ud t + \sigma(t, \mc{X}_t) \ud W_t, \quad \mc{X}_0 \sim \mu_0, \\
    \ud \mc{Y}_t = -h(t, \mc{X}_t, \mc{Z}_t)\ud t + \mc{Z}_t \ud W_t, \quad \mc{Y}_T = g(\mc{X}_T),
    \end{dcases}
\end{equation}
through the relations
$$
    \mc{Y}_t = u(t, \mc{X}_t), \quad \mc{Z}_t = \sigma(t, \mc{X}_t) \transpose \grad_x u(t, \mc{X}_t).
$$
If $\mu_0 = \delta_{x_0}$ is concentrated on a single initial state $x_0$, then the optimal value is given by the value of the BSDE solution at time $0$, {\it i.e.}, $\inf_{\alpha} J(\alpha) = \mc{Y}_0$. If $\tilde b(t, x)$ is chosen to be identically zero, this equality can also be obtained by the comparison principle of BSDEs \cite[Proposition 4.1]{Ca:16}.

In the controlled volatility case, the PDE \eqref{def:control-HJB} is fully nonlinear, and its solution is connected to a solution of the second order BSDE (2BSDE) \cite{cheridito2007second}. If one chooses $\tilde b(t,x)$ and $\Sigma(t,x)$ such that $h$ is determined by
\begin{equation}
    H(t, x, p, q, \hat\alpha(t, x, p, q)) = \tilde b(t, x) \cdot p + h(t, x, p, q) + \half \Tr(\Sigma(t, x)\Sigma(t, x)\transpose q),
\end{equation}
then the solution to the 2BSDE
\begin{equation}\label{def_control_2BSDEreform}
    \begin{dcases}
    \ud \mc{X}_t = \tilde b(t, \mc{X}_t) \ud t + \Sigma(t, \mc{X}_t) \ud W_t, \quad \mc{X}_0 = x_0, \\
    \ud \mc{Y}_t = -h(t, \mc{X}_t, \mc{Y}_t, \mc{Z}_t)\ud t + \mc{Z}_t\transpose \Sigma(t, \mc{X}_t)\ud W_t, \quad Y_T =  g(\mc{X}_T),\\
    \ud \mc{Z}_t = \mc{A}_t \ud t + \Gamma_t \Sigma(t, \mc{X}_t)\ud W_t,  \quad \mc{Z}_T = \grad_x g(\mc{X}_T),
    \end{dcases}
\end{equation}
gives an interpretation of the solution to the PDE \eqref{def:control-HJB} through the relations
$$
    \mc{Y}_t = u(t, \mc{X}_t), \quad \mc{Z}_t =  \grad_x u(t, \mc{X}_t), \quad \Gamma_t = \Hess_x u(t, \mc{X}_t), \quad \mc{A}_t = \mathfrak{L}\grad_x u(t, \mc{X}_t),
$$
where $\mathfrak{L}$ denotes the infinitesimal generator of $\mc{X}$.

\noindent
\textbf{FBSDE approach.} 
The Pontryagin stochastic maximum principle provides the connection to the FBSDE, where the forward and backward equations are coupled. Define the generalized Hamiltonian $\mathcal{H}$ by
\begin{equation}
\label{eq:Hamiltonian-bsde-control}
\mc{H}(t,x,y,z,\alpha) = b(t, x, \alpha) y + \Tr(\sigma(t, x, \alpha)z) + f(t, x, \alpha).
\end{equation}
If the Hamiltonian $\mc{H}$ is convex in $(x, \alpha)$, 
and $(X_t, Y_t, Z_t)$ solves
\begin{equation}\label{def_control_FBSDEreform}
    \begin{dcases}
    \ud {X}_t = b(t, {X}_t, \hat \alpha_t) \ud t +\sigma(t, {X}_t, \hat \alpha_t) \ud W_t, \qquad {X}_0 = x_0,\\
    \ud {Y}_t = -\grad_x \mc{H}(t, {X}_t, {Y}_t,  {Z}_t,  \hat \alpha_t ) \ud t + {Z}_t \ud W_t, \quad Y_T = \partial_x g({X}_T),\\
    \hat{\alpha}_t = \inf_{\alpha} \mc{H}(t, {X}_t, {Y}_t,  {Z}_t, \alpha ),
    \end{dcases}
\end{equation}
then $\hat \alpha$ is the optimal control.  If the value function is smooth enough, then 
\begin{equation}
     Y_t = \grad_x u(t,  X_t), \quad Z_t = \sigma(t, X_t, \hat \alpha_t)\transpose\Hess_x u(t,  X_t).
\end{equation}

\noindent
\textbf{Classical numerical methods. }
To solve stochastic optimal control problems, one typically tries to compute either the solution to the HJB equation~\eqref{def:control-HJB} or the solution to the BSDE in~\eqref{def_control_BSDEreform}. To solve HJB PDEs, classical methods include finite difference schemes~\cite{bonnans2003consistency,oberman2006convergent,jensen2013convergence}, semi-Lagrangian schemes~\cite{debrabant2013semi} and finite element schemes~\cite{beard1997galerkin,boulbrachene2001finite}. We refer to~\cite{feng2013recent} for a review of numerical methods for non-linear second order PDEs (including HJB equations) and to \cite{falcone2016numerical} for a review of finite-difference and semi-Lagrangian schemes for Hamilton-Jacobi type equations. 
For a review of numerical methods for HJB equations arising in finance, we refer to~\cite{forsyth2007numerical}. To solve BSDEs, several numerical methods have also been proposed. Most approaches approximate the backward variable $Y$ using an Euler scheme along a time grid by representing it as a conditional expectation~\cite{bouchard2004discrete,zhang2004numerical}, which can be computed in various ways. Common methods include least-square regression using a set of basis functions~\cite{gobet2005regression,bender2012least}, quantization~\cite{bally2001stochastic,bally2003error} and tree-based approximation~\cite{chevance1997numerical,ma2002numerical}. 
Furthermore, several refinements have been studied, such as higher-order schemes~\cite{zhao2006new,zhao2010stable,chassagneux2014linear}.  We refer to~\cite{chessari2023numerical} for a recent review of numerical methods for BSDEs. We stress that the two approaches (based on PDEs and BSDEs) are connected, see {\it e.g.}~\cite{bouchard2009discrete}. Last, there are other methods, such as the Markov chain approximation method~\cite{kushner2013numerical}. 
Some of these methods have been extended to MFC problems, for which deep learning methods will be discussed in the sequel. For instance~\cite{achdoulauriere2015systemmfc} used a finite-difference scheme for the PDE system, \cite{MR3575615} proposed an augmented Lagrangian method, and \cite{balata2019class} adapted several numerical methods for BSDEs to finite-dimensional MFC problems. 

These classical numerical methods are well understood thanks to the rich background of numerical analysis and most of them have been extensively studied. In many cases, rigorous proofs of convergence have been obtained, for both the schemes ({\it e.g.}, time and space discretizations) and the algorithms. Rates of convergence have also been obtained under suitable assumptions but, to the best of our knowledge, in general these rates suffer from the curse of dimensionality. In practice, Monte Carlo-based methods tend to scale up better than PDE-based methods but remain limited in terms of dimensionality except for specific classes of functions, see {\it e.g.}~\cite{sloan1998quasi} for the provable efficiency of quasi-Monte Carlo methods to compute high-dimensional integrals in particular situations.

For machine learning algorithms introduced in the following sections, if a temporal discretization is needed, we shall consider, for simplicity, a uniform grid $\pi$ on the interval $[0,T]$, {\it i.e.}, a partition $0 = t_0 < t_1 \dots < t_{N_T} = T$, with $t_n - t_{n-1} = \Delta t = T / N_T$.

\subsection{Direct parameterization}\label{sec:control-direct}
We refer to the first class of algorithms as direct parameterization methods, which directly replace the control function with a neural network with appropriate inputs.  Such ideas can be traced back to earlier works such as {\it e.g.}~\cite{psaltis1988multilayered,hunt1992neuralnetworkcontrolsurvey,lehalle1998piecewise,MR2137498}, in which neural networks are used for control and optimal control problems in relatively low dimension and with shallow networks. 

\subsubsection{Global in time approach} 
\label{sec:direct-global-local}

There are two types of direct parameterization methods, depending on how the neural networks get trained.

We start with the method which trains the neural networks used for the control using the whole time horizon at once. More recently, Han and E \cite{han2016deep-googlecitations} were the first to generalize this type of methods to solve problems in high dimensions using recent modern machine learning techniques such as deep neural networks and efficient built-in stochastic gradient descent solvers. This has motivated fruitful studies on high-dimensional control problems. In particular, the algorithms reviewed below in Sections~\ref{sec:sc_delay} and \ref{sec:directMethod} for SC problems with delay and mean field control are both in the spirit of \cite{han2016deep-googlecitations}. More precisely, \cite{han2016deep-googlecitations} solved the following discrete time version of \eqref{def:control-Xt}--\eqref{def:control-cost},
\begin{align}
     &\checkX_{t_{n+1}} = \checkX_{t_n} + b(t_n, \checkX_{t_n}, \alpha_{t_n}) \Delta t + \sigma(t_n, \checkX_{t_n}, \alpha_{t_n}) \Delta \check W_{t_n}, \label{eq:control-Xt-discrete} \\
    &\min_{(\alpha_{t_n})_{n}} \EE\left[\sum_{n=0}^{N_T-1} f(t_n, \checkX_{t_n}, \alpha_{t_n}) \Delta t + g(\checkX_T)\right] \label{eq:control-Xt-cost},
\end{align}
where $\Delta \check W_{t_n} = \check W_{t_{n+1}} - \check W_{t_n}$ are i.i.d random variables with distribution $\mathcal{N}(0, \Delta t)$. Focusing on Markovian controls, they approximate at each time step the control $\alpha_{t_n}$ in \eqref{eq:control-Xt-discrete} by a feedforward neural network $\alpha_{t_n}(\cdot;\theta_n)$ taking inputs $\check X_{t_n}$, where $\theta_n$ denotes all NN's parameters at time $t_n$. We thus obtain the following cost, which is interpreted as a loss that can be minimized by SGD:
\begin{equation}\label{def:control-loss}
    \check J(\theta) = \EE\left[\sum_{n=0}^{N_T-1} f(t_n, \checkX_{t_n}^\theta, \alpha_{t_n}(\check X_{t_n}^\theta; \theta_n)) \Delta t + g(\checkX_T^\theta)\right], 
    \quad \theta = (\theta_n)_{n=0}^{N_T},
\end{equation}
where $(\checkX_{t_n}^\theta)_{n=1}^{N_T}$ follows
\begin{equation}\label{def:control-Xt-nn}
 \checkX_{t_{n+1}}^\theta = \checkX_{t_n}^\theta + b(t_n, \checkX_{t_n}^\theta, \alpha_{t_n}(\check X_{t_n}^\theta; \theta_n)) \Delta t + \sigma(t_n, \checkX_{t_n}^\theta, \alpha_{t_n}(\check X_{t_n}^\theta; \theta_n)) \Delta \check W_{t_n}.
 \end{equation}
In practice, the expected value in \eqref{def:control-loss} is approximated by the following quantity based on Monte Carlo simulations:
\begin{equation}
    L(\theta,S) = \frac{1}{N}\sum_{j=1}^N\left[\sum_{n=0}^{N_T-1} f(t_n, \checkX_{t_n}^{j,\theta}, \alpha_{t_n}(\checkX_{t_n}^{j,\theta}; \theta_n)) \Delta t + g(\checkX_T^{j,\theta})\right],
\end{equation}
where $\{(\check X_{t_n}^{j, \theta})_{n=1}^{N_T}\}_{j = 1, \ldots, N}$ are sample paths of $(\checkX_{t_n}^{\theta})_{n=1}^{N_T}$ in \eqref{def:control-Xt-nn}  using i.i.d. samples of  $(\Delta \check W_{t_n})_{n=1}^{N_T}$. 

Constraints on the states and the controls can be taken into account by adding penalty terms in the loss function as in {\it e.g.} \cite{han2016deep-googlecitations}. 

\begin{remark}\label{rem:global-single-nn}
  Alternatively, one can use a single neural network $(t,x) \mapsto \alpha(t,x;\theta)$ and evaluate it at $(t,x) = (t_n,\check X_{t_n})$, $n = 0, \ldots N_T-1$. In this case, the network can directly capture the time continuity in the control process (if it is continuous). Furthermore, it can be trained for various grids of points in time, and after training, it can be used for arbitrary time points. 
\end{remark}

\subsubsection{Local in time approach}
The second approach has been proposed by Bachouch, Hur\'{e},  Langren\'{e} and Pham in \cite{bachouch2021deepnumerical}. It combines classical dynamic programming (DP) and deep neural networks for approximating the control and possibly the value function. There are two main versions. The first version, called NNContPI, is designed as follows. Assuming that the optimal controls at time steps $t_{n+1}, \ldots, t_{N_T-1}$ are already learnt with neural network parameters $\hat \theta_{n+1}, \ldots, \hat \theta_{N_T-1}$,
the optimal control at time $t_n$ is approximated by a neural network $\alpha_{t_n}(\cdot; \hat \theta_n)$ where the optimized parameters are determined by
\begin{equation}\label{eq:NNContPI}
    \hat \theta_n \in \argmin_{\theta} \EE\Big[f(t_n, \checkX_{t_n}^\theta, \alpha_{t_n}(\checkX_{t_n}^\theta; \theta))\Delta t + \sum_{n'=n+1}^{N_T-1}f(t_{n'}, \checkX_{t_{n'}}^\theta, \alpha_{t_{n'}}(\checkX_{t_{n'}}^\theta; \hat \theta_{n'}))\Delta t + g(\checkX_T^\theta)\Big],
\end{equation}
where $(\checkX_{t_{n'}}^\theta)_{n'=n+1}^{N_T}$ follows \eqref{def:control-Xt-nn} with  $\alpha_{t_n}(\cdot; \theta), \alpha_{t_{n+1}}(\cdot; \hat \theta_{n+1}), \ldots, \alpha_{t_{N_T-1}}(\cdot; \hat \theta_{N_T-1})$ being used. Then, as always, the expected value is approximated by sample paths of $\check X^{j, \theta}$, and the optimal $\hat \theta_n$ is obtained by SGD as described in Appendix~\ref{sec:DLtools}. The second version, termed as Hybrid-now, further approximates the value function ({\it i.e.}, the cost-to-go) at time $t_{n+1}$ using a deep neural network to avoid repeated computation of the last two terms in \eqref{eq:NNContPI}.
More precisely, given the learned approximated value function $ V_{t_{n+1}}(\cdot; \tilde\theta_{n+1})$ at time $t_{n+1}$, the optimal policy at time $t_n$ is determined in a manner that is similar to the NNContPI algorithm
\begin{equation}\label{eq:Hybrid-Now}
    \hat\theta_n \in \argmin_{\theta} \EE[f(t_n, \checkX_{t_n}^\theta, \alpha_{t_n}(\checkX_{t_n}^\theta; \theta))\Delta t +  V_{t_{n+1}}(\checkX_{t_{n+1}}^\theta; \tilde \theta_{n+1})],
\end{equation}
and the value function $V_{t_n}$ at time $t_n$ is then approximated by another neural network $V_{t_n}(\cdot; \tilde \theta_n)$ whose parameters are determined by
\begin{equation}
     \tilde \theta_n \in \argmin_\theta \EE\Big[ |f(t_n, \checkX_{t_n}^\theta, \alpha_{t_n}(\checkX_{t_n}^\theta; \hat\theta_n)) \Delta t +  V_{t_{n+1}}(\checkX_{t_{n+1}}^\theta; \tilde \theta_{n+1}) - V_{t_n}(\check X_{t_n}^\theta; \theta) |^2 \Big].
\end{equation}

\begin{remark}
 The first approach discussed above (global in time, developed by \cite{han2016deep-googlecitations}) learns all the optimal controls $\alpha_{t_n}(\cdot; \hat\theta_n)$, $n = 0, \ldots, N_T-1$ at once, by performing a unique SGD with a loss function that involves the whole time horizon, while the second approach \cite{bachouch2021deepnumerical} learns $\alpha_{t_n}$ sequentially and backwardly, for $n = N_T-1, N_T-2, \ldots, 0$. The first approach may be more efficient memory-wise since it is possible to use a single neural network, as noticed in Remark~\ref{rem:global-single-nn}. In this way, at the end of training, we obtain a function which can be applied at point in time and space. Even if the neural network is trained only on a discrete grid of time steps, we can expect the neural network to generalize relatively well at intermediate times. This is not possible with the second approach, which requires training a different neural network at each time step. On the other hand, the first approach, since it optimizes over the whole time horizon at once, may struggle to capture steep changes in time. Since the gradients are accumulated over all the time steps and the errors are compounded, it may encounter vanishing or exploding gradient problem for large $N_T$ as remarked in \cite{bachouch2021deepnumerical}. 
\end{remark}

\begin{remark}
It is common to use the same architecture for the neural networks $\alpha_{t_n}$ of all the time steps. Even with the same architecture, the parameters can be different, which yields different functions for different time steps.  
\end{remark}

\begin{rem}[Theoretical analysis]
In \cite[Theorems 4.6 and 4.13]{hure2021deepconvergence}, the authors provided the consistency and rates of convergence for the control and the value function subject to universal approximation error of the neural networks.
\end{rem}

\subsection{BSDE-based deep learning algorithms}\label{sec:BSDE}
We now turn our attention to alternative deep learning methods for control problems \eqref{def:control-Xt}--\eqref{def:control-cost}, which are based on solving the associated backward stochastic differential equations (BSDEs) (\emph{cf.} \eqref{def_control_BSDEreform} and \eqref{def_control_2BSDEreform}). To simplify the presentation, we explain the algorithms on a generic BSDE,
\begin{equation}\label{def:BSDE-Yt}
      \ud Y_t = - F(t, X_t, Y_t, Z_t) \ud t + Z_t\transpose \ud W_t, \quad t \in [0,T], \qquad Y_T = G(X_T),
\end{equation}
where $X$ solves the (possibly) coupled forward equation
\begin{equation}\label{def:BSDE-Xt}
   \ud X_t = B(t, X_t, Y_t, Z_t) \ud t + \sigma(t, X_t) \ud W_t, \quad t \in [0,T], \qquad X_0 \sim \mu_0.
\end{equation}

As mentioned in Section~\ref{sec:control_intro}, the connection between optimal control and BSDEs (or FBSDEs) can be established in several ways.   

\subsubsection{Deep backward stochastic differential equation (Deep BSDE) method}
\label{subsubsec:deepBSDE}
The Deep BSDE method was proposed by E, Han and Jentzen in~\cite{MR3736669,HaJeE:18}, which, to the best of our knowledge, were the first works to use deep learning to solve BSDEs in high dimension and has inspired many followup works. This method relies on a stochastic version of a shooting method that is extensively used to solve ODEs. The strategy has been successfully applied to problems in economic contract theory where it is sometimes referred to as Sannikov's trick; see for example \cite{MR1766421,MR2963805,MR3738664}. The idea is to try to guess the initial value $Y_0$ and the $(Z_t)_{t \in [0,T]}$ process so as to meet the terminal condition $Y_T = g(X_T)$. It has originally been proposed in~\cite{MR3736669} for decoupled FBSDE systems, and then extended to fully coupled FBSDEs; see {\it e.g.}~\cite{hanlong2020convergence, ji2020three}. 

To be precise, solving the system \eqref{def:BSDE-Yt}--\eqref{def:BSDE-Xt} is reduced to identifying $(Y_0, (Z_t)_{t \in [0,T]})$ in \eqref{def:BSDE-Yt} which is characterized as the solution to the following control problem, {\it i.e.}, to minimize the cost over $(y_0,z)$ 
$$
    J(y_0,z) = \EE \left[ \| Y_T - g(X_T)\|^2\right],
$$
subject to:
$$
\begin{cases}
     \ud X_t = B(t, X_t, Y_t, z(t,X_t)) \ud t + \sigma(t, X_t) \ud W_t, \qquad t \in [0,T],
     \\
     \ud Y_t = - F(t, X_t, Y_t, z(t,X_t)) \ud t + z(t,X_t)\transpose \ud W_t, \qquad t \in [0,T],
     \\
     X_0 \sim \mu_0, \qquad Y_0 = y_0(X_0).
\end{cases}
$$
Then, as in the methods presented above (see Section~\ref{sec:control-direct}), the control functions $y_0$ and $z$ are replaced by deep neural networks, say $y_0(\cdot; \theta^Y)$ and $z(\cdot; \theta^Z)$ with parameters $\theta^Y$ and $\theta^Z$, respectively. Furthermore, time is discretized using a uniform grid $t_0 < t_1 \dots < t_{N_T} = T$, $t_n - t_{n-1} = \Delta t = T / N_T$. Then the problem becomes to minimize over $\theta = (\theta^Y,\theta^Z)$ the cost 
$$
    \check{J}(\theta) = \EE \left[ \| \check{Y}^{\theta}_T - g(\check{X}^{\theta}_T)\|^2\right],
$$
subject to:
$$
\begin{cases}
     \check{X}^{\theta}_{t_{n+1}} - \check{X}^{\theta}_{t_{n}} = B(t_n, \check{X}^{\theta}_{t_{n}}, \check{Y}^{\theta}_{t_{n}}, z(t_{n},\check{X}^{\theta}_{t_{n}}); \theta^Z) \Delta t + \sigma(t_n, \check{X}^{\theta}_{t_{n}}) \Delta \check{W}_{t_{n+1}}, \qquad n=0,\dots,N_T-1,
     \\
     \check{Y}^{\theta}_{t_{n+1}} - \check{Y}^{\theta}_{t_{n}}  = - F(t_n, \check{X}^{\theta}_{t_{n}}, \check{Y}^{\theta}_{t_{n}}, z(t_{n},\check{X}^{\theta}_{t_{n}}); {\theta^Z}) \Delta t + z(t_{n},\check{X}^{\theta}_{t_{n}}; {\theta^Z})\transpose \Delta \check{W}_{t_{n+1}}, \qquad n=0,\dots,N_T-1,
     \\
     \check{X}^{\theta}_0 = X_0 \sim \mu_0, \qquad \check{Y}^{\theta}_0 = y_{0}(\check{X}^{\theta}_0; \theta^Y),
\end{cases}
$$
where 
$
    \Delta \check{W}_{t_{n+1}} = \check{W}_{t_{n+1}}-\check{W}_{t_{n}}.
$

Finally, the optimization can be carried out by applying the SGD algorithm (Algorithm~\ref{algo:SGD-generic} in Appendix~\ref{sec:SGD-var}) with $L(\theta) = \check{J}(\theta)$ and one sample is 
$
    \xi = (X_0, (\Delta \check{W}_{t_{n}})_{n=1,\dots,N_T}),
$
which is sufficient to simulate $(\check{X}^{\theta}_{t_{n}}, \check{Y}^{\theta}_{t_{n}})_{n=0,\dots,N_T}$.

\begin{remark}
    Instead of using a single neural network for $z$, viewed as a function of $t$ and $x$, another possibility is to put a different neural network at each time step, which is a function of $x$ only. With this approach, each neural network can use fewer parameters since there are fewer inputs. However, possible drawbacks are that: (a) the number of neural networks grows linearly with the number of time steps; and (b) the time dependence is not captured (at least not directly). 
    
    A possible shortcoming of this method is that the optimization is done globally in time: the loss function is computed only after simulating a whole trajectory, and only then can the parameters be updated. For problems with long time horizons or complicated terminal conditions, the method may have difficulty in converging, as has been pointed out in~\cite{hurephamwarin2020deep}. 
\end{remark}

\begin{remark}[Theoretical analysis]
    A posteriori error bounds for the Deep BSDE method have been proved by Han and Long in~\cite[Theorems 1 and 2]{hanlong2020convergence}, extending the results of~\cite{bendersteiner2013posteriori} for uncoupled FBSDEs. They have shown, under suitable conditions, that the numerical error can be bounded by the value of the loss, and that this loss can be made as small as desired provided the approximation capability of the neural network is sufficient. Such results have then been extended to the fully coupled mean-field setting by Reisinger, Stockinger, and Zhang in \cite[Theorems 3.2 and 4.3]{reisinger2020posteriori}.
\end{remark}

Several variants of the Deep BSDE method have been proposed. For example, \cite{ji2020three} considers learning $Y$ as a feedback function of $X$ or using Picard iterations to learn feedback controls based on $(X,Y,Z)$. The Deep BSDE method has been extended to include control problems with mean-field effects~\cite[Secion~4.2]{carmona2019convergence2} and delay~\cite[Section~3.2]{fouque2019deep}, and generalized to stochastic differential games \cite[Section~3.2]{HaHu:19}.

As mentioned in Section~\ref{sec:control_intro}, the solution to a BSDE is closely related to the solution of a semi-linear PDE, which could be derived from an uncontrolled volatility problem. In the case of a fully controlled volatility problem, a similar relation exists, and Beck, E, and Jentzen \cite{beck2019machine} propose the corresponding deep 2BSDE method. It is noteworthy that several works have refined and extended the Deep BSDE method. For instance, \cite{chanwainam2019machine} improved the performance of the Deep BSDE method with specific architectures and training methods, illustrating them on various examples of semi-linear PDEs, including an HJB equation.

\subsubsection{Deep backward dynamic programming (DBDP)}

The DBDP method has been proposed by Hur\'e, Pham and Warin in \cite{hurephamwarin2020deep}, based on ideas similar to the local approach discussed in Section~\ref{sec:direct-global-local}. The main idea is to learn $\check{Y}_{t_n}$ and $\check{Z}_{t_n}$ at each $t_n$ as functions of $\check{X}_{t_n}$ by backward induction in time. So the resolution of the BSDE is decomposed as a sequence of optimization problems that are solved backward in time. This is in contrast with the Deep BSDE method, which goes forward in time, just like the difference between the two algorithms introduced in Section~\ref{sec:control-direct}. 

In DBDP, for each $n$, $\check{Y}_{t_n}$ and $\check{Z}_{t_n}$ are replaced by neural networks, say $y_n(\cdot; {\theta^Y_n})$ and $z_n(\cdot; {\theta^Z_n})$, with possibly different parameters at each time step. Here, one first chooses a sequence of distributions, say $\mu_{t_n}$ from which $\check{X}_{t_n}$ can be sampled for each $t_n$, $n=0,\dots,N_T$. The algorithm proceeds with a backward induction. First, $\theta^Y_{N_T}$ is trained such that $y_{N_T}(\cdot; {\theta^Y_{N_T}}) \approx g(\cdot)$, for example by minimizing
$$
    \check J(\theta^Y) = \EE\left[ \|y_{N_T}(\check{X}_{T}; {\theta^Y}) - g(\check{X}_{T})\|^2 \right],
$$
where $\check{X}_{T} \sim \mu_{N_T}$. Note that the value of $\theta^Z_{N_T}$ is not relevant to the result of the method.
 Then,  the neural networks $y_n(\cdot; {\theta^Y_n})$ and $z_n(\cdot;{\theta^Z_n})$ are trained for $n=N_T-1, N_T-2, \dots, 0$. There are at least two different ways to train these neural networks: 
\begin{itemize}
	\item Version 1: $\theta_n = (\theta^Y_n, \theta^Z_{n})$ is trained to minimize over $\theta = (\theta^Y,\theta^Z)$,
	\begin{align*}
		\check J^1_n(\theta) 
		&= \EE\Big[ \Big\|y_{n+1}(\check{X}^{\theta}_{t_{n+1}}; {\theta^Y_{n+1}}) - y_n(\check{X}_{t_n}; {\theta^Y})   +F(t_n, \check{X}_{t_n}, y_n(\check{X}_{t_n}; {\theta^Y}), z_n(\check{X}_{t_n}; {\theta^Z})) \Delta t 
		\\
		&\qquad\qquad - z_n(\check{X}_{t_n}; {\theta^Z}) \cdot \Delta \check{W}_{t_{n+1}}\Big\|^2\Big],	
	\end{align*}
	with
	$$
        \check{X}^{\theta}_{t_{n+1}} = \check{X}_{t_n} + B(t_n, \check{X}_{t_{n}}, y_n(\check{X}_{t_n}; {\theta^Y}), z_n(\check{X}_{t_{n}}; {\theta^Z})) \Delta t + \sigma(t_n, \check{X}_{t_n})  \Delta \check{W}_{t_{n+1}}, \qquad \check{X}_{t_n} \sim \mu_{t_n}.
    $$
	\item Version 2: $\theta^Y_n$ is trained to minimize over $\theta = \theta^Y$,
	\begin{align*}
		\check J^2_n(\theta) 
		&= \EE\Big[ \Big\|y_{n+1}(\check{X}^{\theta}_{t_{n+1}}; {\theta^Y_{n+1}}) - y_n(\check{X}_{t_n}; {\theta^Y})  + F(t_n, \check{X}_{t_n}, y(\check{X}_{t_n}; {\theta^Y}), \sigma\transpose \mathcal{D}_x y_n(\check{X}_{t_n}; {\theta^Y})) \Delta t 
		\\
		&\qquad\qquad - \mathcal{D}_x y_n(\check{X}_{t_n}; {\theta^Y})\transpose \sigma \Delta \check{W}_{t_{n+1}}\Big\|^2\Big],	
	\end{align*}
	with
	$$
        \check{X}^{\theta}_{t_{n+1}} = \check{X}_{t_n} + B(t_n, \check{X}_{t_{n}}, y_n(\check{X}_{t_n}; {\theta^Y}), \sigma\transpose \mathcal{D}_x y_n(\check{X}_{t_n}; {\theta^Y})) \Delta t + \sigma(t_n, \check{X}_{t_n}) \Delta \check{W}_{t_{n+1}}, \qquad \check{X}_{t_n} \sim \mu_{t_n},
    $$
    where the derivative $\mathcal{D}_x y_n$ represents the numerical differentiation of the neural network $y_n$. In this version, the $Z$ component is directly approximated by the derivative of the neural network for $Y$, so we do not use any parameter $\theta^Z$.
\end{itemize}

Then at each time step $t_n$, the optimization can be carried out by applying the SGD algorithm (see Algorithm~\ref{algo:SGD-generic} in Appendix~\ref{sec:SGD-var}) to the loss $\check J_n^1(\theta)$ or $\check J_n^2(\theta)$, and one sample is
$
    \xi = (\check{X}_{t_n}, \Delta \check{W}_{t_{n+1}}).
$ 
The DBDP method has been used successfully to solve BSDEs associated with semi-linear PDEs and fully non-linear PDEs (\cite{hurephamwarin2020deep} and \cite{pham2021neural}). Several refinements have been studied ({\it e.g.}, higher-order scheme in \cite{chassagneux2022deep}).

\begin{remark} 
In \cite{hurephamwarin2020deep}, the authors also extended this idea to reflected BSDEs that arise in optimal stopping problems and American option pricing in finance \cite[Section~6.5]{PH:09}. However, as mentioned in the introduction, we do not discuss optimal stopping problems in this survey for the sake of brevity.

The difference between the two versions of DBDP is: Version~1 uses independent neural networks for the approximation of $Y_t$ and $Z_t$ in \eqref{def:BSDE-Yt}; while Version~2 only approximates $Y_t$ by neural networks and represents $Z_t$ through auto-differentiating $y_n(\cdot; \theta^Y_n)$. Though the former one has more modeling flexibility, it may also introduce some inconsistency, as in many scenarios, the backward process is a function of time and the forward process $Y_t = u(t, X_t)$ and the adjoint process is indeed the derivative of this function up to a scaling factor: $Z_t = \sigma\transpose(t, X_t) \nabla_x u(t, X_t)$. 

The main advantage of the DBDP method is that it makes use of the time structure of the problem to split it into much simpler problems. Two possible shortcomings of the DBDP method are that: (1) the number of neural networks grows linearly with the number of time steps; and (2) it is not always clear how to choose the sampling distributions $\mu_{t_n}$, which have an impact on the way the neural networks are trained. 
\end{remark}

\begin{remark}[Theoretical analysis]
    Hur\'e, Pham and Warin analyzed in~\cite[Theorems 4.1 and 4.2]{hurephamwarin2020deep} the error and showed, assuming the gradient descent method used to solve the local problems converges to the true optimum, that the approximation error goes to zero as one increases the number of time steps and the number of neurons.
\end{remark}

A similar approach based on dynamic programming and backward induction can be applied to MFC (see Section~\ref{sec:directMethod}), but in this case, the value function is a function of the probability distribution, see e.g.~\cite{MR3343705,MR3258261,MR3631380,carmona2018probabilistic}. Since it is impossible to represent numerically all probability distributions, an approximation must be made. For instance, ~\cite{germain2021deepsets} used empirical distributions combined with symmetric neural networks, and \cite{hanhulong:22} used neural networks to represent directly the function maps evaluated at the optimal state's distribution.  

\subsection{Primal-Dual approaches}

The Deep BSDE method presented above tackles directly a BSDE, without making explicit use of the fact that, in our context, the BSDE comes from an optimal control problem. In this context and under suitable assumptions, we can rely on a dual formulation to introduce primal-dual deep learning methods. Such methods have been studied, {\it e.g.}, in~\cite{henrylabordere2017deep,benderschweizerzhuo2017primal,daveyzheng2020deep} and have applications to several problems in finance. 

Recall that the stochastic optimal control we consider is given by~\eqref{def:control-Xt}--\eqref{def:control-cost}. 
Using the aforementioned methods ({\it e.g.}, Deep BSDE or DBDP), it is generally hard to know how close to being optimal the neural network solution is. This is because we do not know \emph{a priori} the minimal cost. However, we are always sure that these methods provide an upper bound since, given any admissible control $\tilde\alpha$ ({\it e.g.}, in the form of a neural network as in the Deep BSDE method), we can compute $J(\tilde\alpha)$ which is at least as large as the infimum $J^* = \inf_{\alpha} J(\alpha)$. So to claim that $\tilde\alpha$ is almost optimal, it is enough to exhibit a lower bound on $J^*$ that is close to $J(\tilde\alpha)$.

Except in special cases, there is no analytical expression for a good lower bound, and traditional numerical methods might be inefficient if the problem is in high dimension. Fortunately, the optimal value can be computed through a dual problem which is formulated as a maximization problem and hence yields a lower bound.  We refer to {\it e.g.}~\cite{henrylabordere2016dual} for more details.

Formally, the dual problem can be expressed as
$$
    \sup_{\varphi} \EE\left[ \inf_{\alpha} \Phi^{\varphi,\alpha} \right],
$$
where
$$
    \Phi^{\varphi,\alpha} = g(X^\alpha_T) +  \int_0^T f(t, X^\alpha_t, \alpha_t) dt - \int_0^T \varphi(t, X^\alpha_t) \sigma(t, X_t^\alpha) \ud W_t,
$$
and the superscript $\alpha$ in $X_t^\alpha$ emphasizes the dependence on the control $\alpha$. Here the infimum over $\alpha$ is a pathwise optimization since it is inside the expectation. In some cases, it can be solved explicitly. Then we are left with a more standard optimal control problem with feedback control $\varphi$. The latter can be solved for instance using ideas similar to the ones presented in the previous sections such as the direct parameterization or the Deep BSDE method, see Sections~\ref{sec:control-direct} and~\ref{sec:BSDE}. We refer to {\it e.g.}~\cite{henrylabordere2017deep} for more details and numerical examples.

\subsection{PDE-based algorithms}\label{sec:PDE}
To conclude this part of the survey, we discuss PDE-based approaches. As discussed in Section~\ref{sec:control_intro}, the optimal control can be identified by solving the corresponding HJB equation \eqref{def:control-HJB}.

The Deep Galerkin Method (DGM) proposed by Sirignano and Spiliopoulos \cite{SiSp:18} aims at approximating the solution of the parabolic PDE with a deep neural network. In fact, the DGM can tackle PDEs of other (potentially nonlinear) types, with terminal (or initial) condition and boundary conditions. To focus on the main ideas, we present the algorithm for a generic PDE on a spatial domain $\dom \subseteq \RR^d$ and a time interval $[0,T]$,
\begin{equation}
\begin{cases}
     \partial_t u(t, x) + \mc{L}u(t, x) = 0, \quad (t, x) \in [0,T] \times \dom , \\
    u(T, x) = u_T(x), \quad x \in \dom , \\
    u(t, x) = \Gamma(t, x), \quad x \in [0,T] \times \partial \dom.
\end{cases}
\end{equation}
Here $\mc{L}$ is an operator in $x$, possibly nonlinear. A Dirichlet boundary condition is considered although other boundary conditions (Neumann, Robin) can also be treated in this framework.  

The DGM algorithm proposes to replace  $u$ by a deep neural network, denoted by $u(t, x; \theta)$, and minimizes the following loss function
\begin{align}
\label{eq:loss-DGM-HJB-PDE} 
    J(\theta) 
    &= \eta \left\| \partial_t u(t, x;\theta)  + \mc{L}u(t, x; \theta) \right\|^2_{L^2([0,T] \times \dom; \mu_1)} 
    + \eta_I \left\| u(T, x; \theta) - u_T(x) \right\|^2_{L^2(\dom; \mu_2)} \\
    &\quad + \eta_{BC} \left\| u(t, x ;\theta) - \Gamma(t, x) \right\|^2_{L^2([0,T] \times \partial\dom; \mu_3)},
\end{align}
where $\mu_i$, $i=1,2,3$, are probability distributions on the corresponding domains, and $\|f(y)\|^2_{L^2(\mc{Y}; \mu)} = \int_{\mc{Y}} |f(y)|^2 \mu(\ud y)$. The first term is for the PDE residual, the second one is for the terminal condition and the third term is for the boundary condition. The positive constants $\eta, \eta_I$ and $\eta_{BC}$ give more or less importance to each component. The differential operators $\partial_t u(t, x; \theta)$ and $\mc{L}u(t, x; \theta)$ can be computed exactly using automatic differentiation.  However, most of the time, the second derivatives are computationally costly of $\mc{O}(d^2 \times N_{\mathrm{Batch}})$ and third-order  derivatives $\grad_\theta \Hess_x u(t, x; \theta)$ are also needed for SGD algorithms. A fast computation of second derivatives using a Monte Carlo method was proposed in \cite[Section 3]{SiSp:18}.  The squared $L^2$ norm of each term in \eqref{eq:loss-DGM-HJB-PDE} is calculated as the average of the squared value evaluated at random points drawn according to respective probability densities. Then we use SGD (see Algorithm~\ref{algo:SGD-generic}) to minimize the loss function $J$ defined in~\eqref{eq:loss-DGM-HJB-PDE}. One sample is $\xi = ((t,x), x', x'') \in ([0,T] \times \dom) \times \dom \times \partial\dom$, picked according to the distribution $\mu_1 \otimes \mu_2 \otimes \mu_3$. 

When applying the DGM to the HJB equation \eqref{def:control-HJB} associated with the control problem \eqref{def:control-Xt}--\eqref{def:control-cost}, one uses the terminal condition $u_T(x) = g(x)$ and omits the boundary condition if none is present.

\begin{remark}
 
The DGM method can, at least in principle, be applied to a large variety of PDEs. Its flexibility is a key advantage. However, many choices need to be made in practice when implementing this method. First, the sampling distributions $(\mu_i)_{i=1,2,3}$  strongly influence the training process and thus the learned function. Choosing a suitable distribution inside the domain and on its boundary is sometimes not trivial, particularly since one needs to sample from this distribution efficiently. Furthermore, when dealing with a loss function composed of several terms (as {\it e.g.}~\eqref{eq:loss-DGM-HJB-PDE}), balancing the various terms can be challenging but is important to ensure that no term dominates the loss or is obliterated by the other terms. If some weight is too small, the neural network will tend to neglect the corresponding term. On the other hand, if some weight is much larger than needed, it will obfuscate the other terms. Overall, the choice of suitable coefficients seems important to efficiently guide the neural network towards a good local minimum. 

BSDE-based algorithms, as introduced in Section~\ref{sec:BSDE}, also necessitate the selection of a sampling distribution. In the global in time approach, this translates into the choice of the forward process. In the local in time approach, one needs to choose how to sample points during training at each time step. However, these methods do not require domain truncation. Monte Carlo simulations are naturally driven by the realizations of the Brownian motion path.
\end{remark}

\begin{remark}[Theoretical analysis]
    Sirignano and Spiliopoulos showed in~\cite[Theorems 7.1 and 7.3]{SiSp:18}, assuming that the PDE is quasilinear parabolic and has a smooth enough solution, that for any level $\epsilon$, there exists a wide-enough neural network that can make the $L^2$ error $J(\theta)$ smaller than $\epsilon$. They also show, under stronger conditions, how one can obtain a sequence of neural networks converging to the PDE solution. To the best of our knowledge, the convergence of this algorithm remains to be studied for more general PDEs.

\end{remark}

The DGM can be applied to PDE system characterizing solutions to mean-field control problems; see~\cite{carmona2019convergence1} in the ergodic setting. For MFC (see Section~\ref{sec:directMethod}) with nonsmooth costs, Reisinger,  Stockinger and Zhang \cite{reisinger2021fast} developed an iterative algorithm  incorporating the gradient information and the proximal map of the nonsmooth cost. We will explain how to adapt the DGM to PDE systems when discussing applications to MFGs in Section~\ref{sec:MFPDE-deeplearning}.

The DGM has also been extended to deal with path-dependent PDEs in~\cite{saporitoZhang2021PDGM}, which used an LSTM network to capture the dependence on the path, and then passed this information to a feed-forward network in charge of learning the PDE solution. It has also been extended to address HJB equations, solving for both the value function and the optimal control simultaneously and characterizing both as deep neural networks \cite{al2022extensions, barnett2023deep, duarte2024machine}. Other closely related works are the physics-informed neural networks (PINNs)~\cite{raissi2019physics}, which is introduced to approximate solutions to equations arising in physics; and the deep Ritz method \cite{yu2018deep} which solves variational problems that arise from PDEs. The key idea is, here again, to look for the solution to the equation among a class of neural networks. The error analysis has been discussed in \cite{van2022optimally,mishra2022estimates,de2022error,de2022error2,de2022generic,mishra2022estimates2} for various types of PDEs and operators. 

\subsection{Extensions}\label{sec:extension}
\subsubsection{Stochastic control with delay}\label{sec:sc_delay}

In direct parameterization approaches as discussed above, the neural network is used to approximate the optimal control or the value function. The type of inputs depends on the class of controls considered in the optimal control problem. For Markovian controls, which are used in \cite{han2016deep-googlecitations} and \cite{bachouch2021deepnumerical}, the input shall be $\checkX_{t_n}$ for $\alpha_{t_n}(\cdot; \theta_n)$. For open-loop controls, one can naturally consider $(\check W_{t_0}, \ldots, \check W_{t_n})$, while for closed-loop controls, one can use $(\checkX_{t_0}, \ldots, \checkX_{t_n})$ as inputs of the neural network. 
However, the sequence length increases with the number of time steps, which leads to a high computational cost. Furthermore, passing the whole sequence as a single input does not make use of the time structure. For these reasons, other architectures have been considered.

We now illustrate the direct parameterization methods by studying SC problems with delay. Such problems have found many applications, \textit{e.g.}, modeling systems with the aftereffect in mechanics and engineering, biology, and medicine \cite[Chapter 1]{kolmanovskiui1996control}, time-to-build problems in economics \cite{kydland1982time,asea1999time}, modeling the ``carryover'' or ``distributed lag''   advertising effect in marketing \cite{gozzi200513,gozzi2009controlled}, and portfolio selection under the market with memory and delayed responses in finance \cite{oksendal2000maximum,federico2011stochastic,elsanosi2001optimal,li2018portfolio}.

Note that feedforward neural network (FNN) architecture is the most common architecture in deep learning and performs well as a function approximator of Markovian controls. Another popular type of neural networks is recurrent neural networks (RNNs). In a study by Han and Hu \cite{han2020rnn}, it is shown that RNNs have better performance in control problems with a delay effect.

To distinguish between the value of a process at a given time and the portion of trajectory ending at time $t$ with $\delta$ history, for any process $P$, we denote by $\underline P_t = (\underline P_t\bigl(s)\bigr)_{s \in [-\delta, 0]}$ the trajectory of $P$ from time $t-\delta$ to $t$, \emph{i.e.}, $\underline P_t(s) = P_{t+s}$, for $-\delta \leq s \leq 0$. Specifically, we consider a SC problem in which the state process $X$ is characterized by a stochastic delay differential equation (SDDE),
\begin{equation}\label{def:delay-Xt}
\begin{dcases}
    \ud X_t = b(t, \underline X_{t}, \ctrl_t) \ud t + \sigma(t, \underline X_t, \ctrl_t) \ud W_t,  & t \in [0, T],\\
    X_t = \varphi_t, & t \in [-\delta, 0],
\end{dcases}
\end{equation}
and the objective function is given by
\begin{equation}\label{def:delay-cost}
J(\alpha) = \EE\left[\int_0^T f(t, \underline X_t, \ctrl_t) \ud t  + g(\underline X_T)\right].
\end{equation}
Here $\delta \geq 0$ is the fixed delay, and $(\varphi_t)_{t \in [-\delta, 0]}$ is a given process on $[-\delta, 0]$ for the initial condition of $X$, and $X_t$ denotes the value of the state process at time $t$ as usual. The SDDE \eqref{def:delay-Xt} has been well studied in the literature \cite{mohammed1984stochastic,mohammed1998stochastic} (see also Appendix~\ref{supp:SDDE} some preliminaries on SDDEs).

The key difficulty is that, when optimizing over closed-loop controls, one should, in principle, take into account the whole trajectory of the state, which is computationally costly. The authors in \cite{han2020rnn} analyzed this problem with a focus on the deep neural networks' (DNNs) architecture design in order to handle the high dimensionality arising from the delay. 
Without loss of generality, let us consider that the fixed delay $\delta <\infty$ covers $N_\delta$ subintervals, {\it i.e.}, $\delta = N_\delta \Delta t$, 
the partition on $[0,T]$ is extended to $[-\delta, 0]$,
\begin{equation}
 \; -\delta = t_{-N_\delta} \leq t_{-N_\delta+1} \leq \ldots \leq t_0 = 0, \text{ with } t_{n+1} - t_n \equiv \Delta t, \; \forall n = -N_\delta, \ldots, -1.
\end{equation}
and the discretized version of \eqref{def:delay-Xt}--\eqref{def:delay-cost} becomes
\begin{align}
    & \check X_{t_{n+1}} = \check X_{t_n} + b(t_n, \check{\underline{X}}_{t_n}, \ctrl_{t_n}) \Delta t + \sigma(t_n, \check{\underline X}_{t_n}, \ctrl_{t_n}) \Delta \check W_{t_n}, \label{def:Xtdiscrete} \\
    & \inf_{\{\ctrl_{t_n}\}_{n=0}^{N_T-1}} \EE\left[\sum_{n=0}^{N_T-1} f(t_n, \check{\underline X}_{t_n}, \ctrl_{t_n}) \Delta t + g(\check {\underline X}_T)\right], \label{def:costdiscrete}
\end{align}
where for each $n$,  $\check {\underline X}_{t_n}$ represents the discrete path with $N_\delta$ lags and $\Delta \check W_{t_n}$ is the increment in Brownian motions,
\begin{equation}
    \check {\underline X}_{t_n} = (\check X_{t_{n-N_\delta}},\ldots, \check X_{t_n}), \quad \Delta \check W_{t_n} = \check W_{t_{n+1}} - \check W_{t_n}.
\end{equation}
Here $\ctrl_{t_n}$ is a function of time and the past state trajectory $\check{\underline X}_{t_n}$ taking values in $\RR^k$. 
The next two architectures are proposed for approximating $\ctrl_{t_n}$. 

\medskip

\noindent \textbf{Feedforward architecture.}
Motivated by the path-dependent structure of the considered problems (the change of the current state only depends on the history up to lag $\delta$), a natural idea is to approximate $\ctrl_{t_n}$ by a feedforward neural network taking the state history up to lag $\bar{\delta}$ as the input. Note here it could be $\bar{\delta}\neq \delta$ since one may not know the underlying true $\delta$ a priori.
Without loss of generality, one can assume $\bar{\delta}=N_{\bar{\delta}}h~ (N_{\bar{\delta}}\in\mathbb{N}^+$) and define
$\bar{X}_{t_n} \equiv (\check X_{t_{n-N_{\bar{\delta}}}},\ldots, \check X_{t_n})\in \RR^{d\times (N_{\bar{\delta}}+1)}$. Then, the control at time $t_n$ can be approximated by a feedforward fully connected network taking $\bar X_{t_n}$ as input. 

\medskip

\noindent\textbf{Recurrent architecture.} 
Alternatively, the sequence $\bar X_{t_n}$ can be processed by a recurrent neural network (RNN), such as a long-short term memory (LSTM) neural network. The basic principle of an RNN is that the elements of the sequence are processed one by one by applying the same neural network in a recurrent manner, and some information is saved between two applications until the output is calculated. Here, we can use a single RNN for all time steps. At time $t_n$, to produce the value of the control, it uses $\check{X}_{t_n}$ as an input, along with the information saved from the previous time step. We refer to Appendix~\ref{sec:RNNdetails} and Appendix~\ref{sec:LSTMdetails} for more details.

\begin{remark}
Although for both feedforward neural networks and RNN, the input dimensions remain constant as $k$ changes, using the former requires prior knowledge of $\delta$. For the feedforward fully connected neural network, one feeds the discretized state values \eqref{def:Xtdiscrete} of length $N_{\bar\delta} + 1$, so to obtain the best performance, one needs to get a good estimate $\bar\delta$ of $\delta$ first. On the other hand, for the RNN, one only needs to provide the current state value $X_{t_n}$. Notice that in an LSTM all input information up to time $t_n$ is summarized by the $n^{th}$ cell, but if the optimal control depends only on the past up to $\delta$, the forget gates allow to drop out the unneeded information. The exact way information should be dropped is determined by the neural network parameters training. %
Though the authors in \cite{han2020rnn} only experimented with LSTM, other variations of RNNs such as gated recurrent units (GRUs)~\cite{cho2014learning} or peephole LSTM~\cite{gers2002learning} can be considered. 
\end{remark}

Before illustrating the above method, we stress that other methods also exist. For example, \cite{lefebvre2021linear}, a PDE-based approach,  uses a deep learning technique inspired by the physics-informed neural networks (PINNs)  and applies it to solve a mean-variance portfolio selection with execution delay. One may also reformulate the delay problem using anticipated backward stochastic differential equations (ABSDEs) and solve it through BSDE-based methods. Along this line, \cite{fouque2019deep} addressed mean-field control with delay.

\paragraph*{Numerical illustration: a linear-quadratic regulator problem with delay.}  
LQ problems with delay were first investigated by Kolmanovski{\u{\i}} and Sha{\u{\i}}khet \cite{kolmanovskiui1996control}. In the version presented here, the aim is to minimize 
\begin{align}
&\EE_{\varphi} \left[\int_0^T (X_t + e^{\lambda \delta}A_3 Y_t)\transpose Q(t) (X_t + e^{\lambda \delta}A_3 Y_t) + \ctrl_t\transpose R(t)\ctrl_t\ud t \right.+ (X_T + e^{\lambda\delta} A_3Y_T)\transpose G (X_T + e^{\lambda\delta} A_3Y_T)\Bigg], \\
&\text{subject to } \ud X_t = (A_1(t)X_t + A_2(t)Y_t + A_3 Z_t + B(t)\ctrl_t) \ud t + \sigma(t) \ud W_t, \quad t  \in[0,T],\label{eq:lq}
\end{align}
where $X_0 = \varphi \in L^2(\Omega, C([-\delta, 0], \RR^d))$ is a given initial segment, $Y_t = \int_{-\delta}^0 e^{\lambda s} X_{t+s} \ud s$ is the distributed delay and $Z_t = X_{t-\delta}$ is the discrete delay,
$A_1, A_2, Q \in L^\infty([0,T]; \RR^{d\times d})$, $B \in L^\infty([0,T]; \RR^{d \times k})$, $R \in L^\infty([0,T]; \RR^{k\times k})$ are deterministic matrix-valued functions, $\sigma \in L^2([0,T]; \RR^{d \times m})$ is a deterministic matrix-valued function, $A_3, G \in \RR^{d \times d}$ are deterministic matrices. It is assumed that $Q(t), G$ are positive semi-definite and $R(t)$ is positive definite for all $t \in [0,T]$ and continuous on $[0,T]$. To have a tractable solution, a further relation is prescribed,
\begin{equation}\label{lq:parameters}
A_2(t) = e^{\lambda\delta}(\lambda I_d + A_1(t) + e^{\lambda\delta} A_3) A_3,
\end{equation}
where $I_d$ is the identity matrix with rank $d$. This example was studied in
\cite[Section 4]{bauer2005stochastic}, and the main results are also summarized in \cite{han2020rnn}. 

We present results below for a ten-dimensional example. The model parameters are chosen as follows. The dimensions are $d = k = m = 10$, and $\lambda=0.1$.
In~\eqref{eq:lq}, $A_1, A_3, B, \sigma$ are constant coefficient matrices (generated randomly), $Q, R, G$ are constant matrices proportional to identity matrices, and $A_2$ is determined by~\eqref{lq:parameters}. For implementation details, we refer to~\cite{han2020rnn}. The left panel of Figure~\ref{fig:lq_train_curve} displays the total cost of the validation data against training time. The values are averaged every 200 steps. The feedforward model takes the state history as inputs up to lag $\bar{\delta}=\delta$ with $N_{\bar{\delta}}=40$.
The right panel of Figure~\ref{fig:lq_train_curve} displays the optimized cost as a function of the processed lag time $\bar{\delta}$ from $0.2$ to $1$ with step size 0.1, while the actual lag $\delta=1$. If the chosen lag time $\bar{\delta}$ is smaller than the actual lag $\delta$ time, there is a loss of information when the feedforward network processes the data. As expected, we observe that the cost increases as the lag time processed by the feedforward model decreases. A higher optimized cost indicates that the model can only find a strictly sub-optimal strategy due to the lack of information. 
Figure~\ref{fig:lq_path_lstm_shff} displays one sample path (first five dimensions only) of the optimal state $X$ and control $\ctrl$ provided by two neural networks in comparison with the analytical solution, in which the LSTM architecture presents a better agreement. 
The lag time $\bar{\delta}$ processed by the feedforward model is chosen to be the same as $\delta$.
One main drawback of the feedforward model is that it requires to know the true lag time $\delta$ a priori to determine the network's size.

\begin{figure}[!htb]
\centering
\includegraphics[width=0.45\textwidth]{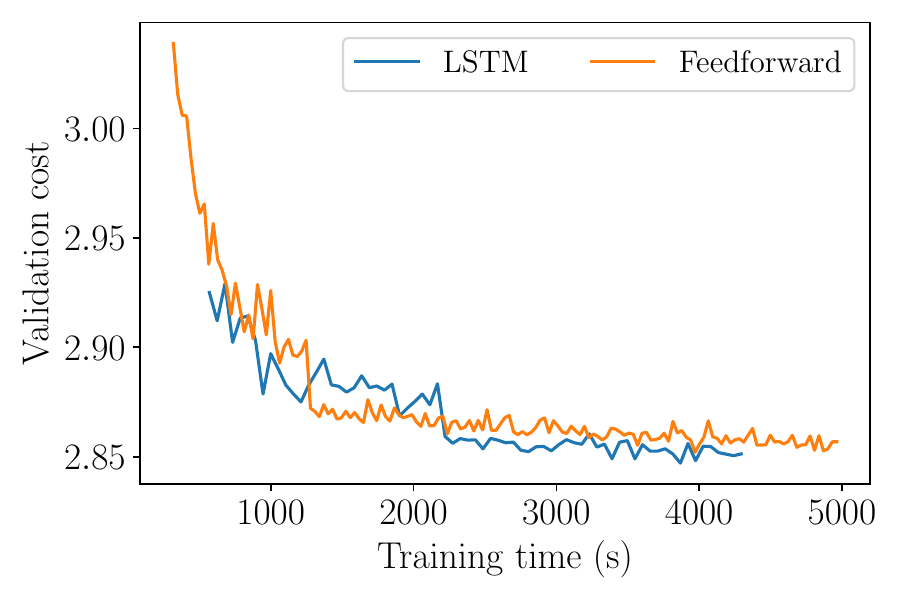}
\includegraphics[width=0.45\textwidth]{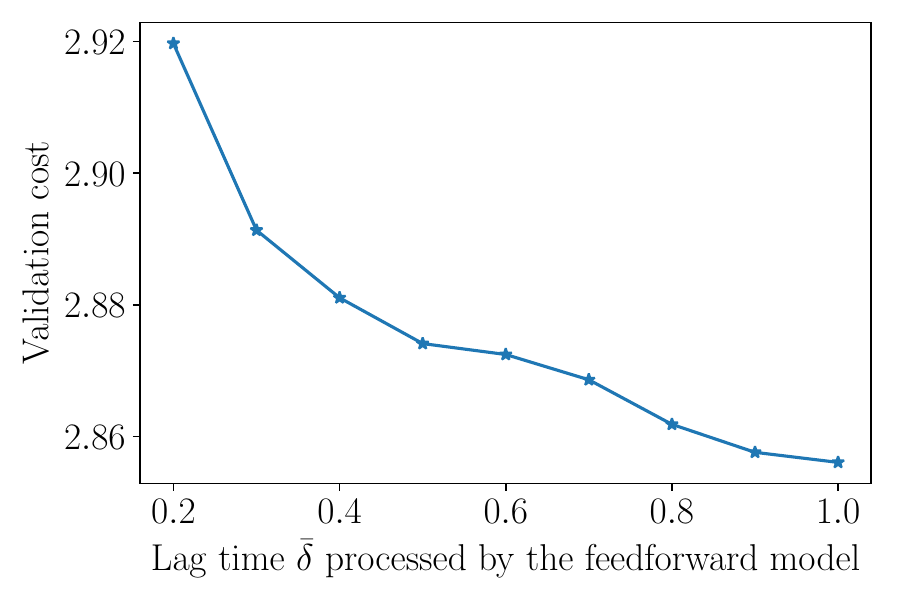}
  \caption{The linear-quadratic regulator problem with delay in Section~\ref{sec:sc_delay}. Left: Training curve of two models in the example of linear-quadratic problem. Right: The effect of lag time $\bar{\delta}$ processed by the feedforward model in the example of the linear-quadratic problem. The lag time $\delta$ in the actual system is $1$. %
  }
  \label{fig:lq_train_curve}
\end{figure}

\begin{figure}[!htb]
\centering
\includegraphics[width=0.98\textwidth, trim = {0, 18em, 0, 18em}, clip]{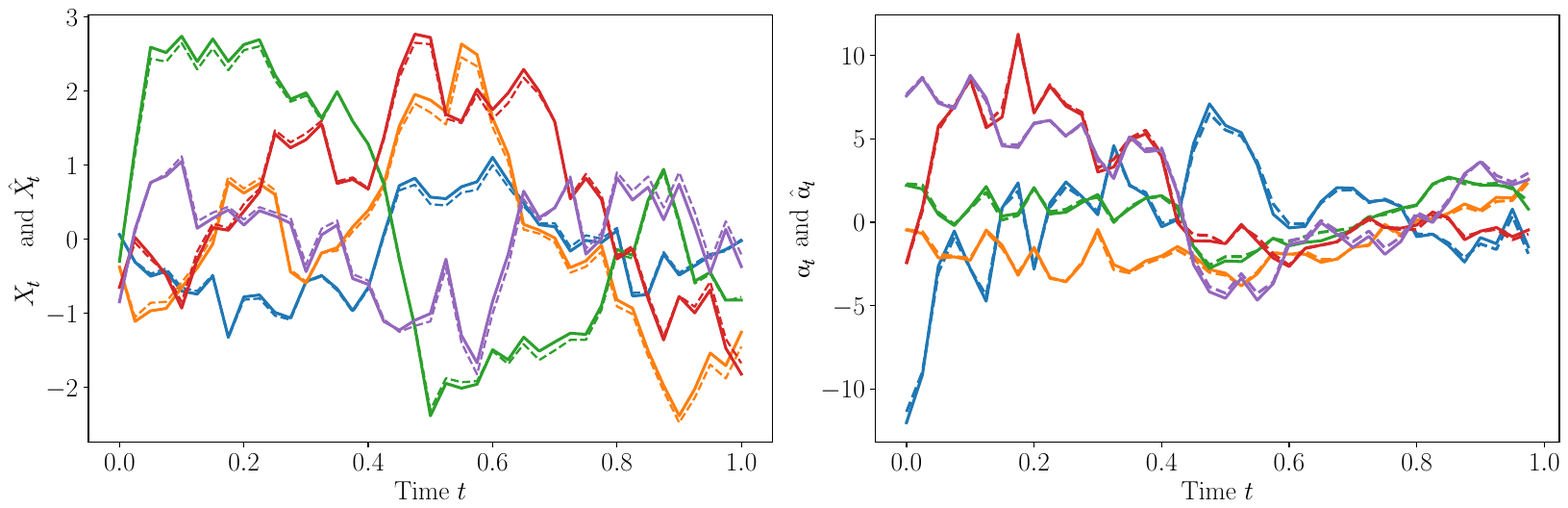}\\
\includegraphics[width=0.98\textwidth,  trim = {0, 18em, 0, 18em}, clip]{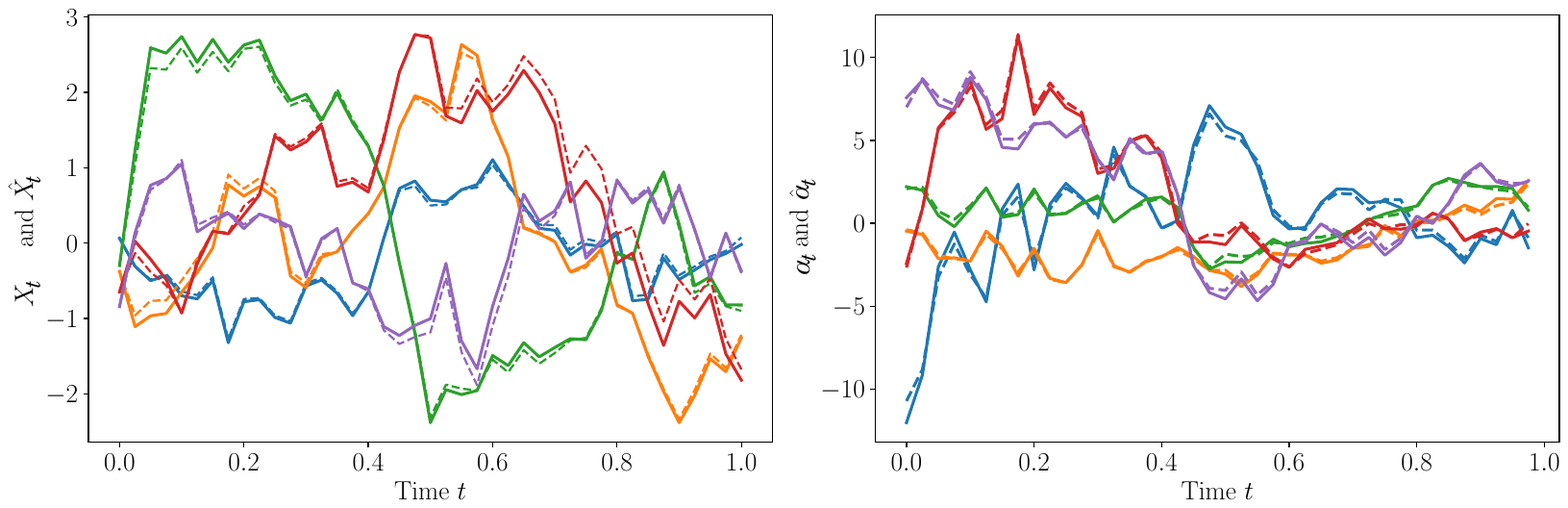}
  \caption{The linear-quadratic regulator problem with delay in Section~\ref{sec:sc_delay}. 
      A sample path of the first 5 dimensions of the state $X_t$ and control $\ctrl_t$ obtained from the LSTM (top) model and FNN (top) model.
      Left: the optimal state process discretized from the analytical solution $(X_t)_i$ (solid lines) and its approximation $(\hat{X}_t)_i$ (dashed lines) provided by the approximating control, under the same realized path of Brownian motion.
      Right: comparisons of the optimal control $(\ctrl_t)_i$ (solid lines) and $(\hat{\ctrl}_t)_i$ (dashed lines). %
  }
  \label{fig:lq_path_lstm_shff}
\end{figure}

\subsubsection{Mean-field type control}\label{sec:directMethod}
We now discuss how direct parameterization methods can be adapted for mean field control (MFC) problems, which are an extension of standard optimal control in which there are mean field interactions. MFC can be interpreted as a social optimum but also has applications in risk, {\it e.g.}, heating or electric loads management~\cite{KIZILKALE20141867,Mathieu2013StateEA},  risk management in finance~\cite{MR2784835} or optimal control with a cost involving a conditional expectation~\cite{achdoulaurierelions2020optimal,MR4133380}, to cite just a few examples. We refer to~\cite{MR3134900} for more background on MFC.

As before, let $\mc{A} \subset \RR^k$ denote the set of admissible actions and let $\mathbb{A}$ be the set of so-called admissible controls. We denote by $\cP_2(\RR^d)$ the set of probability measures on $\RR^d$ which admit a second moment. Let $b, \sigma, f, g$ be Borel-measurable functions,
\begin{equation}
    (b, \sigma, f) : [0,T] \times \RR^d \times \cP_2(\RR^d) \times \mc{A} \to (\RR^d, \RR^{d \times m}, \RR), \quad g : \cP_2(\RR^d) \times \RR^d \to \RR.
\end{equation}
Compared with standard optimal control, see~\eqref{def:oc-bsigmaf}, the functions take as an extra input the population distribution.

\begin{defn}[MFC optimum]
\label{def:mfc-optimum}
A feedback control $\ctrl^* \in \mathbb{A}$ is an optimal control for the MFC problem for a given initial distribution $\mu_0 \in \cP_2(\RR^d)$ if it minimizes
\begin{align}
\label{subchapRCML-num-eq:def-J-MFC}
	J(\alpha): \ctrl \mapsto \EE \left[\int_0^T f(t, X_t, \mu_t, \ctrl_t ) dt + g(X_T, \mu_T) \right],
\end{align}
where $\mu_t$ is the probability distribution of the law of $X_t$, under the constraint that the process $X = (X_t)_{t \ge 0}$ solves the stochastic differential equation of the McKean-Vlasov (MKV) type,
\begin{equation}
\label{subchapRCML-num-eq:dyn-X-general-MFC}
	\ud X_t = b(t, X_t, \mu_t, \ctrl_t) \ud t + \sigma(t, X_t, \mu_t, \ctrl_t) \ud W_t, \qquad t \ge 0,
\end{equation}
with $X_0$ having distribution $\mu_0$. 
\end{defn}
Existence and uniqueness results have been studied in the literature, and we refer the interested reader to {\it e.g.}~\cite{carmona2018probabilistic2} and in particular the assumption  ``Control of MKV Dynamics'' on page 555. Furthermore, under mild conditions, there is an optimal control that is Markovian, so it is sufficient to focus on controls that are functions of time and the individual state. The framework can also encompass problems in which the interactions are through the joint distribution of states and actions and not just the distribution of states. Such problems are sometimes referred to as extended MFC problems. Although the numerical method discussed below can be adapted in a rather straightforward way (as illustrated below in an example below), the theoretical framework is more challenging.

To solve the MFC problem using deep learning, one can for example follow the lines of the global approach described in Section~\ref{sec:direct-global-local} to discretize time and to replace the control by a sequence of neural network functions of the state. However, in contrast with standard OC problems, here, the costs and the dynamics can depend on mean field terms, so we also need to approximate the distribution of the state. To this end, The simplest approach is to use the empirical distribution of a population of particles. 

Following this approach, the original MFC problem is approximated by the following problem: minimize over the neural network parameters $\theta = (\theta_n)_{n=0,\dots,N_T-1}$
\begin{equation}
\label{subchapRCML-num-eq:NN-MFC-cost-totalapprox}
	\check J^{N}(\theta): \theta \mapsto \frac{1}{N} \sum_{i=1}^N \EE \left[\sum_{n=0}^{N_T-1} f(t_n, \check X_{t_n}^{i, \theta}, \check \mu^{N, \theta}_{t_n}, \ctrl_{t_n}(\check X_{t_n}^{i, \theta}; \theta_n) ) \Delta t + g(\check X_T^{i, \theta}, \check \mu^{N,\theta}_T) \right],
\end{equation}
subject to
\begin{equation}
\label{subchapRCML-num-eq:NN-MFC-Nparticles-Deltat}
	\check X_{t_{n+1}}^{i, \theta} = \check X_{t_{n}}^{i, \theta} + b(t_n,  \check X_{t_{n}}^{i, \theta}, \check \mu^{N, \theta}_{t_{n}}, \ctrl_{t_n}(\check X_{t_{n}}^{i, \theta}; \theta_n)) \Delta t + \sigma(t_n, \check X_{t_{n}}^{i, \theta}, \check \mu^{N, \theta}_{t_{n}}, \ctrl_{t_n}(\check X_{t_{n}}^{i, \theta}; \theta_n)) \Delta \check W^i_n, \quad n = 0, \dots, N_T-1,
\end{equation}
where the initial positions $(\check X_0^{i, \theta})_{i=1,\dots,N}$ are i.i.d. with distribution $\mu_0$, the empirical distribution $\mu^{N, \theta}_{t_{n}}$ is 
$$
	\check \mu^{N, \theta}_{t_{n}} = \frac{1}{N} \sum_{j=1}^N \delta_{\check X_{t_{n}}^{j, \theta}}, 
$$ 
and $(\Delta \check W_n^i)_{i=1,\dots,N,\;n=0,\dots,N_T-1}$ denotes i.i.d. random variables with Gaussian distribution $\mathcal N(0, \Delta t)$. Note that this problem is not equivalent to solving an $N$-agent optimal control problem since, in the latter case, the control would in general need to be a function of all the agents' states, while here, we consider only distributed controls, functions of each agent's state. Intuitively, this approximation is justified by the propagation of chaos results~\cite{sznitman1991topics}: as $N\to+\infty$ the particles become independent, and the empirical distribution converges to the distribution of the MKV SDE~\eqref{subchapRCML-num-eq:dyn-X-general-MFC}. See~\cite{lacker2017limit} in the context of MFC problems. 

This technique was proposed in~\cite{carmona2019convergence2} (with a single neural network function of time and space instead of a sequence of neural network functions of space online). Although we focus here on a basic setting for which a simple feedforward fully connected architecture performs well, other architectures may yield better results for problems with more complex time dependencies; see \textit{e.g.}, \cite{fouque2019deep,gomes2022machineprice} for applications with RNNs. See also~\cite{germain2019numerical,agrambakdioksendal2020deep} and the survey~\cite{carmonalauriere2021deepmfgsurvey} for more details.

\begin{rem}[Theoretical analysis]
In \cite[Theorem 3]{carmona2019convergence2} Carmona and Lauri\`{e}re provided a proof of convergence  of the discrete problem to the original MFC problem in the following sense: under suitable assumptions, the difference between the optimal value of this problem $\inf_\theta \check J^{N}(\theta)$ and the optimal value of the original problem $\inf_\ctrl J(\ctrl)$ goes to $0$ as $N_T$, $N$ and the number of parameters in the neural network go to infinity. 
\end{rem}

From here, one can use SGD or one of its variants to optimize $\check J^N$ over the parameters. In contrast with the methods discussed previously for standard OC, here one sample corresponds to one population with $N$ particles. In this context, one sample corresponds to $\xi = ( \check X^{i}_0, (\Delta \check W^i_n)_{n=0,\dots,N_T-1})_{i=1,\dots,N}$, from which we can compute the $N$ (interacting) trajectories and the total empirical cost:
\begin{equation}
\label{subchapRCML-num-eq:NN-MFC-cost-totalapprox-oneS}
	\check J^{\xi, N}(\theta) = \frac{1}{N} \sum_{i=1}^N \left[\sum_{n=0}^{N_T-1} f(\check X_{t_n}^{i, \theta, \xi}, \check \mu^{N, \theta, \xi}_{t_n}, \ctrl_{t_n}(\check X_{t_n}^{i, \theta, \xi}; \theta_n) ) \Delta t + g(\check X_T^{i, \theta, \xi}, \check \mu^{N, \theta, \xi}_T) \right].
\end{equation}
From here, Algorithm~\ref{algo:SGD-generic} in Appendix~\ref{sec:SGD-var} can be applied.

Before illustrating the above method, we would like to note that there are both PDE-based and BSDE-based algorithms for solving MFC problems. These involve reformulating them into HJB and FBSDE systems. As these systems share similarities with those derived from MFG, interested readers are directed to the corresponding sections in Section~\ref{sec:MFG}.
 
\paragraph*{Numerical illustration: a mean-field price impact problem.} We now illustrate the method with a financial application on a problem of optimal execution. The model describes a large group of traders interacting through the price of an asset. The large group of traders has a non-negligible influence on the price, referred to as price impact. For example, when the traders decide to buy at the same moment, the price is driven up, and vice versa when the traders decide to sell simultaneously.

The $N$-agent problem was originally solved as a mean-field game (MFG) by Carmona and Lacker in the weak formulation (\cite{MR3325272}), and revisited in the book of Carmona and Delarue  \cite[Sections 1.3.2 and 4.7.1]{carmona2018probabilistic} in the strong formulation. Here, we focus on the mean-field control setting. In this model, the inventory representative trader is denoted by $X_t$. The control of the trader is denoted by $\ctrl_t$ and corresponds to the trading rate. The dynamics of the inventory are given by
\begin{equation*}
    \ud X_t = \ctrl_t \ud t +\sigma \ud W_t,
\end{equation*}
where $W$ is a standard Brownian motion, and the cost can be rewritten in terms of $X$ only as
\begin{equation*}
J(\ctrl)=\mathbb{E}\Big[ \int_0^T \left(c_{\ctrl}(\ctrl_t)+c_X(X_t)-\gamma X_t \int_{\mathbb{R}} a  \ud \nu^\ctrl_t(a)\right)\ud t +g(X_T)\Big].  
\end{equation*}
Following the Almgren-Chriss linear price impact model, we assume that the functions $c_X$, $c_{\ctrl}$ and $g$ are quadratic. Thus, the cost is of the form
\begin{equation}\label{eq:priceimpact}
J(\alpha)=\mathbb{E}\left[ \int_0^T \left( \frac{c_{\alpha}}{2}{\alpha_t}^2+\frac{c_X}{2}X_t^2-\gamma X_t\int_{\mathbb{R}} a \ud\nu^\ctrl_t(a) \right)\ud t + \frac{c_g}{2}X_T^2\right].
\end{equation}
This model has a semi-explicit solution obtained by reducing the problem to a system of ordinary differential equations (ODEs) as explained in~\cite[Sections 1.3.2 and 4.7.1]{carmona2018probabilistic}.

The deep learning method described above can readily be adapted to solve MFC with interactions through the control's distribution by computing the empirical distribution of controls for an interacting system of $N$ particles. Figure~\ref{fig:ex-price-impact-1} shows the control and the distribution at various time steps. Here, we plot the values of the control represented by the neural network evaluated at samples generated by following the $N$-interacting particle system. We see that the shape is approximately linear at every time step, and furthermore that it coincides with the lines corresponding to the theoretical optimal solution. The distribution, represented by histograms computed using the $N$ particles, starts on the right and moves towards $0$, which can be interpreted as the fact that the traders liquidate their portfolios. We used the parameters: $T=1$, $c_X = 2$, $c_{\ctrl} = 1$, $c_g = 0.3$, $\sigma = 0.5$ and $\gamma = 0.2$.  When using a larger value of $\gamma$, the collective influence of the traders on the price is higher. As can be seen 
in Figure~\ref{fig:ex-price-impact-2} obtained with $\gamma=1$, the traders do not liquidate all the time from $t=0$ until time $t=T$. Instead, they liquidate their portfolio at the beginning and then start buying. We can explain this behavior as follows: since this is a cooperative problem, the traders can use the price impact to drive the price up to increase the value of the stock they own, thus increasing their final reward. 

\begin{figure}[ht]
  \subfloat[]{
  \includegraphics[width=0.45\linewidth]{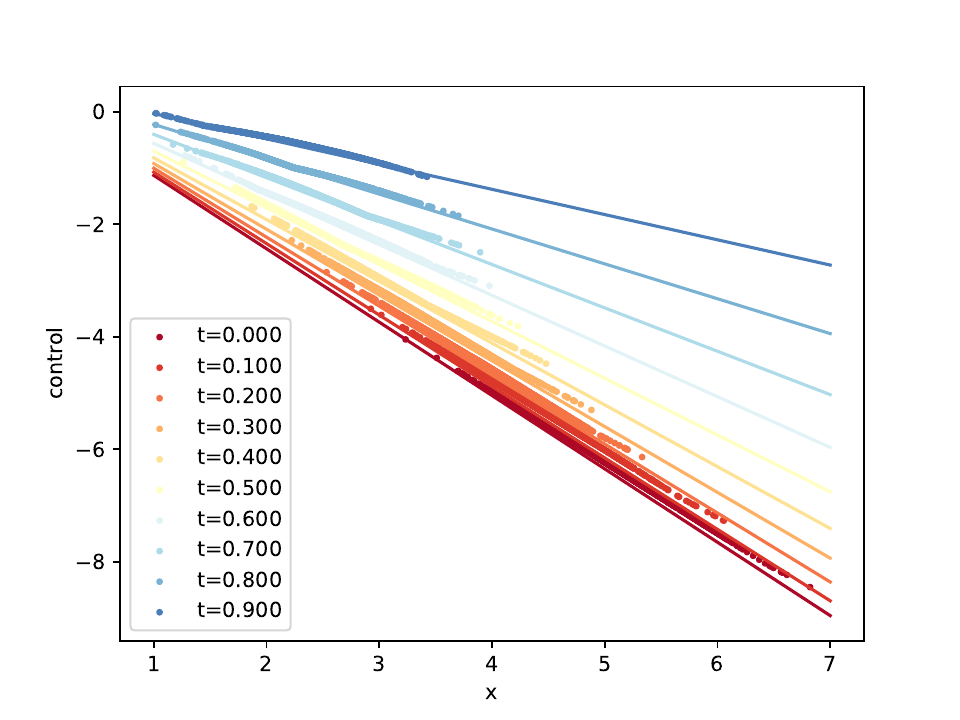}
    }
  \subfloat[]{
  \includegraphics[width=0.45\linewidth]{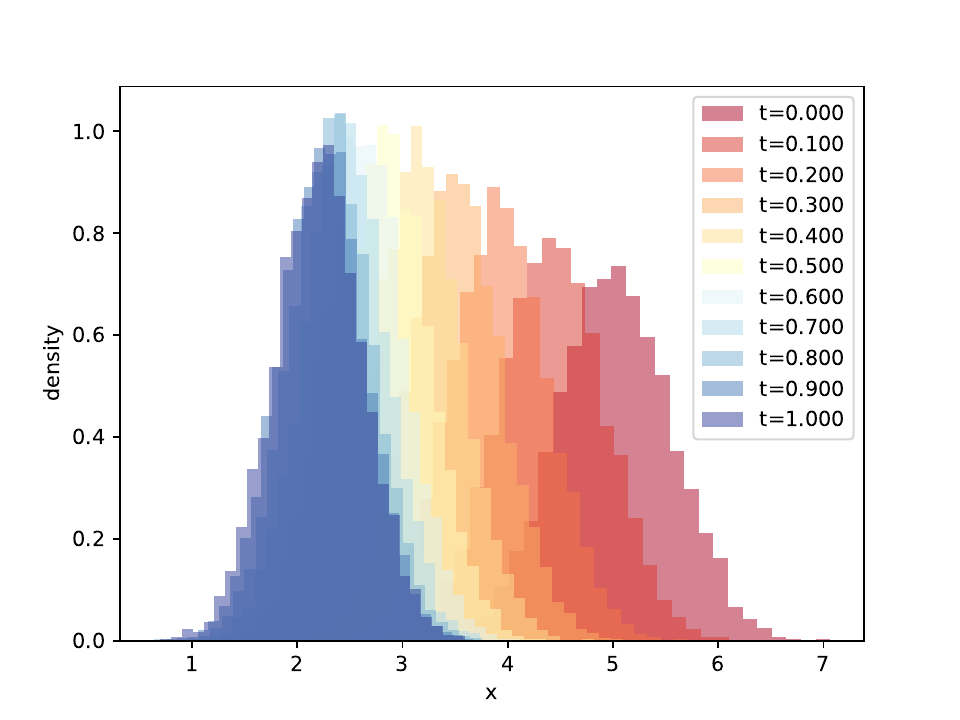}
    }
\caption{Price impact MFC example in Section~\ref{sec:directMethod} solved by direct method. Left: Control learnt (dots) and exact solution (lines). Right: associated empirical state distribution. Here we take $\gamma = 0.2$ in \eqref{eq:priceimpact}. %
}
\label{fig:ex-price-impact-1}
\end{figure}

\begin{figure}[ht]
  \subfloat[]{
  \includegraphics[width=0.45\linewidth]{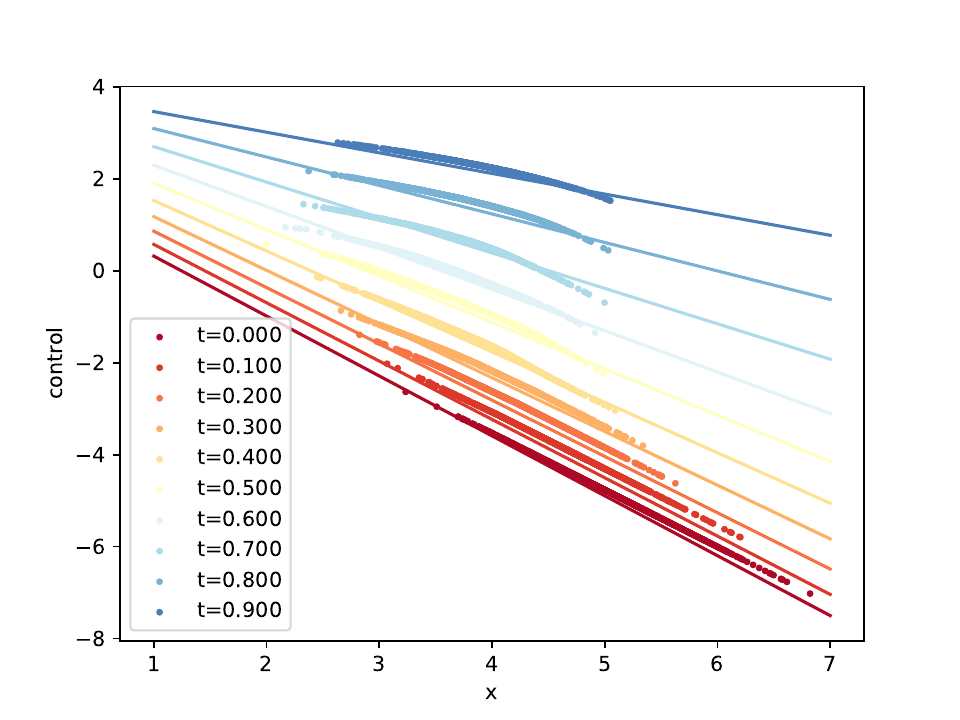}
    }
  \subfloat[]{
  \includegraphics[width=0.45\linewidth]{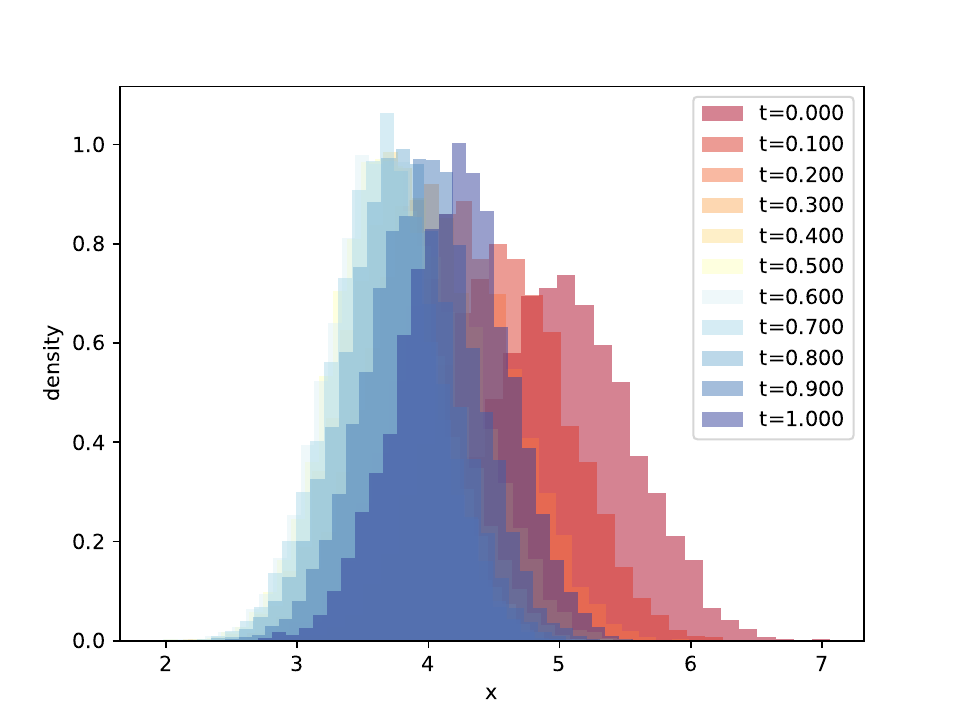}
    }
\caption{Price impact MFC example in Section~\ref{sec:directMethod} solved by direct method. Left: Control learned (dots) and exact solution (lines). Right: associated empirical state distribution. Here, $\gamma = 1$ in \eqref{eq:priceimpact}. %
}
\label{fig:ex-price-impact-2}
\end{figure}

\section{Stochastic Differential Games}\label{sec:SDG}

The previous section is dedicated to studying how a single agent makes strategic decisions in a random environment. This section will focus on stochastic differential games, which model and analyze the interactions between multiple rational agents in a dynamical random system; see \textit{e.g.}~\cite{bacsar1986tutorial,isaacs1999differential,yong2014differential,carmona2016lectures}, and \cite[Chapter 16]{bensoussan2018estimation} for more background on differential games. In games, an important concept is the so-called Nash equilibrium, which refers to a situation in which no one has an incentive to deviate. Finding a Nash equilibrium is one of the core problems in noncooperative game theory. However, the computation of Nash equilibria is extremely time-consuming and memory-demanding for large populations of players \cite{DaGoPa:2009}. 

When the number of player $N$ becomes extremely large, the recently developed theory of mean-field game (MFG)  \cite{HuMaCa:06,HuCaMa:07,LaLi1:2006,LaLi:2007} provides an approximation to $N$-player Nash equilibria if individuals interact symmetrically \cite{CaDe:13}. Despite the huge model reduction from modeling $N$ players to 
one representative player interacting with the population distribution, the MFG itself is still hard to solve numerically if this representative player's state is in a high-dimensional space, or if its evolution involves delay features, common noise, or complicated constraints.

A rich literature on game theory has been developed to study the consequences of strategies on interactions between a large group of rational ``agents'', {\it e.g.}, system risk caused by inter-bank borrowing and lending \cite{CaFoSu:15}, price impacts imposed by agents' optimal liquidation \cite{cardaliaguet2017mean}, market price from the monopolistic competition, and optimal investment with relative performance concerns \cite{LaZa:17,lacker2020many,hu2021n}. This makes it crucial to develop efficient theory and fast algorithms for computing Nash equilibria of $N$-player stochastic differential games and MFGs.

Various numerical methods for differential games have been proposed in the literature. In particular, some of the classical methods for stochastic optimal control problems mentioned at the end of Section~\ref{sec:control_intro} have been extended to treat games, since each player needs to solve an optimal control problem. Differential games typically leads to systems of HJ or HJB equations, which can be tackled using classical numerical schemes for PDEs ({\it e.g.} ~\cite{falcone2006numerical,bardi1999stochastic}). The case of pursuit-evation games has attracted particular attention and leads to Hamilton-Jacobi-Isaacs equations, which can also be treated using classical methods such as semi-Lagrangian schemes ({\it e.g.}~\cite{bardi1999numerical}). Nash equilibria can also be characterized using systems of FBSDEs, which can be solved numerically by adapting existing methods for FBSDEs ({\it e.g.}~\cite{exarchos2019stochastic}). The Markov chain approximation method for optimal control problems has also been extended to games in various settings~\cite{kushner2002numerical,kushner2004numerical,kushner2007numerical}. Besides classical differential games, MFGs have also attracted a growing interest from the numerical viewpoint. Numerical methods to solve the PDE system typically rely on a finite-difference scheme introduced by Achdou {\it et al.}~\cite{MR2679575,MR3097034}; extensions include~\cite{achdouKobeissi2020mfgcfinitediff} for MFGs of control and~\cite{camilli2021approximation} for MFGs-type systems with time-fractional derivatives. Semi-Lagrangian schemes have also been investigated~\cite{MR3148086,MR3392626,MR3828859}. A scheme based on Markov chain approximation has been proposed in~\cite{bayraktar2018numerical}. When the MFG is of variational form, optimization methods can be applied; see~\cite{MR3395203,MR3772008,MR3731033,BricenoAriasetalCEMRACS2017}.

As discussed in the previous sections, deep neural networks with many layers have been shown to solve efficiently stochastic control problems.  
These deep learning techniques also brought the possibility of computing equilibria in high-dimensional games, although extra difficulties arise since in general Nash equilibria are fixed point problems and are not optimization problems. 

Next, we review deep learning algorithms for solving stochastic (moderately large) $N$-player games, and mean-field games as an approximation for extremely large $N$-player games. In particular, the strategy of {\it deep fictitious play}, which integrates fictitious play (a learning scheme in game theory \cite{Br:49,Br:51}) into deep neural network designs, will be used frequently to develop parallelizable deep learning algorithms for computing Nash equilibria of stochastic differential games.

\subsection{$N$-player stochastic games}\label{sec:Nplayer}

\subsubsection{Theoretical background}

We first  consider a stochastic differential game with $N$ players indexed by $i \in \mathcal{I} =  \{1, 2, \ldots, N\}$. Each player has a state process $X_t^i \in \RR^d$ and takes an action $\alpha_t^i$ in the action set $\mc{A} \subset \RR^k$. In the context of games, we will use interchangeably the terms \emph{control} and \emph{strategy}. The dynamics of the controlled state process $X^i$ on $[0,T]$ are given by
\begin{equation}\label{def_Xt_general}
\ud X_t^i = b^i(t, \bm{X}_t, \bm{\alpha}_t) \ud t + \sigma^i(t, \bm{X}_t, \bm{\alpha}_t) \ud W_t^i + \sigma^0(t, \bm{X}_t, \bm{\alpha}_t) \ud W_t^0, \quad X_0^i = x^i, \quad i \in \mc{I},
\end{equation}
where $\bm{W} =[W^0, W^1,\ldots, W^N]$ are $(N+1)$ $m$-dimensional independent Brownian motions, 
and we shall call $W^i$ the individual noises and $W^0$ the common noise. Although at the $N$-player level, adding an extra common noise process is equivalent to choosing $N$ correlated Brownian motions, the current formulation is more convenient when passing the limit $N \to \infty$.
$(b^i, \sigma^i)$ are deterministic functions: $[0,T] \times \RR^{dN} \times \mc{A}^N \to \RR^d \times \RR^{d\times m}$. 
The $N$ dynamics are coupled since all states $\bm{X}_t = [X_t^1, \ldots, X_t^N]$ and all the controls\footnote {Although in the literature of financial mathematics, one usually models $b^i$ and $\sigma^i$ to only depend on player $i$'s own action, it is common in the economics literature that player $i$'s state is also influenced by others' actions. For instance, $\alpha_t^i$ and $X_t^i$ can be respectively the price set by the company and the quantity it produces. To be general, we include this feature in our model, yielding \eqref{def_Xt_general}. Note that by choosing $b^i$ and $\sigma^i$ properly, one can reduce it to the simpler case where each player controls her \emph{private} state through $\alpha^i_t$.}  $\bm{\alpha}_t = [\alpha_t^1, \ldots, \alpha_t^N]$ affect the drifts $b^i$ and diffusion coefficients $\sigma^i$. Given a set of strategies $(\bm{\alpha}_t)_{t \in [0,T]}$, the cost associated to player $i$ is of the form
\begin{equation}\label{def:DFP-open-cost}
J^i(\bm{\alpha}) = \EE\left[\int_0^T f^i(t, \bm X_t, \bm\alpha_t) \ud t + g^i(\bm X_T)\right],
	\end{equation}
where the running cost $f^i: [0,T] \times \RR^{dN}\times \mc{A}^N \to \RR$ and terminal cost $g^i: \RR^{dN} \to \RR$ are deterministic measurable functions. 

Player $i$ chooses $(\alpha_t^i)_{t \in [0,T]}$ to minimize her cost $J^i(\balpha)$ within a set $\mathbb{A}^i$ of admissible strategies. As in Section~\ref{sec:SCP}, the set $\mathbb{A}^i$ usually describes the measurability and integrability of $\alpha_t^i$. Different measurabilities correspond to different information structures available to the players, leading to different types of solutions to the game. In the noncooperative setting, the notion of optimality is the so-called \emph{Nash equilibrium} (NE), and the three main types are open-loop NE ($\bm{W}_{[0,t]}$),  closed-loop NE ($\bm{X}_{[0,t]}$), and closed-loop in feedback form NE ($\bm{X}_t$) also known as Markovian Nash equilibrium.  In this section, we will focus on the open-loop case (Section ~\ref{sec:open-loop}) and the closed-loop in feedback form case (Section~\ref{sec:MarkovianNE}).

\paragraph*{Markovian Nash equilibrium \emph{vs.} open-look Nash equilibrium}
While both equilibria represent strategic solutions in multi-agent decision-making scenarios, they differ in terms of their approach to modeling and predicting the actions of rational agents. Open-loop Nash equilibria are characterized by agents making decisions without considering the current state or actions of other agents, relying solely on pre-determined strategies. These strategies are only adapted to the evolution of the background noises. On the other hand, Markovian Nash equilibria involve agents that base their actions on the current state of the system, as well as the actions of other agents. The ease of solving open-loop and Markovian equilibria depends on the specific problem and the system's complexity. In general, it is easier to analyze Nash equilibria in games with open-loop controls compared to Markovian equilibria. This is because open-loop strategies are independent of other agents' actions.  FBSDEs are particularly well suited for open-loop controls. However, 
an open-loop policy is in principle a function of all the past noises, which can be hard to approximate numerically and it can also make it hard to interpret the agents' decisions. 
In contrast, Markovian controls are easier to interpret since they are functions of the current state, and they also allow to use tools such as the dynamic programming principle and PDEs. Therefore, PDE-based machine learning algorithms are efficient in Markovian control scenarios. For further discussion based on a linear-quadratic model, we refer readers to \cite[Section~5.4]{Hu2:19}.

\paragraph*{Notations } Before proceeding to the open-loop case, we first summarize some commonly used notations as below. Given a probability space $(\Omega, \MCF, \PP)$, we consider 
\begin{itemize}
\item $\bm{W} = [W^0, W^1, \ldots, W^N]$, a $(N+1)$-vector of $m$-dimensional independent Brownian motions;
\item $\mathbb{F} = \{\MCF_t, 0\leq t \leq T\}$, the augmented filtration generated by $\bm{W}$;
\item $\mathbb{H}^{2}_T(\RR^d)$, the space of all progressively measurable $\RR^d$-valued stochastic processes $\alpha: [0,T] \times \Omega \to \RR^d$ such that $\ltwonorm{\alpha}^2 = \EE[\int_0^T \abs{\alpha_t}^2 \ud t ] < \infty$;

\item $\bm{\alpha} = [\alpha^1, \alpha^2, \ldots, \alpha^N]$, a strategy profile, {\it i.e.}, a collection of all players' strategies. The notation $\bm{\alpha}^{-i} = [\alpha^1, \ldots, \alpha^{i-1}, \alpha^{i+1}, \ldots, \alpha^{N}]$ with a negative superscript means the strategies of all players excluding player $i$'s. If a non-negative superscript $\mt{k}$ appears ({\it e.g.}, $\bm{\alpha}^\mt{k}$), this $N$-tuple stands for the strategies from stage $\mt{k}$.  $\bm{\alpha}^{-i,\mt{k}} =  [\alpha^{1,\mt{k}}, \ldots, \alpha^{i-1,\mt{k}}, \alpha^{i+1,\mt{k}}, \ldots, \alpha^{N,\mt{k}}]$ is a $(N-1)$-tuple representing strategies excluding player $i$ at stage $\mt{k}$. We use the same notations for other stochastic processes ({\it e.g.}, $\bm{X}^{-i}, \bm{X}^\mt{k}$).

\end{itemize}

In the rest of this subsection, we will describe how ideas developed for stochastic optimal control can be adapted to solve $N$ player games. We will discuss the direct parameterization, the BSDE-based approach and the PDE-based approach.

\subsubsection{Direct parameterization}
\label{sec:open-loop}

The idea of the direct parameterization approach discussed in Section~\ref{sec:control-direct} can be readily adapted to games, once the information structure ({\it i.e.}, class of controls) for each player is chosen. We have already discussed the direct parameterization for closed-loop feedback controls in Section~\ref{sec:control-direct}. In this section, we will focus on open-loop controls and explain how the method can be implemented. 

In the open-loop setting, each player's control $\alpha^i$ is in the space $\mathbb{A} = \mathbb{H}^2_T(\mc{A})$. 
\begin{defn}[Open-loop Nash equilibrium]
    A strategy profile $(\hat{\bm{\alpha}}_t)_{t \in [0,T]} = (\hat \alpha^{1}_t, \ldots, \hat\alpha^{N}_t)_{t \in [0,T]} \in \mbA^N$ is called an open-loop Nash equilibrium if
\begin{equation}\label{def_Nash}
\forall i \in \mc{I} \text{ and } \beta^i \in \mbA, \quad J^i(\hat{\bm{\alpha}}) \leq J^i(\beta^i, \hat{\bm{\alpha}}^{-i}),
\end{equation}
where $\hat{\bm{\alpha}}^{-i}= [\hat \alpha^{1}, \ldots, \hat \alpha^{i-1}, \hat \alpha^{i+1}, \ldots, \hat \alpha^{N}] \in \mbA^{N-1}$ is the equilibrium strategies of all the players except the $i$-th one.
\end{defn}

Hu proposed in~\cite{Hu2:19} the idea of deep fictitious play (DFP for short), which consists of decoupling the $N$-player game into $N$ individual decision problems using \emph{fictitious play}, and then solving these $N$ individual problems using deep neural networks. Fictitious play was first introduced by Brown for static games \cite{Br:49,Br:51}, and was recently adapted to the mean-field setting by \cite{cardaliaguet2015learning,hadikhanloo2018ph,BrCa:18}. It is a simple yet important learning scheme in game theory for finding Nash equilibria. Deep learning provides efficient tools for solving the decoupled yet still high-dimensional optimization problems. 

The DFP algorithm starts with an initial strategy profile $(\balpha_t^0)_{t \in [0,T]} \in \mbA^N$, which can be interpreted as the initial belief of all players, and updates $\mt{K}$ times: $\balpha_t^0 \to \balpha_t^1 \to \ldots \to \balpha_t^\mt{k} \to \ldots \to \balpha_t^\mt{K}$ or until some stopping criterion is satisfied. At the beginning of stage $\mt{k}+1$, $\bm{\alpha}^\mt{k}$ is observable by all players. Player $i$ then considers that other players are going to reuse the same strategy as they used in the previous iteration, namely $\bm{\alpha}^{-i,\mt{k}}$, and looks for the best response. Then, player $i$ faces an optimization problem
\begin{equation}\label{def_J_SFP}
    \inf_{\beta^i \in \mbA} J^{i}(\beta^i;\bm{\alpha}^{-i, \mt{k}}), \quad 
    J^{i}(\beta^i; \bm\alpha^{-i, \mt{k}}) = \EE\left[\int_0^T f^i(t,\bm{X}_t^{\balpha}, (\beta^i, \bm{\alpha}^{-i, \mt{k}})) \ud t + g^i(\bm X_T^{\balpha})\right],
\end{equation}
where $\bm X_t^{\balpha} = [X_t^{1,{\balpha}}, X_t^{2,{\balpha}}, \ldots, X_t^{N,{\balpha}}]$ are state processes controlled by $(\beta^i, \bm\alpha^{-i,\mt{k}})$,
\begin{align}\label{eq:DPF-open-Xt}
\ud X_t^{\ell,\balpha}= b^\ell(t, \bm{X}_t^{\balpha}, (\beta^i, \bm{\alpha}^{-i, \mt{k}})) \ud t  +  \sigma^\ell(t, \bm{X}_t^{\balpha}, (\beta^i, \bm{\alpha}^{-i, \mt{k}})) \ud W_t^\ell   +  \sigma^0(t, \bm{X}_t^{\balpha}, (\beta^i, \bm{\alpha}^{-i, \mt{k}})) \ud W_t^0, \; 
\end{align}
with $X_0^{\ell,\alpha} = x^\ell$, for all $\ell \in \mc{I}$. Denote by $\alpha^{i, \mt{k}+1}$ the minimizer in \eqref{def_J_SFP},
\begin{equation}\label{def_SFP_opt-i}
    \alpha^{i,\mt{k} +1} = \argmin_{\beta^i \in \mbA^i} J^{i}(\beta^i; \bm{ \alpha}^{-i, \mt{k}}), \quad \forall i \in \mc{I}, \; \mt{k} < \mt{K}.
\end{equation}
We assume $\alpha^{i,\mt{k}+1}$ exists and is unique throughout this section. More precisely, $\alpha^{i,\mt{k}+1}$ is player $i$'s optimal strategy at stage $\mt{k}+1$ when the other players' controls are \eqref{def_Xt_general} evolve according to $\alpha^{j,\mt{k}}$, $j \neq i$. All players find their best responses simultaneously, which together form $\bm{\alpha}^{\mt{k}+1}$.

\begin{rem}
	Note that the above learning process is different than the usual simultaneous fictitious play \cite{Br:49,Br:51}, where the belief is described by the average over strategies played in previous iterations of the algorithm: $\frac{1}{\mt{k}}\sum_{\mt{k}'=1}^\mt{k} \bm{\alpha}^{-i,\mt{k}'}$. 
\end{rem}

In general, one can not expect that the players' actions $\balpha^\mt{k}$ always converge. However, if the sequence $\{\bm{\alpha}^{\mt{k}}\}_{\mt{k}=1}^\infty$ ever admits a limit, denoted by $\bm\alpha^{\infty}$, one would expect it to form an open-loop Nash equilibrium under suitable assumptions. Intuitively, in the limiting situation, when all other players are using strategies $\alpha^{j,\infty}$, $j \neq i$, by some stability argument, player $i$'s optimal strategy to the control problem \eqref{def_J_SFP} should be $\alpha^{i,\infty}$, meaning that she will not deviate from $\alpha^{i,\infty}$, which makes $(\alpha^{i,\infty})_{i=1}^N$ an open-loop equilibrium.

\begin{rem}[Theoretical analysis]
    In \cite[Theorem 3.1]{Hu2:19}, Hu proved that for linear-quadratic games, under appropriate conditions, the family $\{\bm \alpha^n\}_{n \in \NN}$  converges to an open-loop Nash equilibrium of the original problem \eqref{def_Xt_general}-\eqref{def:DFP-open-cost}. Moreover, the limit, denoted by $\bm\alpha^\infty$, is independent of the choice of initial belief $\bm\alpha^0$. 
\end{rem}

When seeking an open-loop equilibrium, admissible controls in the optimization problem~\eqref{def_SFP_opt-i} for player $i$ are  $\bW$-adapted. After a time discretization, the control at time $t_n$ can be expressed as a function of $\bW$ at the previous time steps. So, the problem is by nature high-dimensional, which motivates the use of deep neural networks. Similarly to the direct approach reviewed in Section~\ref{sec:control-direct}, Hu~\cite{Hu2:19} solved each player's control problems \eqref{def_J_SFP} at stage $\mt{k}$ by using a feedforward fully connected network to directly parameterize $\alpha^{i, \mt{k}}$. 
When seeking an open-loop equilibrium, the $N$ optimization problems~\eqref{def_SFP_opt-i} for $i=1,\dots,N$ can be solved in parallel, as $\alpha^{i, \mt{k}}$ are meant to be $\bW$-adapted thus are not influenced by changes in the other players' states $\bX_t^{-i}$. More precisely, to solve the optimization problem in~\eqref{def_SFP_opt-i}, each $\beta_{t_n}^{i, \mt{k}}$, $n=0,\dots,N_T-1$, is implemented by
\begin{equation}\label{def:NN-OL}
    \beta_{t_n}^{i, \mt{k}} \sim \beta_{t_n}^{i, \mt{k}}(\check\bX_0, \check\bW_{t_1}, \ldots, \check\bW_{t_n}; \theta_n^i),
\end{equation}
which maps $\RR^{N(n+1)}$ to $\RR$, 
and the optimal strategy at the $\mt{k}^{th}$ stage, denoted by $(\alpha_{t_n}^{i, \mt{k}}(\cdot; \theta_{n}^{i, \mt{k}}))_{n = 0}^{N_T-1}$, is determined by minimizing a discretized version of \eqref{def_J_SFP}, {\it i.e.},
\begin{equation}\label{eq:OL-cost}
   \{\theta_n^{i, \mt{k}}\}_{n=0}^{N_T-1} \in \argmin_{\{\theta_n^i\}_{n=0}^{N_T-1}} \EE\left[\sum_{n=0}^{N_T - 1} f^i(t_n, \check\bX_{t_n}^\theta, (\beta^{i, \mt{k}}_{t_n}(\check\bX_0, \check\bW_{t_1}, \ldots, \check\bW_{t_n}; \theta_n^i), \balpha_{t_n}^{-i, \mt{k}-1})) \Delta t + g^i(\check\bX_T^\theta)\right],
\end{equation}
where $(\check \bX_{t_n}^\theta)_{n=1}^{N_T}$ follows an Euler scheme corresponding to \eqref{eq:DPF-open-Xt} with controls $(\beta_{t_n}^{i,\mt{k}}, \balpha_{t_n}^{-i, \mt{k}-1})_{n=0}^{N_T-1}$ being used. The pseudo-code is given in Algorithm~\ref{algo:DFP-open} in Appendix~\ref{supp:pseudocode}.

\paragraph*{Numerical illustration: a linear-quadratic systemic risk problem with $24$ players.} \label{par:LQ-systemic-risk-finite-N} 

We illustrate the DFP algorithm on the LQ systemic risk game presented in Section~\ref{sec:intro-LQsysrisk} and originally introduced in \cite{CaFoSu:15}. For explicit formulas describing the open-loop Nash equilibrium, we refer to~\cite[Section~3.1]{CaFoSu:15}. 

In the numerical illustration, the number of players is set to be $N = 24$, and  the time steps is set at $N_T = 20$, after observing the maximum relative errors, defined by
\begin{equation}
    \max_{i \in \mathcal{I}} \frac{J^i(\alpha^{i, \mt{k}}; \balpha^{-i, \mt{k}-1}) - J^i(\hat \balpha)}{J^i(\hat \balpha)},
\end{equation}
did not increase too much from $N_T = 50$ to $N_T = 20$. The initial positions for the $i^{th}$ player is $x_0^i = 0.5i$, and results are presented in Figures~\ref{fig:N24_cost_traj}-\ref{fig:N24control}. Some key features that have been observed: the maximum of relative error drops below $3\%$ after ten iterations; the average error of estimated trajectories are convex/concave functions of time $t$; the standard deviation of estimated error aggregates from steps to steps. In fact, the convexity/concavity with respect to time $t$ is caused by two factors: the propagation of errors, which produces an increase in error mean; and the existence of terminal cost, which puts more weight on $X_T$ than $X_t, t \in (0,T)$, resulting in a better estimate of $X_T$ and a decreasing effect. 

\begin{figure}[h t p]
	\begin{tabular}{cc}
		\includegraphics[width=0.45\textwidth]{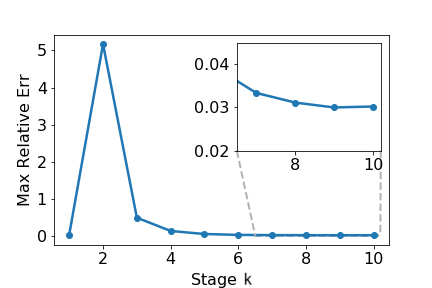}	&	
		\includegraphics[width=0.45\textwidth]{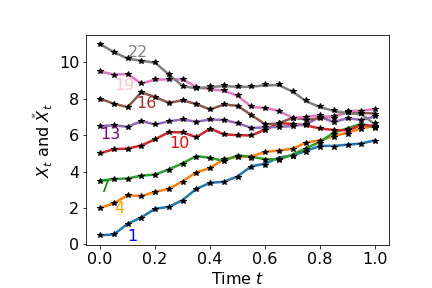}
	\end{tabular}
	\caption{Comparisons of cost functions and optimal trajectories for $N=24$ players in the linear quadratic systemic risk problem in Section~\ref{sec:open-loop}. Left: the maximum relative errors of the cost functions for 24 players; Right: for the sake of clarity, the comparison of optimal trajectories is only presented for the $1^\text{st}$, $4^\text{th}$, $7^\text{th}$, $10^\text{th}$, $13^\text{th}$, $16^\text{th}$, $19^\text{th}$ and $22^\text{th}$ players, where the solid lines are given by the closed-form solution and the stars are computed by deep fictitious play. %
 }\label{fig:N24_cost_traj}
\end{figure}

\begin{figure}[h t b]
	\begin{tabular}{ccc}
		\includegraphics[width=0.3\textwidth]{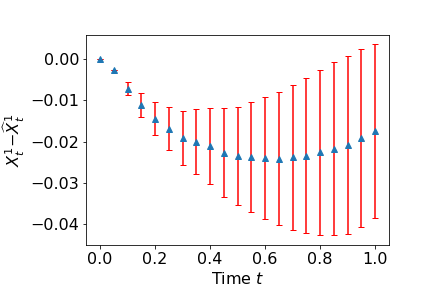} &
		\includegraphics[width=0.3\textwidth]{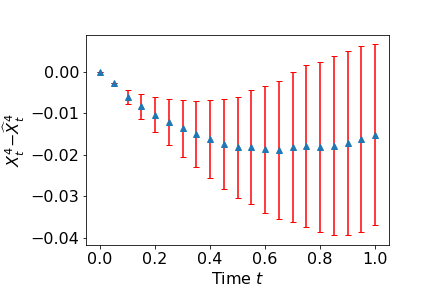}	&	
		\includegraphics[width=0.3\textwidth]{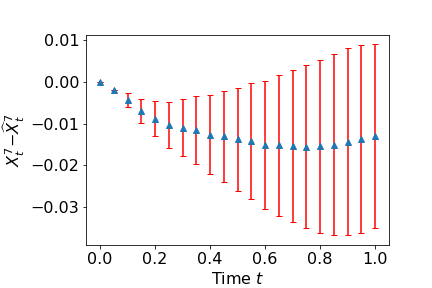} \\
		\includegraphics[width=0.3\textwidth]{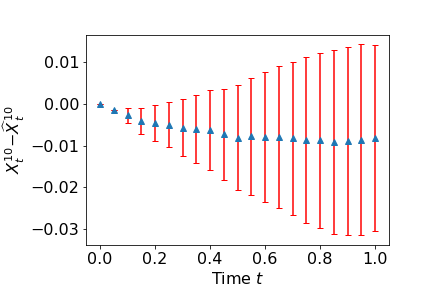} &
		\includegraphics[width=0.3\textwidth]{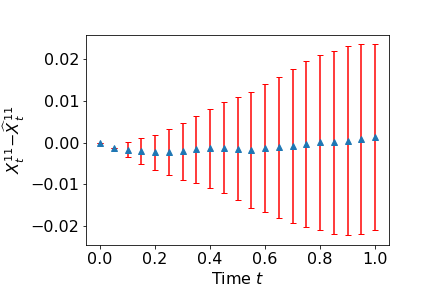} &
		\includegraphics[width=0.3\textwidth]{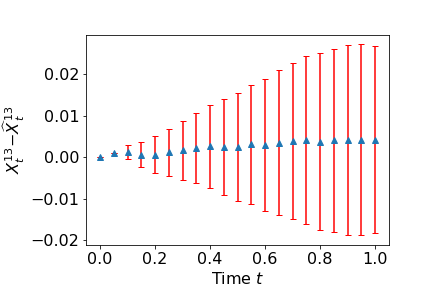} \\
		\includegraphics[width=0.3\textwidth]{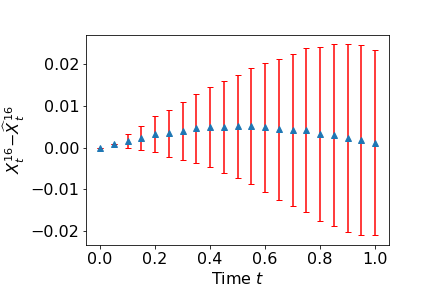} &
		\includegraphics[width=0.3\textwidth]{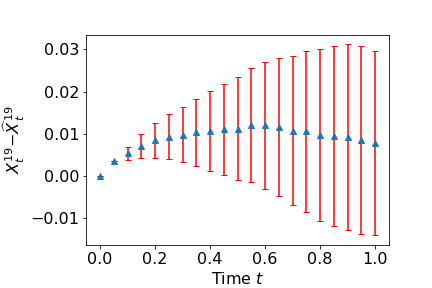} &
		\includegraphics[width=0.3\textwidth]{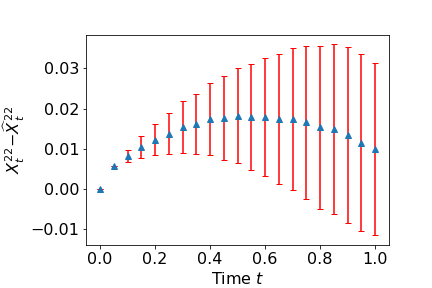} \\
	\end{tabular}
	\caption{Comparisons of trajectories for $N=24$ players in the linear quadratic game in Section~\ref{sec:open-loop} using the learnt equilibrium strategy profile. For the sake of clarity, we only show the mean (blue triangles) and standard deviation (red bars) of trajectories errors for the $1^\text{st}$, $4^\text{th}$, $7^\text{th}$, $10^\text{th}$, $11^\text{th}$, $13^\text{th}$, $16^\text{th}$, $19^\text{th}$ and $22^\text{th}$ player, respectively. The results are based on a total of $65536$ sample paths. They show that deep fictitious play provides a relatively uniformly good accuracy. %
 }\label{fig:N24traj}
\end{figure}

\begin{figure}[h t b]
	\begin{tabular}{ccc}
		\includegraphics[width=0.3\textwidth]{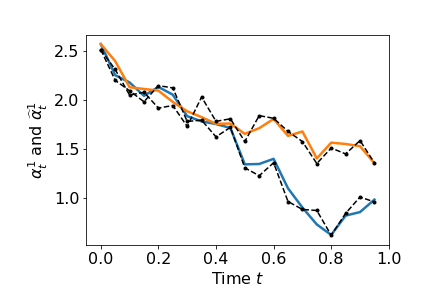} &
		\includegraphics[width=0.3\textwidth]{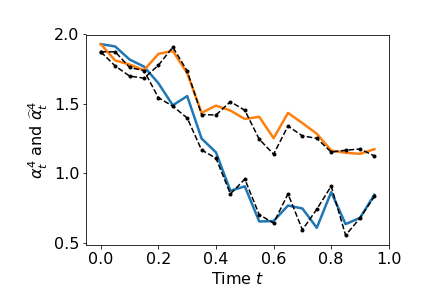}	&	
		\includegraphics[width=0.3\textwidth]{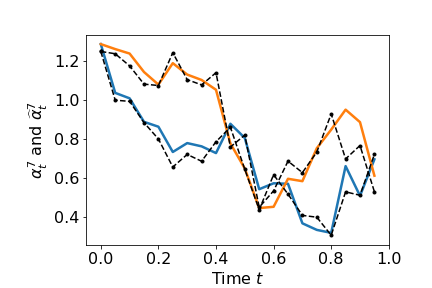} \\
		\includegraphics[width=0.3\textwidth]{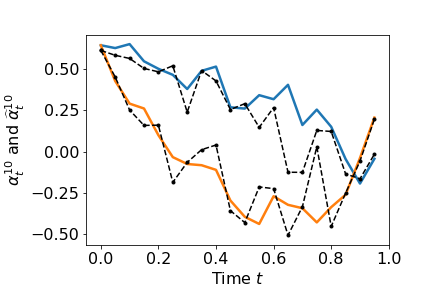} &
		\includegraphics[width=0.3\textwidth]{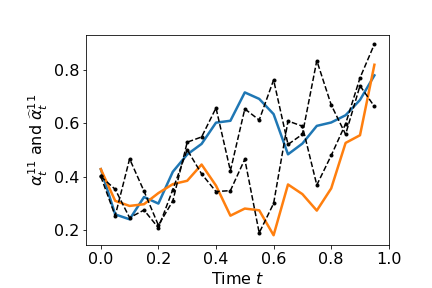} &
		\includegraphics[width=0.3\textwidth]{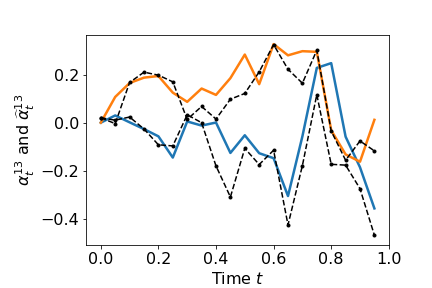} \\
		\includegraphics[width=0.3\textwidth]{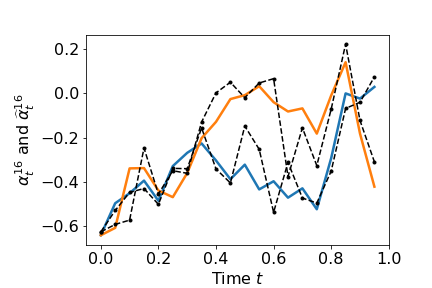} &
		\includegraphics[width=0.3\textwidth]{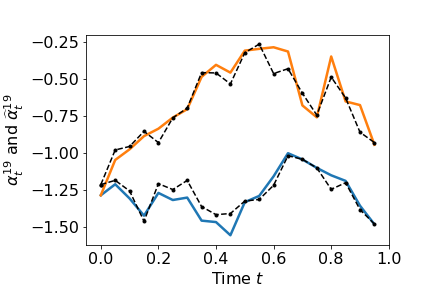} &
		\includegraphics[width=0.3\textwidth]{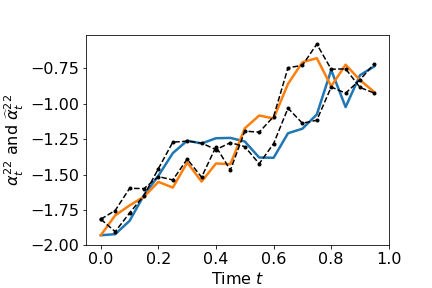} 
	\end{tabular}
	\caption{Comparisons of optimal controls for $N=24$ players in the linear quadratic game in Section~\ref{sec:open-loop}. For clarity, we only show two sample paths of optimal controls for the $1^\text{st}$, $4^\text{th}$, $7^\text{th}$, $10^\text{th}$, $11^\text{th}$, $13^\text{th}$, $16^\text{th}$, $19^\text{th}$ and $22^\text{th}$ player, respectively. The solid lines are optimal controls given by the closed-form solution, and the dotted dash lines are computed by deep fictitious play. %
 }\label{fig:N24control}
\end{figure}

In the original paper, the author also gave numerical experiments with $N = 5$ and $N = 10$ players. To better illustrate that the algorithm can mitigate the curse of dimensionality, the performance comparison is made across different $N$. In particular, the author computed the error $\max_{i \in \mc{I}}\max_{n \leq N_T}|X_{t_n}^i - \check X_{t_n}^{i, \theta}|$,
where $X$ denotes the state process following the open-loop Nash equilibrium, while $\check X^\theta$ is the deep fictitious play counterpart. The error is $1.09\times 10^{-2}$ for $N = 5$, $1.49\times 10^{-2}$ for $N = 10$ and $2.08\times 10^{-2}$ for $N = 24$.

\subsubsection{BSDE-based deep learning algorithms}
\label{sec:MarkovianNE}

We now turn our attention to the deep learning approach based on BSDE presented in Section~\ref{sec:BSDE} and discuss how it can be adapted to $N$-player games. Here again, various classes of information structures could be discussed, as long a system of BSDEs can be derived to characterize the Nash equilibrium. For the sake of the presentation, we will focus here on closed-loop feedback controls.

A Markovian strategy for player $i$ is a measurable function of $(t, \bX_t)$. Given a Markovian strategy profile  $\balpha_t$, \eqref{def_Xt_general} rewrites as
\begin{equation}
    \ud X_t^i = b^i(t, \bm{X}_t, \bm{\alpha}(t, \bX_t)) \ud t + \sigma^i(t, \bm{X}_t, \bm{\alpha}(t, \bX_t)) \ud W_t^i + \sigma^0(t, \bm{X}_t, \bm{\alpha}(t, \bX_t)) \ud W_t^0, \quad X_0^i = x^i_0, \quad i \in \mc{I}.
\end{equation}
Each player aims to minimize the cost 
\begin{equation}\label{def:DFP-closed-cost}
J^i(\bm{\alpha}) = \EE\left[\int_0^T f^i(t, \bm X_t, \bm\alpha(t, \bm X_t)) \ud t + g^i(\bm X_T)\right].
\end{equation}
The dependence of $b^i$ and $\sigma^i$ on $\balpha$ leads to a problem with stronger coupling than the open-loop setting. The problem is thus harder to solve, both theoretically and numerically. 

\begin{defi}
A strategy profile $\hat{\bm \alpha} = (\hat \alpha^{1}, \ldots, \hat \alpha^{N})$ consisting of Markovian strategies is called a Markovian Nash equilibrium, if 
\begin{equation}
\forall i \in \mc{I}, \emph{ and any Markovian strategy $\beta^i$} , \quad J^i(\hat {\bm \alpha}) \leq J^i(\hat \alpha^{1}, \ldots,\hat \alpha^{i-1}, \beta^i, \hat \alpha^{i+1}, \ldots, \hat \alpha^{N}),
\end{equation}
where on the right-hand side, $(\hat \alpha^1, \ldots, \hat \alpha^{i-1}, \beta^i, \hat \alpha^{i+1}, \ldots, \hat \alpha^N)(t, \bX_t)$ is used in \eqref{def_Xt_general} to solve for $\bX_t$. 
\end{defi}

In this case, due to the fact that each player's control is a function of the other players' states, any change in $\bX_t$ will cause a change in $\balpha_t$. Consequently, merely adapting the algorithm in Section~\ref{sec:open-loop} will not be efficient nor parallelizable. In the Markovian setting, finding a Nash equilibrium can be reduced to solving $N$ coupled HJB equations. To this end, we define $u^i(t, \bm x)$ as the value function of player $i$.  
For the sake of notation clarity, we present the discussion based on the ``vectorized'' system,
\begin{equation}\label{def_Xt_markovian}
\ud \bm X_t = b(t, \bm X_t, \bm \alpha(t, \bm X_t)) \ud t + \Sigma(t, \bm X_t, \balpha(t, \bX_t)) \ud \bm W_t, \quad \bm X_0 = \bm x_0,
\end{equation}
where $\bX$, $b(t, \bx, \balpha(t, \bx))$ and $\bW$ are vectorizations of $X^i, b^i, W^i$ respectively, and $\Sigma(t, \bx, \balpha(t, \bx))$ is matrix-valued given by
\begin{equation}
   \bX_t = \begin{bmatrix}
  X_t^1 \\
  X_t^2\\
  \vdots\\
  X_t^N
\end{bmatrix}, \quad
b = \begin{bmatrix}
     b^1\\
     b^2\\
     \vdots\\
     b^N
\end{bmatrix}, \quad 
\bW_t =\begin{bmatrix}
     W_t^0\\
     W_t^1 \\
     \vdots\\
     W_t^N
\end{bmatrix}, \quad
\Sigma = \begin{bmatrix}
     \sigma^0 & \sigma^1\\
     \sigma^0 & \sigma^2 \\
     \vdots & \vdots \\
     \sigma^0 & \sigma^N
\end{bmatrix}.
\end{equation}
Using the dynamic programming principle, the HJB system reads 
\begin{align}
\label{def_HJB}
\begin{dcases}
    \partial_t u^i(t,\bx) + G^i(t,\bm x,\hat{\bm \alpha}(t,\bx), \nabla_{\bx} u^i(t,\bx), \text{Hess}_{\bx} u^i(t,\bx))  = 0,\\
u^i(T,\bx) = g^i(\bx), \quad i \in \mc{I},
    \\
    \hat\alpha^i(t,\bx) = \inf_{\alpha^i \in \mc{A}} G^i(t,\bm x, (\alpha^i, \hat{\bm \alpha}^{-i}(t,\bx)), \nabla_{\bx} u^i(t,\bx), \text{Hess}_{\bx} u^i(t,\bx)),
\end{dcases}
\end{align}
where $G^i$ is given by 
\begin{equation}\label{def:Gi}
    G^i = G^i(t,\bm x, \bm \alpha, p, q) =  b(t,\bm x, \bm \alpha) \cdot p + f^i(t, \bm x, \bm \alpha) + \half \text{Tr}(\Sigma\Sigma\transpose(t, \bx, \balpha) q).
\end{equation}
Note that when $N = 1$, $G^i$ coincides with the $H$ function defined in \eqref{def:control-H}. When $N>1$, the equation for $u^i$ in \eqref{def_HJB} depends on the equilibrium controls of the other players, which themselves depend on the value functions of these players. Hence, the $N$ equations are coupled.

Han and Hu proposed in~\cite{HaHu:19} a version of the DFP algorithm for finding Markovian NE. Similar to ideas reviewed in Section~\ref{sec:open-loop}, the optimization problems for the $N$ players are decoupled by fictitious play in the following sense. Starting with a guess of the solution $\bm\alpha^0 \in \mathbb{A}$, the algorithm iteratively updates $\bm\alpha^{\mt{k}+1}$: at iteration $\mt{k}+1$, $\bm \alpha^\mt{k}$ is observed by all players, and player $i$'s decision problem is 
\begin{equation}\label{def_J_fictitious}
\inf_{\alpha^i \in \mathbb{A}^i} J^i(\alpha^i; \bm \alpha^{-i,\mt{k}}),
\end{equation}
where $J^i$ is defined in \eqref{def:DFP-open-cost}, and the state process $\bm X_t$ is given in \eqref{def_Xt_markovian} with $\bm \alpha$ replaced by $(\alpha^i, \bm \alpha^{-i,\mt{k}})$. This problem is decoupled from the problems solved by the other players at the same iteration. The optimal strategy, assuming it exists, is denoted by $\alpha^{i, \mt{k}+1}$.  
Problems \eqref{def_J_fictitious} for all $i \in \mc{I}$ are solved simultaneously using $\bm \alpha^{-i, \mt{k}}$, and the optimal responses together form $\bm \alpha^{\mt{k}+1}$. Thanks to the Markovian structure, problem \eqref{def_J_fictitious} corresponds to solving the HJB equation 
\begin{equation}
\label{eq:FP_PDE}
    \partial_t u^{i,\mt{k}+1} + \inf_{\alpha^i \in \mc{A}^i} G^i(t, \bm x, (\alpha^i, \bm \alpha^{-i,\mt{k}}(t, \bm x)), \nabla_{\bx} u^{i, \mt{k}+1}, \text{Hess}_{\bx}u^{i, \mt{k}+1}, u^{i, \mt{k}+1}) = 0,
\end{equation}
with the terminal condition $u^{i,\mt{k}+1}(T,\bx) = g^i(\bx)$. In the uncontrolled volatility problem (as {\it e.g.} in~\cite{HaHu:19}), $G^i$ does not depend on the $q$ variable and hence the PDE is semi-linear. Let us assume that the minimizer $\argmin_{\alpha^i \in \mc{A}^i} G^i(t, \bm x, \bm \alpha, p, s)$ exists, is unique, and can be computed as an explicit function of other arguments:  $(t, \bm x, (\alpha^j)_{j \neq i}, p, s) \mapsto \alpha^{i,\mt{k}+1,*}(t, \bm x, (\alpha^j)_{j \neq i}, p, s)$. Then, the new control for player $i$ is given by: 
\begin{equation}\label{def_alphaast}
    \alpha^{i, \mt{k}+1}(t, \bm x) = \alpha^{i,\mt{k}+1,*}(t, \bm x, \bm \alpha^{-i,\mt{k}}(t, \bm x), \nabla_{\bx} u^{i, \mt{k}+1}(t, \bx), u^{i, \mt{k}+1}(t, \bx)).
\end{equation} 
Solving~\eqref{eq:FP_PDE} for all $i\in \mc{I}$ completes one stage in the loop of fictitious play. 

To solve the $N$ decoupled semi-linear PDEs numerically, one can for example interpret their solution through BSDEs via the non-linear Feynman-Kac formula~\eqref{def_control_BSDEreform}:
\begin{empheq}[left=\empheqlbrace]{align}
    &\bX_t^{i, \mt{k}+1} = \bx_0 + \int_{0}^{t}\tilde b^{i}(s,\bX_s^{i, \mt{k}+1}; \balpha^{-i, \mt{k}}(s, \bX_s^{i, \mt{k}+1}))\, \mathrm{d}s + \int_{0}^{t}\Sigma(s,\bX_s^{i, \mt{k}+1})\, \mathrm{d}\bW_s, \label{eq:BSDE_forward} \\
    &Y_t^{i, \mt{k}+1} = g^i(\bX_T^{i, \mt{k}+1}) \nonumber \\ &\hspace{5em} + \int_{t}^{T}h^i(s,\bX_s^{i, \mt{k}+1}, Y_s^{i, \mt{k}+1}, \bZ_s^{i, \mt{k}+1};  \balpha^{-i, \mt{k}}(s, \bX_s^{i, \mt{k}+1}))\, \mathrm{d}s - \int_{t}^{T}(\bZ_s^{i, \mt{k}+1})\transpose\, \mathrm{d}\bW_s, \label{eq:BSDE_backward} 
\end{empheq}
with $ Y_t^{i, \mt{k}+1} = u^{i, \mt{k}+1}(t, \bX_t^{i, \mt{k}+1})$ and $\bZ_t^{i, \mt{k}+1} = \Sigma(t, \bX_t^{i, \mt{k}+1})\transpose\nabla_{\bx} u^{i, \mt{k}+1}(t, \bX_t^{i, \mt{k}+1})$, where $\tilde b^i$ and $h^i$ are functions such that~\eqref{eq:FP_PDE} can be rewritten as
\begin{align}
     \partial_t u^{i,\mt{k}+1} + \half \text{Tr}(\Sigma\transpose \text{Hess}_{\bx} u^{i,\mt{k}+1} \Sigma) &+ \tilde b^i(t, \bx; \balpha^{-i, \mt{k}})\cdot \nabla_{\bx} u^{i, \mt{k}+1} 
     \notag 
     \\ 
     &+h^i(t, \bx,  u^{i, \mt{k}+1}, \Sigma\transpose\nabla_{\bx} u^{i, \mt{k}+1}; \balpha^{-i, \mt{k}})=0. \label{eq:FP_PDE_explicit}
\end{align}
We stress once again that $\balpha^{-i,\mt{k}}$ are known functions when solving the BSDE,  which is why the BSDEs for $i \in \mc{I}$ are decoupled. They can thus be solved in parallel, for example using the deep BSDE algorithm (see Section~\ref{subsubsec:deepBSDE}).

For each player $i$, the BSDE \eqref{eq:BSDE_forward}--\eqref{eq:BSDE_backward} is then solved using the deep BSDE algorithm (see Section~\ref{subsubsec:deepBSDE}). More precisely, 
let us consider the minimization problem (for less cumbersome notations, the subscript $t_n$ in $\bX, \bY, \bZ$ has been replaced by $n$, and the superscript $\mt{k}$ that denotes the index of fictitious play is dropped for simplicity),  
\begin{align}
&\inf_{\psi_0\in \cN_0^{i'},~\{\phi_n\in \cN_n^i\}_{n=0}^{N_T-1} } \EE|g^i(\check\bX_{N_T}^{i}) - \check \bY_{N_T}^{i}|^2, \label{eq:disc_objective}\\
&s.t.~~ \check\bX_0^{i}=\bx_0, \quad \check \bY_0^{i} =\psi_0(\check \bX_0^{i}), \quad \check \bZ_{n}^{i}=\phi_n(\check \bX_{n}^{i}), \quad n=0,\dots,N_T-1\notag \\
&\qquad \check\bX_{n+1}^{i} = \check \bX_{n}^{i} + \tilde b^i(t_n,\check\bX_{n}^{i}; \balpha^{-i}(t_n, \check\bX_n^{i}))\Delta t +\Sigma(t_n,\check\bX_{n}^{i})\Delta \bW_{t_n}, \label{eq:disc_X_path} \\
 &\qquad \check \bY_{n+1}^{i} = \check \bY_{n}^{i} - h^i(t_n,\check \bX_{n}^{i},\check \bY_{n}^{i},\check \bZ_{n}^{i}; \balpha^{-i}(t_n, \check\bX_n^{i}))\Delta t + (\check\bZ_{n}^{i})\transpose\Delta \bW_{t_n}, \label{eq:disc_Y_path}
\end{align}
where we recall that
  $\Delta t = t_{n+1} - t_n, \quad \Delta \bW_{t_n} = \bW_{ t_{n+1}}-\bW_{t_n}.$
Here $\cN_0^{i'}$ and $\{\cN_n^i\}_{n=0}^{N_T-1}$ are hypothesis spaces of player $i$ related to deep neural networks. The goal of the optimization is to find optimal deterministic maps $\psi_0^{i,\ast}, \{\phi_n^{i,\ast}\}_{n=0}^{N_T-1}$ such that the loss function is small. The pseudo-code of the proposed deep fictitious play algorithm is summarized in Algorithm \ref{def_algorithm1} in Appendix~\ref{supp:pseudocode}.

\begin{rem}
In additional to the cost functional considered in \eqref{def:DFP-open-cost}, in \cite{HaHu:19} Han and Hu also considered the risk-sensitive minimization problem,
\begin{equation}\label{def_J'}
J^i(\bm \alpha) = \EE\left[\xi_i \exp\left\{\xi_i\left(\int_0^T f^i(t, \bm X_t, \bm \alpha(t, \bm X_s)) \ud t + g^i(\bm X_T)\right)\right\} \right],
\end{equation}
where $\xi_i$ is a parameter characterizing how risk-averse or risk-seeking player $i$ is. This flexibility allows one to model much broader classes of games that accommodate the players' attitudes to risk. In this case, the $G^i$ function defined in \eqref{def:Gi} becomes
\begin{equation}
     G^i = G^i(t,\bm x, \bm \alpha, p, q, s) =  b(t,\bm x, \bm \alpha) \cdot p + \xi_i sf^i(t, \bm x, \bm \alpha) + \half \text{Tr}(\Sigma\Sigma\transpose (t, \bx, \balpha) q ).  
\end{equation}
Numerical examples on risk-sensitive problems are also presented therein. 
\end{rem}

\begin{rem}[Theoretical analysis]
In \cite{han2020convergence}, Han, Hu and Long provided a theoretical foundation for the DFP algorithm for Markovian Nash equilibria with the objective~\eqref{def:DFP-open-cost}. They proved the convergence to the true Nash equilibrium if the decoupled sub-problems \eqref{def_J_fictitious}--\eqref{eq:FP_PDE} are solved exactly and repeatedly \cite[Theorem 2]{han2020convergence}. They also gave a posteriori error bound on the numerical error on the deep BSDE algorithm \cite[Theorem 3]{han2020convergence},  identified the $\eps$-Nash equilibrium produced by DFP \cite[Theorem 4]{han2020convergence}, and analyzed the numerical performance of the algorithm on the original game \cite[Theorem 5]{han2020convergence}. 
\end{rem}

The DFP algorithm has then been extended in \cite{chen2021large} with the Scaled Deep Fictitious Play (SDFP) algorithm. There, the authors integrated the importance sampling and invariant layer embedding into DFP. Then focusing on the homogeneous agent problem, they utilized the symmetry and showed numerical experiments with up to $N = 3,000$ agents. Along with the idea of combining fictitious play and deep learning, \cite{han2021deepham} recently studied the heterogeneous agent model in macroeconomics. In addition to parameterizing the control (as proposed in \cite{han2016deep-googlecitations}), they also use neural networks to approximate the value function to reduce further the computational complexity of evaluating expectations of the type \eqref{eq:OL-cost}. \cite{EliePerolatLauriereGeistPietquin-2019_AFP-MFG} solved the mean-field games in model-free setting and analyzed its convergence by the error analysis of
each learning iteration step, in analogy with how the convergence of RL algorithms reduces to the aggregation of
repeated supervised learning approximation errors \cite{farahmand2010error,scherrer2015approximate}. \cite{cont2022dynamics} studied the interactions among market makers in dealer markets via solving a stochastic differential game of intensity control with partial information, and provided insights into tacit collusion due to market price interactions.

\paragraph*{Numerical illustration: the linear-quadratic systemic risk example revisited.}\label{sec:LQsysrisk-revisited} 
Here we revisit the example introduced in Section~\ref{sec:intro-LQsysrisk} and studied in Section~\ref{sec:open-loop}, but focus on the Markovian Nash equilibrium.

To describe the model in the form of \eqref{def_Xt_markovian}, we concatenate the log-monetary reserves $X_t^i$ of $N$ banks to form $\bX_t=[X_t^1,\dots,X_t^N]\transpose$.
The associated drift term and diffusion term are defined as
\begin{equation}
\label{eq:example1_driftdiff}
b(t, \bm x, \balpha)=
[a(\bar x - x^1) + \alpha^1, 
\ldots, 
a(\bar x - x^N) + \alpha^N
]\transpose\in \RR^{N\times1}
, \quad \bar x = \frac{1}{N}\sum_{i=1}^N x^i, \end{equation}
\begin{equation}
\label{eq:example1_sigma}
\Sigma(t, \bm x)=
\begin{bmatrix}
     \sigma\rho & \sigma\sqrt{1-\rho^2} & 0 & \cdots & 0\\
     \sigma\rho & 0 & \sigma\sqrt{1-\rho^2} & \cdots & 0\\
     \vdots & \vdots & \vdots & \ddots & \vdots\\
     \sigma\rho & 0 & 0 & \cdots &\sigma\sqrt{1-\rho^2}
\end{bmatrix}\in \RR^{N\times (N+1)},
\end{equation}
and $\bW_t = (W_t^0, \ldots, W_t^N)$ is $(N+1)$-dimensional. Recall the running and terminal costs that player $i$ aims to minimize
\begin{equation}
f^i(t, \bm{x},\balpha) = \half (\alpha^i)^2 - q \alpha^i(\bar x - x^i) + \frac{\eps}{2}(\bar x - x^i)^2,  \quad g^i(\bm{x}) = \frac{c}{2}(\bar x - x^i)^2.
\end{equation} 
The closed-form Nash equilibrium is detailed in \cite[Sections 3.2-3.3]{CaFoSu:15}.

The coupled HJB system corresponding to this game reads
\begin{equation}\label{eq_HJBexample1}
\partial_t u^i + \inf_{\alpha^i}\left\{\sum_{j=1}^N [a(\bar x - x^j) + \alpha^j] \partial_{x^j}u^i + \frac{(\alpha^i)^2}{2} - q\alpha^i(\bar x - x^i) + \frac{\eps}{2}(\bar x - x^i)^2 \right\} 
+ \half \text{Tr}(\Sigma\transpose \text{Hess}_{\bx} u^i \Sigma) = 0,
\end{equation} 
with the terminal condition $u^i(T,\bm x) = \frac{c}{2}(\bar x - x^i)^2$, $i \in \mc{I}$. The minimizer in the infimum gives a candidate of the optimal control for player $i$: 
$\alpha^{i}(t, \bm x) = q(\bar x - x^i) - \partial_{x^i}u^i(t, \bm x).$  
Plugging it back into the $i^{th}$ equation yields a PDE of form \eqref{eq:FP_PDE_explicit},
\begin{multline}\label{eq:example1}
    \partial_t u^i + \half \text{Tr}(\Sigma\transpose \text{Hess}_{\bx} u^i \Sigma) + a(\bar x - x^i)\partial_{x^i}u^i + \sum_{j\neq i} [a(\bar x - x^j) + \alpha^{j}(t,\bm x)] \partial_{x^j}u^i \\
 + \frac{\eps}{2}(\bar x - x^i)^2 - \half(q(\bar x - x^i) - \partial_{x^i} u^i)^2= 0,
\end{multline}
where $\alpha^j$ with $j\neq i$ are considered exogenous for player $i$'s problem, and are given by the best responses of the other players from the previous stage. To be precise, $\tilde b^i$ and $h^i$ in \eqref{eq:FP_PDE_explicit} are defined as
\begin{align}
\tilde b^i(t, \bx; \balpha^{-i})&=
[
a(\bar{x} - x^1) + \alpha^1, 
\ldots,
a(\bar{x} - x^i),
\ldots,
a(\bar{x} - x^N) + \alpha^N
]\transpose,\\
h^i(t, \bx,  y, \bz;  \balpha^{-i})&=\frac{\eps}{2}(\bar x - x^i)^2 - \half(q(\bar x - x^i) - \frac{z^i}{\sigma\sqrt{1-\rho^2}})^2,
\end{align}
where $\bz=(z^0, z^1,\dots, z^N)\in\RR^{N+1}$.

Figures~\ref{fig:LQMFG_err}--\ref{fig:LQMFG_path} show the performance of the DFP algorithm on a ten-player game, using the parameter,
\begin{equation}\label{def_parameters}
    a = 0.1,\quad  q = 0.1, \quad c = 0.5,\quad  \eps = 0.5,\quad  \rho = 0.2, \quad \sigma = 1, \quad T = 1.
\end{equation}
The relative squared error (RSE) is defined  by 
{\footnotesize
\begin{align*}
    \text{RSE} = \frac{\sum_{\substack{i \in \mathcal{I}\\1 \leq  j \leq J}} \left(u^i(0, \bm x_{t_0}^{(j)}) - \widehat u^i(0,\bm x_{t_0}^{(j)})\right)^2}{\sum_{\substack{i \in \mathcal{I}\\ 1 \leq j \leq J}} \left(u^i(0, \bm x_{t_0}^{(j)}) - \bar u^i\right)^2},
    \;\text{or }\; 
    \text{RSE} = \frac{\sum_{\substack{i \in \mathcal{I}\\ 0 \leq n \leq N_T-1 \\ 1 \leq j \leq J}} \left(\nabla_{\bx} u^i(t_n, \bm x_{t_n}^{(j)}) - \nabla_{\bx} \widehat u^i(t_n,\bm x_{t_n}^{(j)})\right)^2}{\sum_{\substack{i \in \mathcal{I}\\ 0 \leq n \leq N_T-1\\1 \leq j \leq J}} \left(\nabla_{\bx} u^i(t_n, \bm x_{t_n}^{(j)})- \overline {\nabla_{\bx} u}^i\right)^2},
\end{align*}
}where  $\hat{u}^i$ is the prediction from the neural networks,
and $\bar u^i$ (resp. $\overline {\nabla_{\bx} u}^i$) is the average of $u^i ~(\emph{resp.~} \nabla_{\bx}u^i)$ evaluated at all the indices $j, n$. 
To compute the relative error, $J=256$ ground truth sample paths $\{\bm x_{t_n}^{(j)}\}_{n=0}^{N_T-1}$ are generated using an Euler scheme based on  \eqref{def_Xt_markovian}\eqref{eq:example1_driftdiff}\eqref{eq:example1_sigma} and the true optimal strategy. Note that the superscript ${(j)}$ here does not mean the player index, but the $j^{th}$ path for all players. 

In particular, Figure~\ref{fig:LQMFG_err} compares the relative squared error
as $N_{\text{SGD\_per\_stage}}$ varies from 10 to 400. The convergence of the learning curves with small $N_{\text{SGD\_per\_stage}}$ asserts that each individual problem does not need to be solved so accurately. Furthermore, the fact that the performances are similar under different $N_{\text{SGD\_per\_stage}}$ with the same total budget of SGD updates suggests that the algorithm is insensitive to the choice of this hyperparameter. %
The final relative squared errors of $u$ and $\nabla u$ averaged from three independent runs of deep fictitious play are 4.6\% and 0.2\%,  respectively. 
Figure~\ref{fig:LQMFG_path} presents one sample path for each player of the optimal state process $X_t^i$ and the optimal control $\alpha_t^i$ \emph{vs.} their approximations $\hat{X}_t^i, \hat{\alpha}_t^i$ provided by the optimized neural networks.

\begin{figure}[!ht]
    \centering
    \includegraphics[width=\figwidth]{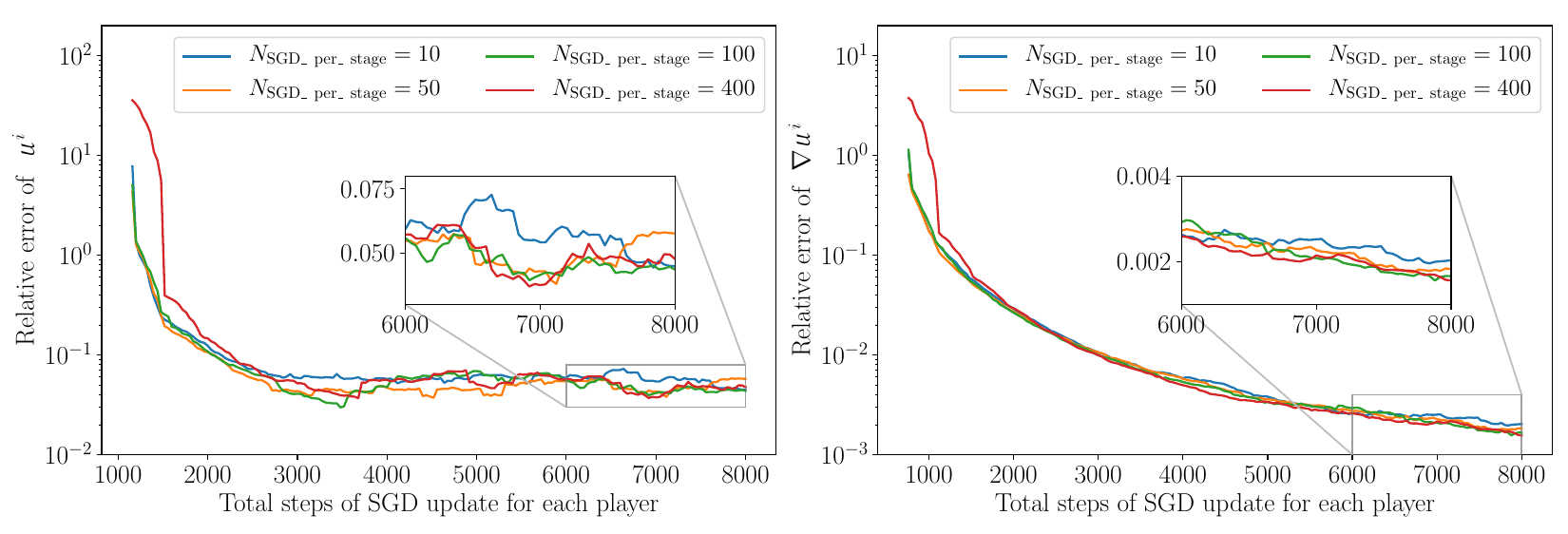}
    \caption{Linear-quadratic systemic risk example in Section~\ref{sec:MarkovianNE}. The relative squared errors of $u^i$ (left) and $\nabla u^i$ (right) along the training process of deep fictitious play for the inter-bank game. 
    The relative squared errors of $u^i(0, \check \bX_0^{i})$ and $\{\nabla u^i(t_n, \check \bX_n^{i})\}_{n=0}^{N_{T-1}}$ are evaluated. 
    }
    \label{fig:LQMFG_err}
\end{figure}

\begin{figure}[!ht]
    \centering
    \includegraphics[width=\figwidth]{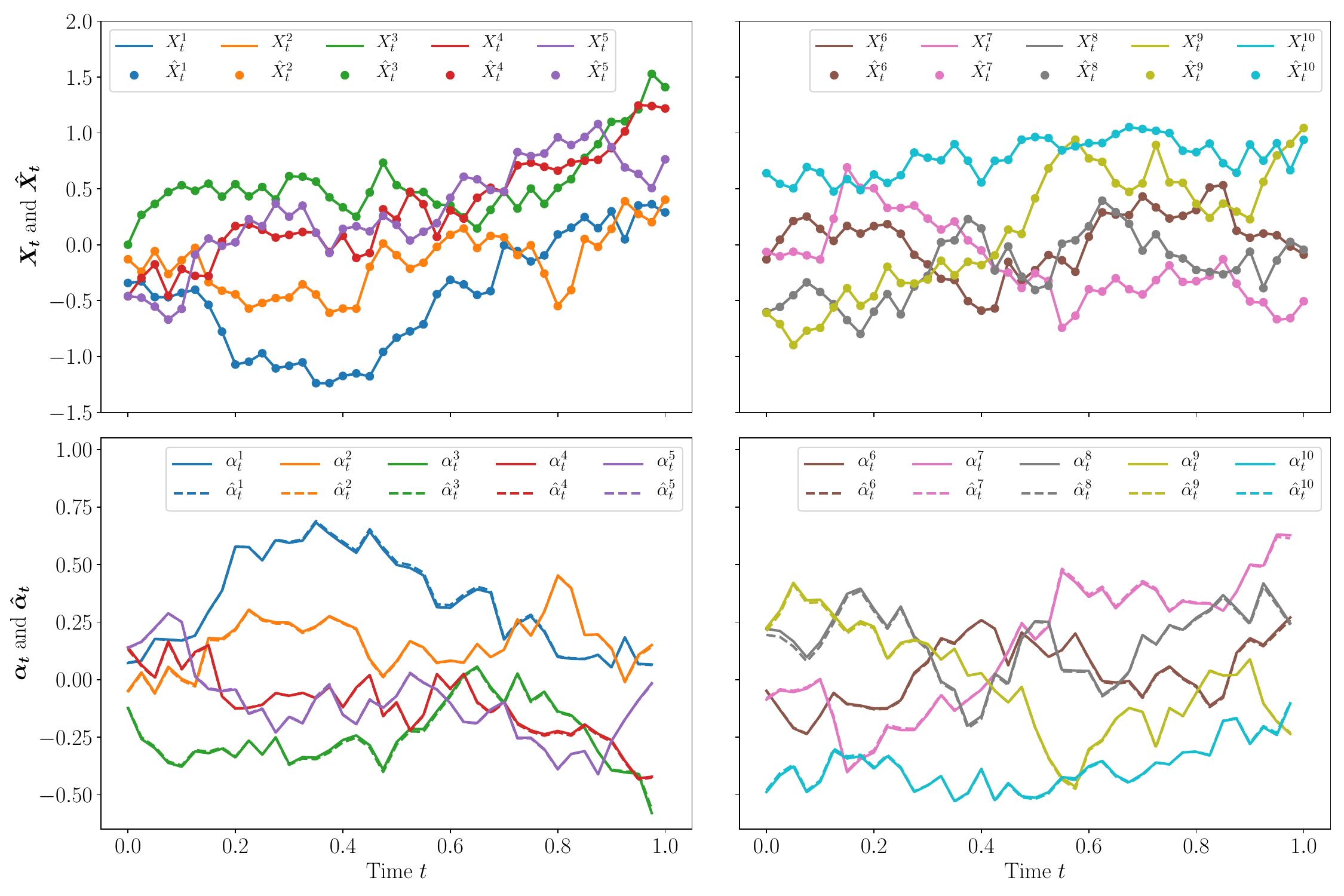}
    \caption{Linear-quadratic systemic risk example in Section~\ref{sec:MarkovianNE}. A sample path for each player of the inter-bank game with $N=10$. Top: the optimal state process $X_t^i$ (solid lines) and its approximation $\hat{X}_t^i$ (circles) provided by the optimized neural networks, under the same realized path of Brownian motion. Bottom: comparisons of the strategies $\alpha_t^i$ and $\hat{\alpha}_t^i$ (dashed lines). %
    } 
    \label{fig:LQMFG_path}
\end{figure}

\subsubsection{PDE-based deep learning algorithms}

For the sake of completeness, we briefly discuss how PDE-based methods can be adapted to solve $N$-player Nash equilibria. As mentioned earlier, such equilibria can be characterized using PDE systems, as seen in~\eqref{def_HJB}. We recall that this system stems from the combination of the HJB equations for every player's control problem. They are coupled since the value function of each player depends on the other players' actions. Instead of translating this system into an FBSDE system as discussed in Section~\ref{sec:MarkovianNE}, one can directly solve the HJB equations using for instance the DGM method presented in Section~\ref{sec:PDE} for stochastic optimal control problems. The value function for each player is replaced by a neural network. Then, the loss is computed based on the residuals of all the PDEs at random points in space and time. The application is rather straightfoward so we refrain from repeating the DGM algorithm. Technical details in the implementation (\textit{e.g.}, the neural network architecture or the distribution used to sample points) should be discussed on a case-by-case basis.

This approach has been used for instance in~\cite[Section 5]{al2022extensions} by Al-Aradi {\it et al.} to solve a system of HJB equations arising in a finite-player game modeling systemic risk introduced by Carmona {\it et al.}~\cite{MR3325083}. This example has already been discussed in the numerical illustration provided in Section~\ref{sec:open-loop}.

\subsection{Mean-field games}\label{sec:MFG}

Mean-field games, introduced independently by Lasry and Lions in \cite{LaLi1:2006,LaLi2:2006,LaLi:2007} and by Huang, Malham\'{e} and Caines in \cite{HuMaCa:06,HuCaMa:07}, provide a paradigm to approximate the solutions of stochastic games with very large number of players. The approximation relies on two main assumptions: anonymity, which means that the players interact only through the population's empirical distribution, and indistinguishability, which means that $(b^i, \sigma^i, f^i, g^i)$ are the same for all $i$. Then, one can pass (at least formally) to the limit by letting the number of players grow to infinity. In the asymptotic problem, the influence of each individual player on the rest of the population vanishes and the Nash equilibrium can be characterized by studying the problem posed to a representative player. Under suitable assumptions, it can be shown that solving this limiting problem provides an approximate equilibrium for the finite-player game.  
We refer to the notes~\cite{Cardaliaguet-2013-notes} and the books \cite{MR3134900,carmona2018probabilistic,carmona2018probabilistic2} as well as the references therein for further background on mean-field games.  Mean field games are usually studied using the notion of mean field Nash equilibrium. As in finite-player games, intuitively this notion corresponds to a situation in which no player can benefit from unilateral deviation. This is different from the solution to MFC problems discussed in Section~\ref{sec:directMethod}, which can be interpreted as a social optimum, \textit{i.e.}, a situation in which the players collectively choose their control in order to minimize the social cost.

Numerical methods for such games have been developed using mostly traditional techniques such as finite difference schemes~\cite{achdoucapuzzo2010mean,achdoucamillicapuzzo2012mean}, semi-Lagrangian schemes~\cite{carlinisilva2014fully,carlinisilva2018discretization}, or methods based on probabilistic approaches~\cite{chassagneuxcrisandelarue2019numerical,angiuli2019cemracs}. See, {\it e.g.}, \cite{achdoulauriere2020meansurvey,AMSnotesLauriere} for recent surveys. However, similarly to control problems or finite-player games, these methods do not scale well in terms of dimensionality, and in particular, they are not very suitable for problems with delay or with common noise. This motivates the development of deep learning methods, some of which are described in~\cite{carmonalauriere2021deepmfgsurvey}. In the sequel, we describe the theoretical framework of mean field games and then survey recent deep learning methods.

\subsubsection{Theoretical background}
\label{sec:mfg-theoretical-background}
We start by defining the notion of MFG, and then we discuss how equilibria can be characterized in terms of PDEs, BSDEs, and the so-called master equation. 

\paragraph{Definition of the problem}
Going back to problem~\eqref{def:DFP-open-cost}, let us assume that $b^i, \sigma^0, \sigma^i, f^i, g^i$ depend on the rest of the population's states and actions in an anonymous way, {\it i.e.}, there exists functions $b, \sigma^0, \sigma, f, g$ such that 
\begin{equation}
\begin{aligned}
&b^i(t, \bm{X}_t, \bm{\alpha}_t) = b(t, X^i_t, \nu^N_t, \alpha^i_t), \quad \sigma^0(t, \bm{X}_t, \bm{\alpha}_t) = \sigma^0(t, X^i_t, \nu^N_t, \alpha^i_t), \quad \sigma^i(t, \bm{X}_t, \bm{\alpha}_t) = \sigma(t, X^i_t, \nu^N_t, \alpha^i_t), \\ & f^i(t, \bm{X}_t, \bm{\alpha}_t) = f(t, X^i_t, \nu_t^N, \alpha^i_t), \quad g^i(\bm{X}_T) = g(X_T^i, \mu_T^N),
\end{aligned}
\end{equation}
where $\nu^N_t = \frac{1}{N} \sum_{j=1}^N \delta_{(X^j_t, \alpha^j_t)}$ is the empirical state-action distribution of the population and $\mu^N_t = \frac{1}{N} \sum_{j=1}^N \delta_{X^j_t}$  is its first marginal, which corresponds to the state distribution. We keep the same notation $\sigma^0$ for simplicity, with a slight abuse of notation. 

Then the cost associated to a strategy profile $\bm{\alpha}$ is defined as
\begin{equation}\label{def:MFG-cost}
J^i(\bm{\alpha}) = \EE\left[\int_0^T f(t, X^i_t, \nu_t^N, \alpha^i_t) \ud t + g(X_T^i, \mu_T^N)\right],
	\end{equation}
where the processes $X^j$, $j=1,\dots,N$, solve the SDE system
\begin{equation}%
\ud X_t^j = b(t, X^j_t, \nu^N_t, \alpha^j_t) \ud t + \sigma(t, X^j_t, \nu^N_t, \alpha^j_t) \ud W_t^j + \sigma^0(t, X^j_t, \nu^N_t, \alpha^j_t) \ud W_t^0, \quad X_0^j \sim \mu_0, \quad j \in \mc{I},
\end{equation}
where the initial positions are i.i.d., with $\nu$ and $\mu$ being as above the flows of empirical state-action and empirical state distributions. The influence of a given player on the dynamics and the cost of another player occurs only through the empirical distribution flow $\nu^N$. So when $N$ increases, the influence of each player decreases. By symmetry, we can expect that in the limit it is sufficient to study the problem for a single representative player.

To formulate the MFG, let  $\nu=(\nu_t)_{0\le t\le T}$ be a stochastic distribution flow adapted to the filtration generated by $W^0$, which is interpreted at the evolution of the population's state-action configuration. Let $\alpha$ be an open-loop control. A representative player's dynamics are given by,
\begin{align}
\label{eq:dyn-X-general-MFG}
\begin{dcases}
    \ud X_t^{\nu, \ctrl} 
    = b(t, X_t^{\nu, \ctrl}, \nu_t, \ctrl_t) \ud t + \sigma(t, X_t^{\nu, \ctrl}, \nu_t, \ctrl_t) \ud W_t 
	 + \sigma^0(t, X_t^{\nu, \ctrl}, \nu_t, \ctrl_t) \ud W^0_t, \quad t \ge 0
	 \\
	 X^{\nu, \ctrl}_0 \sim \mu_0,
\end{dcases}
\end{align}
where $W$ is a standard $m$-dimensional Brownian motion independent of $W^0$. For a representative player, the cost associated to using the control $\alpha$ when the population is given by the distribution flow  $\nu=(\nu_t)_{0\le t\le T}$ is defined as
\begin{align}
\label{subchapRCML-num-eq:def-J-MFG}
	J^{MFG}(\ctrl; \nu) = \EE \left[\int_0^T f(t, X_t^{\nu, \ctrl}, \nu_t, \ctrl_t ) \ud t + g(X_T^{\nu, \ctrl}, \mu_T) \right],
\end{align}
under the constraint that the  process $X^{\nu, \ctrl} = (X_t^{\nu, \ctrl})_{t \ge 0}$ solves the SDE~\eqref{eq:dyn-X-general-MFG}.

\begin{defn}[Mean-field Nash equilibrium] 
\label{def:MFGeq-finitehorizon}
Consider the MFG problem introduced above. A pair $(\hat{\nu},\hat{\ctrl})$ consisting of a stochastic flow $\hat{\nu}=(\hat \nu_t)_{0\le t\le T}$ of probability measures in $\cP_2(\RR^d)$ adapted to the common noise filtration and an open-loop control $\hat{\ctrl} = (\hat{\ctrl}_t)_{t \in [0,T]}$ is a mean-field Nash equilibrium if it
satisfies the following two conditions
\begin{enumerate}
	\item $\hat{\ctrl}$ minimizes $J^{MFG}(\cdot; \hat{\nu})$;
	\item For all $t \in [0,T]$, $\hat{\nu}_t$ is the probability distribution of $(X_t^{\hat{\nu}, \hat{\ctrl}}, \hat{\ctrl}_t)$ conditioned on $W^0$.
\end{enumerate}
\end{defn}
Note that, in the first condition, $\hat\nu$ is fixed when an infinitesimal agent performs their optimization. The second condition ensures that if all the players use the control $\hat{\ctrl}$, the law of their individual states and actions is indeed $\hat{\nu}$. 
The original formulation of MFGs~\cite{MR2295621} considers interactions through the state distribution only. MFGs with interactions through the joint distribution of state and actions, as presented in the above definition, are sometimes referred to as extended MFGs or MFGs of controls, see {\it e.g.}, \cite{GoPaVo:14,GoVo:16,cardaliaguet2017mean,kobeissi2022classical,laurieretangpi2022convergence}.

\begin{remark}
    For a given model (cost functions and dynamics), one can look for a mean field Nash equilibrium (Definition~\ref{def:MFGeq-finitehorizon}) or a mean field social optimum (Definition~\ref{def:mfc-optimum}). As already mentioned, the first notion corresponds to a situation in which the agents selfishly minimize their individual cost, while the second notion corresponds to a situation in which the agents cooperate to minimize the social cost. From the mathematical viewpoint, the key difference is that in Definition~\ref{def:MFGeq-finitehorizon}, the distribution is fixed when one looks for an optimal control, while in Definition~\ref{def:mfc-optimum}, the control directly influences the distribution. As a consequence an MFG (Nash equilibrium) problem is a fixed point problem, while an MFC (social optimum) is an optimization problem (or more precisely an optimal control problem for McKean-Vlasov dynamics). In general, the two solutions are different, which leads to a notion of price of anarchy. See e.g.~\cite{graber2016linear,MR3968548,cardaliaguet2019efficiency,AMSnotesLauriere,angiulifouquelauriere2020unified} for more details and comparisons of MFG and MFC solutions.
\end{remark}

Next, we review several ways to characterize the mean field Nash equilibrium concept using analytical and probabilistic techniques.

\vskip 6pt
\paragraph{PDE system}
\label{sec:MFG-PDE-system}

For simplicity, let us assume that there is no common noise. We assume that there exists an equilibrium, and we denote by $\hat{\nu} = (\hat{\nu}_t)_{t \ge 0}$ the associated mean-field flow of distributions. When considering Markovian controls, one can define the value function $u$ by
\begin{equation}
    u(t, x) = \inf_\alpha \EE\left[\int_t^T f(s, X_s, \hat{\nu}_s, \alpha_s) \ud s + g(X_T, \hat \mu_T) \vert X_t = x\right]. 
\end{equation}
As in standard OC problems,  
under suitable conditions, $u(t, x)$ solves the HJB equation:
\begin{equation}
    \label{def:control-HJB-MFG}
    \begin{dcases}
   	\partial_t u (t,x)  + \min_{\alpha \in \mc{A}} H(t, x, \hat{\nu}_t, \grad_x u(t,x), \Hess_x u(t, x), \alpha) = 0,
   	\\
   	u(T, x) = g(x,  \hat{\mu}_T),
    \end{dcases}
\end{equation}
where
\begin{equation}
\label{def:control-H-MFG}
    H(t, x, \nu, p, q, \alpha) = b(t, x, \nu, \alpha) \cdot p + \half \Tr(\sigma(t, x, \nu, \alpha)\sigma(t,x, \nu, \alpha)\transpose q) + f(t, x, \nu, \alpha).
\end{equation} 
If \eqref{def:control-HJB-MFG} has a classical solution, then the optimal control is given by
\begin{equation*}
	\hat{\ctrl}(t, x)
	= \ctrl(t, x, \hat{\nu}_t, \grad_x u(t,x), \Hess_x u(t,x)),
\end{equation*}
where
\begin{equation*}
	\ctrl(t, x, \nu, p, q)
	= 
	\argmin_{\alpha \in \mc{A}} H(t, x, \nu, p, q, \alpha).
\end{equation*}
The consistency condition for the equilibrium mean-field flow is equivalent to: the state distribution flow $\hat{\mu} = (\hat{\mu}_t)_{t \ge 0}$ solves the following Kolmogorov-Fokker-Planck (KFP) PDE
\begin{equation}
    \label{eq:MFG-general-KFP}
    \begin{dcases}
        \displaystyle \partial_t \hat{\mu}(t,x)  - \sum_{i,j} \frac{\partial^2}{\partial_{x_i}\partial_{x_j}}\left( \hat{D}_{i,j}(t,x)\hat{\mu}(t,x)\right) + \diver\Bigl( \hat{\mu}(t,x) \hat{b}(t,x)\Bigr) = 0,
        \\
        \hat{\mu}(0) = \mu_0,
    \end{dcases}
\end{equation}
where
\begin{equation}
    \label{eq:def-KFP-hat-D-b}
    \hat{D}(t,x) = \frac{1}{2} \sigma(t, x, \hat{\nu}_t, \hat{\alpha}(t,x))\sigma(t,x, \hat{\nu}_t, \hat{\alpha}(t,x))\transpose, \qquad \hat{b}(t,x) = b(t,x,\hat{\nu}_t,\hat{\ctrl}(t,x)),
\end{equation}
and the state-action distribution $\hat{\nu}_t$ at time time $t$ is the push forward of $ \hat{\mu}_t$ by $(I_d, \hat\alpha(t,\cdot))$, which we will denote by 
$
    \hat{\nu}_t =  \hat{\mu}_t \circ (I_d, \hat\alpha(t,\cdot))^{-1}. 
$
The forward-backward PDE system \eqref{def:control-H-MFG}--\eqref{eq:MFG-general-KFP} characterizes the mean field Nash equilibrium. We refer to {\it e.g.}, \cite{kobeissi2022classical} for the existence of classical solutions to such PDE systems under suitable assumptions. 

\begin{remark}
MFC problems also give rise to analogous forward-backward PDE systems, except that the solution $u$ of the backward equation is not interpreted as a value function of an optimal control problem but rather as an adjoint state. We refer to \cite{MR3134900,achdoulauriere2015systemmfc} for more details. The KFP equation remains the same, but the HJB equation has one extra term reflecting the fact that the whole population performs the optimization simultaneously.
\end{remark}

In the presence of common noise, the HJB and KFP equations become stochastic. We will not discuss this system in the sequel, and refer the interested readers to \cite{peng1992stochastichjb} for the derivation of stochastic HJB equations and~\cite{carmonadelarue2014mastereqlarge,cardaliaguetdelaruelasrylions2019master} for stochastic HJB-KFP systems arising in MFG (with the state distribution only).

\vskip 6pt
\paragraph{FBSDE system}

We now review the characterization of MFG equilibria using BSDEs. As for standard OC (see Section~\ref{sec:control_intro}), BSDEs can be used to characterize the value function or its gradient. For simplicity, we assume that there is no common noise. We further assume that the volatility of the idiosyncratic noise is uncontrolled, in which case $\hat \ctrl$ is independent of $\Hess_x u(t,x)$ and the PDE \eqref{def:control-HJB-MFG} becomes semi-linear:
\begin{multline}
	\partial_t u (t,x)  + \half \Tr(\sigma(t, x, \hat{\nu}_t)\sigma(t,x, \hat{\nu}_t)\transpose \Hess_x u(t, x)) + b(t, x, \hat{\nu}_t,	\hat{\ctrl}(t, x, \hat{\nu}_t, \grad_x u(t,x))) \cdot  \grad_x u(t,x) \\
	+ f(t, x, \hat{\ctrl}(t, x, \hat{\nu}_t, \grad_x u(t,x))) = 0.
\end{multline}
Suppose that there exist functions $\mu(t, \nu, x)$ and $h(t, x, \nu, z)$ such that 
\begin{multline*} 
    \tilde{b}(t, \hat{\nu}_t, x) \cdot \grad_x u(t,x) + h(t, x, \hat{\nu}_t, \sigma(t,x) \transpose \grad_x u(t,x)) 
    \\= b(t, x, \hat{\nu}_t, \hat{\ctrl}(t, x, \hat{\nu}_t, \grad_x u(t,x))) \cdot  \grad_x u(t,x) 
+ f(t, x, \hat{\nu}_t, \hat{\ctrl}(t, x, \grad_x u(t,x))).
\end{multline*}
Then the non-linear Feynman-Kac formula (see \cite{PaPe:90}) gives the following BSDE interpretation of $u(t, x)$:
\begin{equation}\label{def_MF_BSDEreform}
    \begin{dcases}
    \ud \mc{X}_t = \tilde{b}(t, \hat{\nu}_t, \mc{X}_t) \ud t + \sigma(t, \hat{\nu}_t, \mc{X}_t) \ud W_t, \quad \mc{X}_0 \sim \mu_0, \\
    \ud \mc{Y}_t = -h(t, \hat{\nu}_t, \mc{X}_t, \mc{Z}_t)\ud t + \mc{Z}_t \ud W_t, \quad \mc{Y}_T = g(\mc{X}_T, \hat{\mu}_T),
    \end{dcases}
\end{equation}
by the relation
$$
\mc{Y}_t = u(t, \mc{X}_t), \quad \mc{Z}_t = \sigma(t, \hat{\nu}_t, \mc{X}_t) \transpose \grad_x u(t, \mc{X}_t).
$$
Moreover, the optimal value is given by $\mathbb{E}[\mc{Y}_0] = \mathbb{E}[u(0, \mc{X}_0)]$. This BSDE characterizes the value function for a representative player given the mean field flow $\hat\nu$. Then, the consistency condition reads: 
$$
    \hat{\nu}_t = \Law\left(\mc{X}_t, \hat{\ctrl}(t, \mc{X}_t, \hat{\nu}_t, (\sigma(t, \hat{\nu}_t, \mc{X}_t) \transpose)^{-1}\mc{Z}_t)\right).
$$

In the controlled volatility case, the PDE \eqref{def:control-HJB-MFG} is fully nonlinear, and its solution is connected to a solution of the 2BSDE, see~\cite{cheridito2007second} and Section~\ref{sec:control_intro}.  

The Pontryagin stochastic maximum principle provides the connection to the FBSDE. Define the generalized Hamiltonian $\mc{H}$ by
\begin{equation}
\label{eq:Hamiltonian-bsde-control-MFG}
\mc{H}(t,x,\nu,y,z,\alpha) = b(t, x, \nu, \alpha) y + \Tr(\sigma\transpose(t, x, \nu, \alpha)z) + f(t, x, \nu, \alpha).
\end{equation}
If the Hamiltonian $\mc{H}$ is convex in $(x, \alpha)$, 
and $(X_t, Y_t, Z_t)$ solve
\begin{equation}
    \begin{dcases}
    \ud X_t = b(t, X_t, \hat{\nu}_t, \hat \alpha_t) \ud t +\sigma(t, X_t, \hat{\nu}_t, \hat \alpha_t) \ud W_t, \qquad X_0 \sim \mu_0,\\
    \ud Y_t = -\grad_x \mc{H}(t, X_t, \hat{\nu}_t, Y_t,  Z_t,  \hat \alpha_t ) \ud t + Z_t \ud W_t, \quad Y_T = \partial_x g(X_T, \hat{\mu}_T),
    \end{dcases}
\end{equation}
such that $\hat \alpha$ minimizes $\mc{H}$ along $(X_t, \hat{\nu}_t, Y_t, Z_t)$, then $\hat \alpha$ is the optimal control. If the value function is smooth enough, then 
\begin{equation}
     Y_t = \grad_x u(t,  X_t), \quad Z_t = \sigma(t, X_t, \hat{\nu}_t, \hat \alpha)\transpose\Hess_x u(t,  X_t).
\end{equation}
In this case, the consistency condition for the equilibrium mean field flow $\hat{\nu}$ reads
$$
    \hat{\nu}_t = \Law\left(X_t, \hat{\ctrl}(t, X_t, \hat{\nu}_t, (\sigma(t, \hat{\nu}_t, X_t) \transpose)^{-1} Y_t, Z_t)\right).
$$
When there is common noise, the FBSDE system becomes
\begin{equation}
    \begin{dcases}
    \ud X_t = b(t, X_t, \hat{\nu}_t, \hat \alpha_t) \ud t +\sigma(t, X_t, \hat{\nu}_t, \hat \alpha_t) \ud W_t
    +\sigma^0(t, X_t, \hat{\nu}_t, \hat \alpha_t) \ud W^0_t, \qquad X_0 \sim \mu_0,\\
    \ud Y_t = -\grad_x \mc{H}(t, X_t, \hat{\nu}_t, Y_t,  Z_t,  Z^0_t,  \hat \alpha_t ) \ud t + Z_t \ud W_t + Z^0_t \ud W^0_t, \quad Y_T = \partial_x g(X_T, \hat{\mu}_T),
    \end{dcases}
\end{equation}
where, compared with~\eqref{eq:Hamiltonian-bsde-control-MFG},  the definition of $\mc{H}$ includes an extra term $\Tr({\sigma^0}\transpose(t, x, \nu, \alpha)z^0)$.

\vskip 6pt

\begin{remark} 
MFC problems also lead to analogous FBSDE systems. In the absence of common noise, Pontryagin's maximum principle is derived for instance in~\cite{MR3045029} and~\cite{acciaio2019extendedmfc} when the interactions are through the state or the state-action distributions, respectively. This leads to a BSDE with an extra term accounting for the variation of the distribution during the optimization of the control.
\end{remark}

\vskip 6pt 
All the above systems are particular cases of the following generic system of FBSDEs of McKean-Vlasov type (MKV FBSDE for short)
\begin{equation}
\label{eq:MKV-FBSDE-general}
\left\{
\begin{aligned}
	\ud X_t
	=
	\,& B\left(t, X_t, \Law(X_t, Y_t, Z_t |W^0), Y_t, Z_t, Z^0_t \right) \ud t
		\\
		&\qquad
		+  \sigma(t, X_t, \Law(X_t, Y_t, Z_t |W^0), Y_t, Z_t, Z^0_t) \ud W_t 
		\\
		&\qquad + \sigma^0(t, X_t, \Law(X_t, Y_t, Z_t |W^0), Y_t, Z_t, Z^0_t) \ud W^0_t,
	\\
	\ud Y_t
	=
	\,& - F\Big(t, X_t, \Law(X_t, Y_t, Z_t |W^0), Y_t, \sigma\transpose(t, X_t, \Law(X_t, Y_t, Z_t |W^0), Y_t, Z_t, Z^0_t) Z_t, 
	\\
	&\quad\qquad{\sigma^0}\transpose(t, X_t, \Law(X_t, Y_t, Z_t |W^0), Y_t, Z_t, Z^0_t) Z^0_t \Big) \ud t
		\\
		&\qquad
		+ Z_t \ud W_t + Z^0_t \ud W^0_t,
    \\
    \Law(X_0) &= \mu_0, \qquad Y_T = G(X_T, \Law(X_T|W^0)).
\end{aligned}
\right.
\end{equation}

\begin{remark}
When there is no common noise, $W^0$ and $Z^0$ are dropped, and the system becomes
\begin{equation}
\label{eq:MKV-FBSDE-general-no-CN}
\left\{
\begin{aligned}
	\ud X_t
	=
	& B\left(t, X_t, \Law(X_t, Y_t, Z_t), Y_t, Z_t \right) \ud t
			+  \sigma(t, X_t, \Law(X_t, Y_t, Z_t), Y_t, Z_t) \ud W_t ,
	\\
	\ud Y_t
	=
	& - F\Big(t, X_t, \Law(X_t, Y_t, Z_t), Y_t, \sigma\transpose(t, X_t, \Law(X_t, Y_t, Z_t), Y_t, Z_t) Z_t\Big) \ud t
			+ Z_t \ud W_t,
    \\
    \Law(X_0) &= \mu_0, \qquad Y_T = G(X_T, \Law(X_T)).
\end{aligned}
\right.
\end{equation}
When the interactions are not through the state-action distribution but through the state distribution only, $\Law(X_t, Y_t, Z_t)$ is reduced to $\Law(X_t)$.
\end{remark}

\paragraph{Master equation}
\label{sec:background-master-eq}

As mentioned earlier, in the PDE system~\eqref{def:control-HJB-MFG}--\eqref{eq:MFG-general-KFP}, $u$ plays the role of the value function of a representative player when the rest of the population is at equilibrium. This function depends explicitly on $t$ and $x$ but, intuitively, a player's value function can also depend on the population distribution. When there is no common noise, this distribution evolves in a deterministic way, so knowing $\mu_0$ and $t$ as well as the control used by the population (which is the equilibrium control $\hat\alpha$, assuming the population is at equilibrium) is enough to recover $\hat\mu(t)$, {\it e.g.}, by solving the corresponding KFP equation~\eqref{eq:MFG-general-KFP}. However, we can make this dependence explicit by considering a function $\cU: [0,T] \times \RR^d \times \cP(\RR^d) \to \RR$ such that
\begin{equation}
\label{eq:masterfield-to-u}
    \cU(t,x,\hat\mu(t)) = u(t,x),
\end{equation}
where $\hat\mu = (\hat\mu(t))_t$ is the mean-field equilibrium distribution flow. This correspondence is even more useful when common noise influences the dynamics of the players. In this case, $u(t,x)$ is a random variable whereas $\cU$ is still a deterministic function and the lefthand side of~\eqref{eq:masterfield-to-u} is random only due to $\hat\mu(t)$. This function $\cU$ has been instrumental in proving the convergence of finite-player Nash equilibria towards mean field Nash equilibria, see~\cite{cardaliaguetdelaruelasrylions2019master} for more details.

It turns out that, under suitable conditions, $\cU$ satisfies the PDE that we will present below, introduced by Pierre-Louis Lions and called the Master equation. It involves partial derivatives with respect to the probability measure argument in $\cU$. We say that a function $F: \cP(\RR^d) \to \RR$ is $\mathcal{C}^1$ if there exists a continuous map $\displaystyle\frac{\delta F}{\delta \mu}: \cP(\RR^d) \times \RR^d \to \RR$ such that, for any $\mu,\mu' \in \cP(\RR^d)$,
$$
\lim_{s \to 0^+} \frac{F((1-s)\mu + s\mu') - F(\mu)}{s} = \int_{\RR^d} \frac{\delta F}{\delta \mu}(\mu,y) d(\mu'-\mu)(y).
$$
The derivative $\displaystyle\frac{\delta F}{\delta \mu}$ is sometimes referred to as the flat derivative. 
If $\displaystyle\frac{\delta F}{\delta \mu}$ is of class $\mathcal{C}^1$ with respect to the second variable, the intrinsic derivative $\partial_\mu F:\cP(\RR^d) \times \RR^d \to \RR$ is defined by
$$
    \partial_\mu F(\mu,y) = \partial_y \frac{\delta F}{\delta \mu}(\mu,y).
$$
We will write $\partial_\mu F(\mu)(y)$ instead of $\partial_\mu F(\mu,y)$. For more details, we refer to the lectures of Pierre-Louis Lions~\cite{PLL-CDF}, as well as \cite{Cardaliaguet-2013-notes} and \cite[Chapter 5]{carmona2018probabilistic}.

We can now present the Master equation. To the best of our knowledge, the theory has not yet been developed for the general MFG model described above. We thus consider the case in which the volatility is not controlled, and the interactions are only through the state distribution instead of the state-action distribution. For the sake of brevity, we omit the derivation and refer to {\it e.g.}~\cite[Section 4.4]{carmona2018probabilistic2}. The Master equation is the following backward PDE, posed on the space $[0,T] \times \RR^d \times \cP_2(\RR^d)$,
\begin{align*}
    &\partial_t \cU(t,x,\mu)
    \\
    &\quad +b(t,x,\mu,\hat\alpha(t, x, \mu, \partial_x \cU(t,x,\mu))) \cdot \partial_x \cU(t,x,\mu)
    \\
    &\quad + \int_{\RR^d} b(t,v,\mu,\hat\alpha(t,v,\partial_x \cU(t,v,\mu))) \cdot \partial_\mu \cU(t,x,\mu)(v) d \mu(v)
    \\
    &\quad + \frac{1}{2} \Tr \left[ (\sigma \sigma\transpose + \sigma^0 (\sigma^0)\transpose)(t,x,\mu) \partial_{xx}^2 \cU(t,x,\mu)\right]
    \\
    &\quad + \frac{1}{2} \int_{\RR^d} \Tr \left[ (\sigma \sigma\transpose + \sigma^0 (\sigma^0)\transpose)(t,v,\mu) \partial_v\partial_\mu \cU(t,x,\mu)(v)\right] d \mu (v)
    \\
    &\quad + \frac{1}{2} \int_{\RR^{2d}} \Tr \left[ (\sigma \sigma\transpose + \sigma^0 (\sigma^0)\transpose)(t,v,\mu) \partial_\mu^2 \cU(t,x,\mu)(v,v')\right] d \mu (v) d \mu (v')
    \\
    &\quad + \int_{\RR^d} \Tr \left[ (\sigma^0(t,x,\mu) (\sigma^0)\transpose)(t,v,\mu) \partial_x\partial_\mu \cU(t,x,\mu)(v)\right] d \mu (v)
    \\
    &\quad + f(t,x,\mu,\hat\alpha(t, x, \mu, \partial_x \cU(t,x,\mu)))
    =0,
\end{align*}
for $t \in [0,T]$, $x \in \RR^d$ and $\mu \in \cP_2(\RR^d)$, and with the terminal condition: for every $x \in \RR^d$ and $\mu \in \cP_2(\RR^d)$,
$$
    \cU(T,x,\mu) = g(x,\mu).
$$
For more details on the analysis of this PDE, we refer the interested reader to the monographs~\cite{cardaliaguetdelaruelasrylions2019master}, \cite[Chapters 4 to 7]{carmona2018probabilistic2}, and~\cite{ChassagneuxCrisanDelarue_Master} concerning the existence of classical solutions under suitable conditions.

\subsubsection{Direct parameterization}
\label{sec4_MFG_with_CN}

As discussed earlier (Section~\ref{sec:control-direct}), the direct parameterization approach for optimal control (Section~\ref{sec:open-loop}) can be extended to finite-player games. It can also be extended to mean-field games by updating alternatively the control (using the direct parameterization method) and the mean-field instead of updating individually all the player's controls as in finite-player games. This idea can be applied to various classes of controls (\textit{e.g.}, open-loop and closed-loop ones). To avoid repetition, we discuss below the application of this approach to a class of MFGs with common noise, which forces to use a more general class of controls.

As seen in Definition~\ref{def:MFGeq-finitehorizon}, a mean-field equilibrium is a standard control problem (corresponding to the first item in the definition) plus a fixed point problem (corresponding to the second item). Motivated by MFG models with common noise,   
Min and Hu \cite{MinHu:21} proposed an algorithm called Sig-DFP utilizing the concept of signature in rough path theory \cite{LyonsTerryJ2007DEDb} and fictitious play from game theory \cite{Br:49,Br:51}. Signature is used to accurately represent the conditional distribution of the state given the common noise, and fictitious play is used to solve the fixed-point problem and identify the equilibrium \cite{cardaliaguet2015learning}.

For a path $x:[0,T]\to \RR^d$, the $p$-variation is defined by 
\begin{equation}
    \|x\|_{p} = \left( \sup_{D\subset[0,T]} \sum_{n=0}^{r-1} \|x_{t_{n+1}}-x_{t_n}\|^p \right)^{1/p},
\end{equation} 
where $D \subset [0,T]$ denotes a partition $0 \leq t_0 < t_1 < \ldots < t_r \leq T$. Let $T((\R^d))=\bigoplus_{k=0}^\infty (\R^d)^{\bigotimes k}$ be the tensor algebra. Let  $\mathcal{V}^p([0,T], \R^d)$ be the space of continuous mappings from $[0,T]$ to $\R^d$ with finite $p$-variation, equipped with norm $\|\cdot\|_{\mathcal{V}^p}=\|\cdot\|_{\infty}+\|\cdot\|_{p}$. 

\begin{defn}[Signature]
Let $X\in \mathcal{V}^p([0,T], \R^d)$ such that the following integral is well defined. The signature of $X$, denoted by $S(X)$, is the element of $T((\R^d))$ defined by $S(X) = (1, X^1, \cdots, X^k \cdots)$ with
\begin{equation}\label{def:signature}
    X^k = \int_{0<t_1<t_2<\cdots<t_k<T} \ud X_{t_1}\otimes\cdots\otimes \ud X_{t_k}.
\end{equation}
Denoting by $S^M(X)$ the truncated signature of $X$ of depth $M$, {\it i.e.}, $S^M(X) = (1, X^1, \cdots, X^M)$ which has the dimension $\frac{d^{M+1}-1}{d-1}$.
\end{defn}
In the current setting, $X$ is a semi-martingale, thus equation \eqref{def:signature} is understood in the Stratonovich sense. The signature has many nice properties, including the following ones. First, it characterizes paths uniquely up to the tree-like equivalence, and the equivalence is removed if at least one dimension of the path is strictly increasing \cite{boedihardjo2014signature}. Therefore, in practice one usually augments the original path $X_t$ with the time dimension, {\it i.e.}, working with $\hat{X}_t = (t, X_t)$ since $S(\hat{X})$ characterizes paths $\hat{X}$ uniquely. Second, terms in the signature present a factorial decay property \cite{lyons2002system}, which implies that a path can be well approximated with just a few terms in the signature ({\it i.e.}, a small $M$). Last, As a feature map of sequential data, the signature has a universality property \cite{NEURIPS2019_deepsig}, which is summarized below.

    Let $p\ge 1$ and $f: \mathcal{V}^p([0,T], \R^d)\to \R$ be a continuous function. For any compact set $K\subset \mathcal{V}^p([0,T], \R^d)$, if $S(x)$ is a geometric rough path (see \cite[Definition 3.13]{LyonsTerryJ2007DEDb} for a detailed definition) for any $x\in K$, then for any $\eps >0$, there exists a linear functional $l$ in the dual space of  $\in T((\RR^d))$ such that
    \begin{equation}
        \sup_{x\in K} |f(x) - \langle l, S(x) \rangle| <\epsilon.
    \end{equation}

Motivated by the unique characterization of $(W_s^0)_{s \in [0,t]}$ by $S(\hat W_t^0)$ and the factorial decay property, one can approximate 
\begin{equation}
  \nu_t \equiv \mc{L}(X_t, \alpha_t \vert \mc{F}_t^0 ) = \mc{L}(X_t, \alpha_t \vert S(\hat W_t^0)), \label{eq:propbysig}
\end{equation}
by $\mc{L}(X_t, \alpha_t \vert S^M(\hat W_t^0))$, for $\hat W^0_t = (t, W^0_t)$. 
In particular, if the mean-field interaction is through moments $\bar\nu_t = \EE[\iota(X_t, \alpha_t) \vert \F_t^0]$, for some measurable function $\iota$, the approximation can be arbitrarily accurate for sufficiently large $M$, see \cite[Lemma 4.1]{MinHu:21}. Then \cite{MinHu:21} proposed to use the approximation
\begin{equation}\label{eq:mfg_cn_linearregression}
    \bar\nu_t \approx \langle \tilde l, S^M(\hat W_t^0) \rangle, \; \text{ where } \tilde l =  \argmin_{\bm\beta} \|\bm y - \bm X \bm\beta\|^2,  \;
    \bm y = \{\iota(X_t(\omega_i), \alpha_t(\omega_i))\}_{i=1}^N,\;    \bm X = \{S^M(\hat{W}^{0}_{t}(\omega_i))\}_{i=1}^N,
\end{equation}
where $\omega_i$ denotes the $i^{th}$ sample path. The rationale behind this approximation is the universality of signatures and the interpretation of ordinary linear regression: the least square minimization gives the best possible prediction of $\EE[\bm y | \bm X]$ using linear relations. Once $\tilde l$ is obtained, the prediction on an unseen common noise is efficient: $\bar\nu_t(\tilde\omega)\approx \langle \tilde l, S^M(\hat{W}^0_{t}(\tilde\omega)) \rangle$ for any $\tilde \omega$ and $t$.

Then finding the mean-field equilibrium is broken down into the following steps. We start with an initial value $\bar\nu^{(0)}$. Then, we solve the standard control problem given $\bar\nu^{(0)}$ in \eqref{subchapRCML-num-eq:def-J-MFG} in the spirit of \cite{han2016deep-googlecitations}. From here, we approximate $\bar\nu^{(1)}$ via signature using \eqref{eq:mfg_cn_linearregression}, i.e., compute $\tilde l^{(1)}$. These steps are repeated until convergence. The update of $\bar\nu_t$ from step to step is done by averaging $\tilde l^{(n)}$. The Sig-DFP algorithm consists of repeatedly solving \eqref{eq:dyn-X-general-MFG}--\eqref{subchapRCML-num-eq:def-J-MFG} for a given $\bar\nu$ using deep learning in the spirit of \cite{han2016deep-googlecitations}, and passing the obtained $\bar\nu$ to the next iteration by using signatures. A flowchart illustrating the ideas is given in Figure~\ref{fig:algo}.

\begin{figure}
    \centering
    \includegraphics[width = 0.7\textwidth]{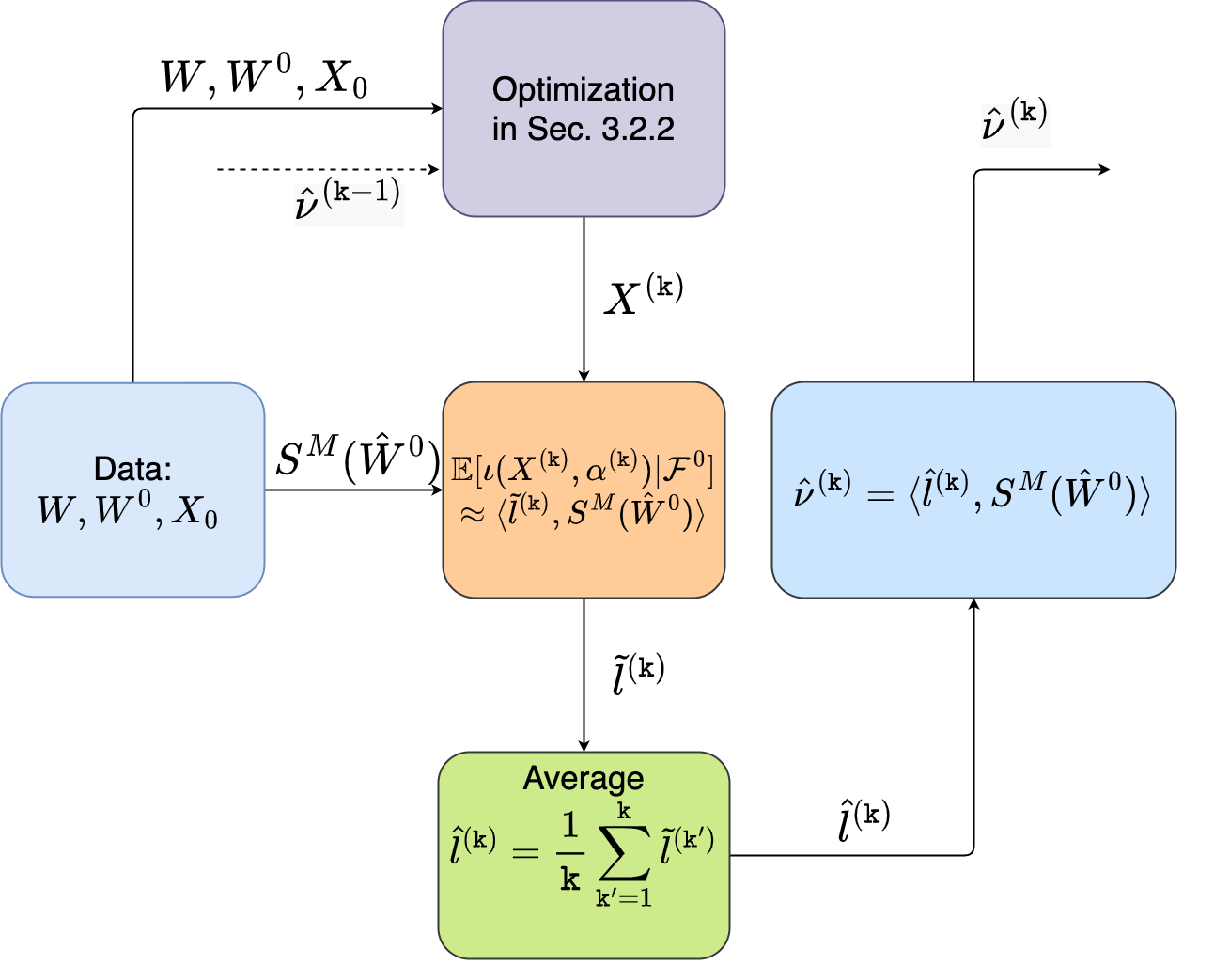}
    \caption{Flowchart of one iteration in the Sig-DFP Algorithm. Input: idiosyncratic noise $W$, common noise $W^0$, initial position $X_0$ and vector $\hat{\nu}^{(\mt{k}-1)}$ from the last iteration. Output: vector $\hat{\nu}^{(\mt{k})}$ for the next iteration. %
    }
    \label{fig:algo}
\end{figure}

More precisely, at each step, given a proxy $\hat \nu^{(\mt{k-1})}$ of the equilibrium distribution $\hat \nu$, the problem \eqref{eq:dyn-X-general-MFG}--\eqref{subchapRCML-num-eq:def-J-MFG} becomes a standard stochastic control and is solved by using the direct parameterization approach reviewed in Section~\ref{sec:direct-global-local}: the loss function will be the discretized version of \eqref{subchapRCML-num-eq:def-J-MFG},  $\check X$ will be follow the Euler scheme of \eqref{eq:dyn-X-general-MFG} with $\nu$ replaced by $\hat \nu^{(\mt{k}-1)}$ and so do $f$ and $g$, and the control $\alpha_{t_n}$ is parameterized by a neural network of the following form
\begin{equation}
     \alpha_{t_n}= \alpha(t_{n}, \check X_{t_n}, \hat\nu^{(\mt{k}-1)}_{t_n}; \theta),
\end{equation}
which takes $\hat\nu^{(\mt{k}-1)}_{t_n}$ as an extra input on top of $(t_{n}, \check X_{t_n})$. 
The optimizer $\theta^\ast$ obtained in this way gives $\alpha_{t_n}^{(\mt k)}$, with which the optimized state process paths are simulated. The conditional law, denoted by $\nu^{(\mt k)}$, is approximated using signatures via \eqref{eq:mfg_cn_linearregression}. This finishes one iteration of fictitious play. Denote by $\tilde\nu^{(\mt k)}$ the approximation of $\nu^{(\mt k)}$, we then pass $\tilde\nu^{(\mt k)}$ to the next iteration via updating $\hat\nu^{(\mt k)} = \frac{1}{\mt k}\tilde\nu^{(\mt k)} + \frac{\mt k-1}{\mt k}\hat \nu^{(\mt k-1)}$ by averaging the coefficients obtained in \eqref{eq:mfg_cn_linearregression}. We summarize it in Algorithm \ref{alg:sig-dfp} in Appendix~\ref{supp:pseudocode}; see~\cite[Appendix B]{MinHu:21} for the implementation details. We remark that signatures can also be useful for generating multimodal data \cite{MinHuIch:23}.

\begin{rem}[Theoretical analysis]
In \cite[Theorems 4.1 and 4.2]{MinHu:21} Min and Hu provided a proof of convergence of this algorithm showing that, under suitable assumptions, the difference between the $\mt{k}^{th}$ iteration solution and the mean-field equilibrium can be made arbitrarily small, provided that $\mt{k}$ is sufficiently large and $\nu^{(\mt k)}$ can be approximated sufficiently well by truncated signatures.
\end{rem}

\paragraph*{Numerical illustration: MFG of optimal consumption and investment.}\label{sec:MFG-sig} 

We consider an extended heterogeneous MFG proposed by \cite{lacker2020many}, where agents interact via both states and controls. The setup is similar to \cite{LaZa:17} except for including consumption and using power utilities. Each agent's type is characterized by a random vector $\zeta=(\xi, \delta, \theta, b, \sigma, \sigma^0, \epsilon)$, and the optimization problem reads 
\begin{equation}
        \sup_{\pi, c} \EE\biggl[ 
    \int_0^T U(c_t X_t(\Gamma_t m_t)^{-\theta}; \delta)\ud t + \epsilon U(X_T m^{-\theta}_T;\delta)
    \biggr],     
 \label{def: InvestConsumpValue}
\end{equation}
where 
$U(x; \delta) = \frac{1}{1-\frac{1}{\delta}}x^{1-\frac{1}{\delta}}$, $\delta\ne 1$, is the power utility function, 
$X_t$ follows
\begin{equation}
    \ud X_t = \pi_t X_t(b \ud t + \sigma \ud W_t + \sigma^0 \ud W_t^0) - c_tX_t \ud t,  \label{def: InvestConsumpSDE}
\end{equation}
and $X_0 = \xi$. The processes $\Gamma_t = \exp\EE[\log c_t |\mcF^0_t]$ and $m_t = \exp\EE[\log X_t |\mcF^0_t]$ are the mean-field interactions from the control and state processes. Two constraints are posed: $X_t \geq 0$, $c_t \geq 0$.

The interpretation of this problem is as follows. There are infinitely many agents trade in a common investment horizon $[0,T]$, each invests between a bond (with constant return rate $r$) and a private stock with dynamics $\ud S_t/S_t = b \ud t + \sigma \ud W_t + \sigma^0\ud W_t^0$, and consume $c_t$ of his wealth at time $t$. The portion of wealth into $S_t$ is denoted by $\pi_t$. Assuming $r \equiv 0$ without loss of generality, the wealth process reads \eqref{def: InvestConsumpSDE}. Then each agent aims to maximize his utility of consumption plus his terminal wealth compared to his peers' averages $\Gamma_t$ and $m_t$. To relate it to the formulation \eqref{eq:dyn-X-general-MFG}--\eqref{subchapRCML-num-eq:def-J-MFG}, $\ctrl \equiv (\alpha^1, \ctrl^2) = (\pi, c)$ will be a 2D control with the constraint $\alpha^2_t \geq 0$, $b(t, x, \nu, \alpha) = b\alpha^1 x  - \alpha^2 x$, $\sigma(t, x, \nu, \alpha) = \sigma \alpha^1 x$, $\sigma^0(t, x, \nu, \alpha) = \sigma^0 \alpha^1 x$, $f = -U$ and $g = -U$. The explicit solutions are derived in \cite{lacker2020many} and summarized in \cite[Appendix D]{MinHu:21}.

For this experiment, we use truncated signatures of depth $M=4$. The optimal controls $(\pi_t, c_t)_{0\le t\le 1}$ are parameterized by two neural networks $\pi(\cdot; \theta)$ and $c(\cdot; \theta)$, each with three hidden layers of size 64 and taking $(\zeta, t, X_t, m_t, \Gamma_t)$ as inputs due to the nature of heterogeneous extended MFG. Due to the extended mean-field interaction term $\Gamma_t$, we will propagate two conditional distribution flows, {\it i.e.}, two linear functionals $\hat {l}^{(\mt{k})}, \hat {l}_c^{(\mt{k})}$ during each iteration of fictitious play. Instead of estimating $m_t, \Gamma_t$ directly, we estimate $\EE[\log X_t|\mcF^0_t], \EE[\log c_t|\mcF^0_t]$ by $\langle \hat {l}^{(\mt{k})}, S^4(W_t^0) \rangle$, $\langle \hat{l}_c^{(\mt{k})}, S^4(W_t^0) \rangle$ and then take exponential to get $m_t, \Gamma_t$. To ensure the non-negativity condition of $X_t$, we evolve $\log X_t$ and then take exponential to get $X_t$. For optimal consumption, $c(\cdot; \theta)$ is used to predicted $\log c_t$ and thus $\exp c(\cdot; \theta)$ gives the predicted $c_t$. With 600 iterations of fictitious play and a learning rate of 0.1 decaying by a factor of 5 for every 200 iterations, the relative $L^2$ errors for $\pi_t, c_t, m_t, \Gamma_t$ are 0.1126, 0.0614, 0.0279, 0.0121, respectively. Figure \ref{fig:OCI} compares $X$ and $m$ to their approximations, and plots the maximized utilities. Further comparison with the existing literature, different choices of truncation $M$, and the ability to deal with higher $m_0$ are also discussed in \cite{MinHu:21}.

\begin{figure}[hbt!]
    \centering
    \subfloat[$X_t$]{
         \includegraphics[width=0.32\columnwidth]{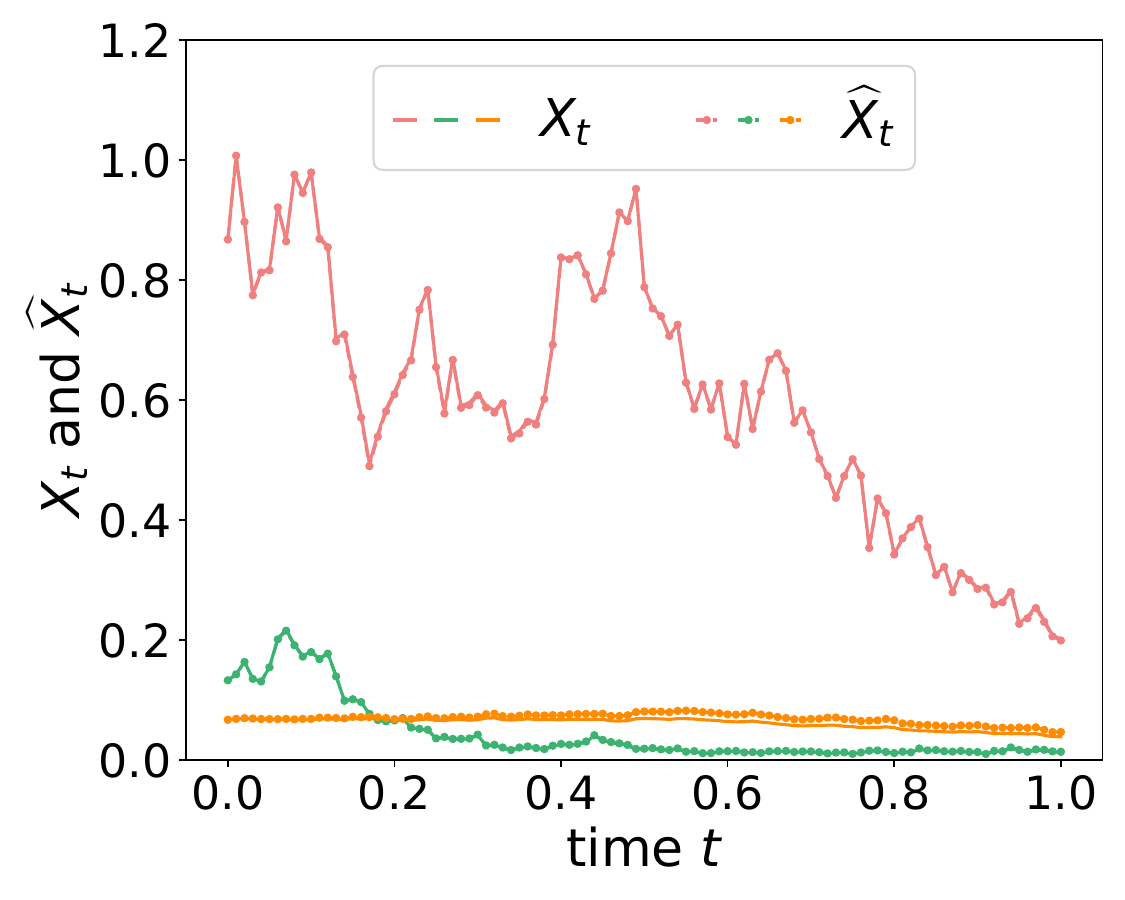}
     }
     \subfloat[$m_t =  \exp\EE(\log X_t |\mcF^B_t)$]{
         \includegraphics[width=0.32\columnwidth]{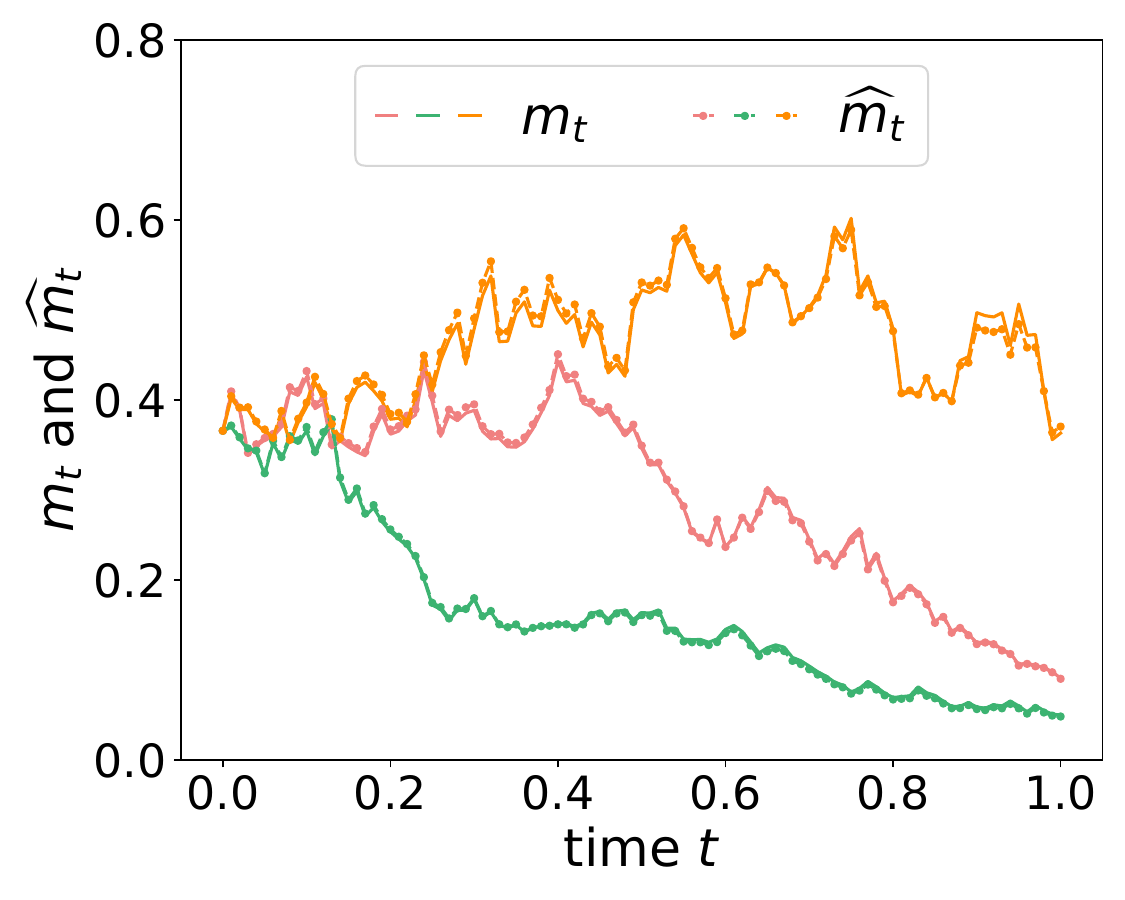}
     }
     \subfloat[Maximized Utility]{
         \includegraphics[width=0.32\columnwidth, trim = {0em 15em 0 10em}]{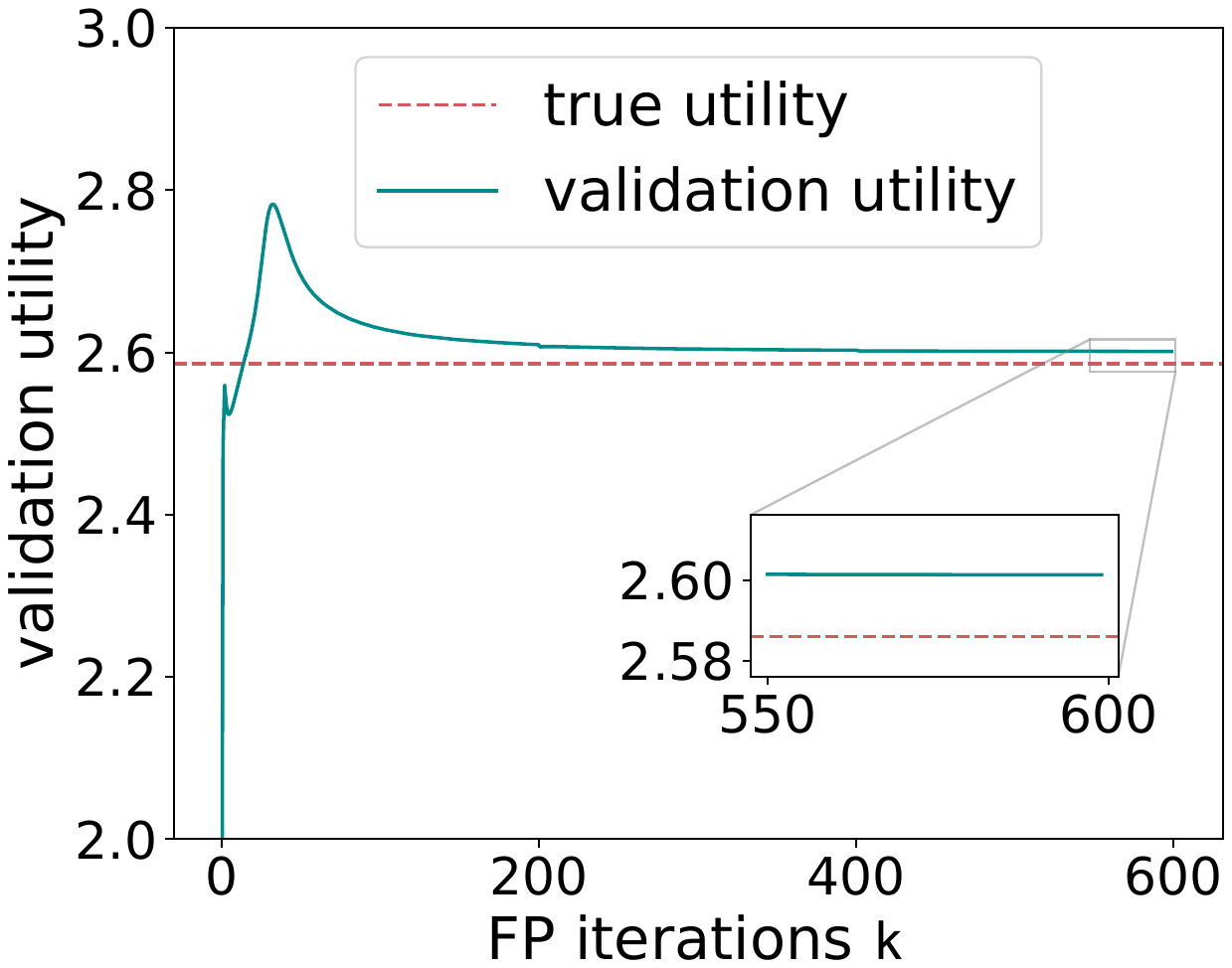}}
    
    \caption{MFG of optimal consumption and investment in Section~\ref{sec4_MFG_with_CN}. Panels (a) and (b) give three trajectories of $X_t$ and $m_t=\exp\bigl(\EE(\log X_t |\mcF^0_t)\bigr)$ (solid lines) and their approximations $\hat X_t$ and $\hat m_t$ (dashed lines) using different $(X_0, W, W^0)$ from validation data. Panel (c) shows the maximized utility computed using validation data over fictitious play iterations. Parameter choices are: $\delta \sim U(2, 2.5), b \sim U(0.25, 0.35), \sigma\sim U(0.2, 0.4), \theta, \xi \sim U(0,1), \sigma^0\sim U(0.2, 0.4)$, $\epsilon\sim U(0.5, 1)$. 
    }
    \label{fig:OCI}
\end{figure}

\subsubsection{BSDE-based deep learning algorithms}
\label{sec4_MFGFBSDE}

We now explain how to adapt the Deep BSDE method introduced in~\cite{MR3736669} and reviewed in Section~\ref{subsubsec:deepBSDE}  to mean-field FBSDEs. We recall that the principle of the method is to use neural networks to approximate $Y_0$ and $(Z_t)_{t \in [0,T]}$ and to train the neural network parameters by relying on Monte Carlo samples until the terminal condition is approximately matched. In the mean-field setting, the same idea can be used to solve forward-backward systems of McKean-Vlasov (MKV) SDEs;  see~\cite{fouque2019deep,carmona2019convergence2,germain2019numerical, hanhulong:22}.

Let us consider the FBSDE~\eqref{eq:MKV-FBSDE-general-no-CN} in the absence of common noise, with interactions through the state distribution only, and uncontrolled volatility. We rewrite the problem as: minimize over $y_0 : \RR^d \to \RR^d$ and $z: \RR_+ \times \RR^d \to \RR^{d \times m}$ the cost functional
$$
	J(y_0, z) = \EE \left[ \, \left| Y^{y_0,z}_T - G(X^{y_0,z}_T, \Law(X^{y_0,z}_T)) \right|^2 \, \right],
$$
where $(X^{y_0,z}, Y^{y_0,z})$ solves
\begin{equation}
\label{eq:MKV-FBSDE-2forward}
\begin{cases}
\ud X^{y_0,z}_t
	=
	B\left(t, X^{y_0,z}_t, \Law(X^{y_0,z}_t), Y^{y_0,z}_t \right) \ud t
		+  \sigma(t, X_t) \ud  W_t, \quad t \ge 0,
	\\
	\ud  Y^{y_0,z}_t
	=
	- F\left(t, X^{y_0,z}_t, \Law(X^{y_0,z}_t), Y^{y_0,z}_t, \sigma\transpose(t, X_t) z(t, X^{y_0,z}_t)  \right) \ud t
		+ z(t, X^{y_0,z}_t) \ud  W_t, \quad t \ge 0, 
  \\
  X^{y_0,z}_0 \sim \mu_0, \qquad Y^{y_0,z}_0 = y_0 (X_0^{y_0,z}). 
\end{cases}
\end{equation}
The above problem is an MFC problem if we view $(X^{y_0,z}_t,Y^{y_0,z}_t)$ as state and $(y_0,z)$ as control. Under suitable conditions, the optimally controlled process $(X,Y)$ solves the MKV FBSDE system~\eqref{eq:MKV-FBSDE-general-no-CN} and vice versa.

Then, to be able to implement this method, we can proceed similarly to the method described in Section~\ref{sec:directMethod}. The mean-field distribution can be approximated by an empirical distribution based on a finite population of interacting particles. Furthermore, the controls $y_0$ and $z$ can be replaced by neural networks, say $y_{\theta}$ and $z_{\omega}$ with parameters $\theta$ and $\omega$ respectively. Time can be discretized using for instance an Euler-Maruyama scheme. We thus obtain a new optimization problem over finite-dimensional parameters that can be adjusted using SGD.

\begin{remark}[Theoretical analysis]
    Motivated by numerical schemes and in particular the above adaptation of the deep BSDE method \cite{MR3736669,HaJeE:18} to MKV FBSDEs, Reisinger, Stockinger and Zhang
    in~\cite{reisinger2020posteriori} analyzed a posteriori error for approximate solutions based on a discrete time scheme and a finite population of interacting particles \cite[Theorems 3.2 and 4.3]{reisinger2020posteriori}. \cite{hanhulong:22} proposed a deep learning method for computing MKV FBSDEs with a general form of mean-field interactions, and proved that the convergence of the numerical solution obtained by the proposed method to the true solution is free of the curse of dimensionality \cite[Theorem 3.9]{hanhulong:22} by using a special class of integral probability metrics previously developed in \cite{hanhulong:21}.
\end{remark}

 Although we focus here on the continuous-state space setting, the same strategy can be applied to finite-state MFGs;  see, \textit{e.g.},~\cite{aurell2022optimalstackelberg,aurell2022finitegraphon}.

\paragraph*{Numerical illustration: a linear-quadratic mean-field game in systemic risk problems.} 
\label{sec:MFG-sysrisk-example-deepBSDE-CN}
We now consider the MFG version of the systemic risk model introduced in Section~\ref{sec:intro-LQsysrisk} which has been studied in Section~\ref{par:LQ-systemic-risk-finite-N} and revisited in Section~\ref{sec:LQsysrisk-revisited}. This MFG mean-field game has been analyzed in~\cite{CaFoSu:15}. Given a mean-field flow $\mu = (\mu_t)_{t\in[0,T]}$, the log-monetary reserves of a typical bank evolves according to the dynamics:
\begin{equation}
\label{eq:ex1_dynamics-mfg}
    \ud X_t = [a(\bar \mu_t - X_t)   + \alpha_t ] \ud t  + \sigma \left(\rho \ud W_t^0 + \sqrt{1-\rho^2} \ud W_t\right). %
\end{equation}
The standard Brownian motions $W^0$ and $W$ are independent, in which $W$ stands for the idiosyncratic noises and $W^0$ denotes the systemic shock, which is an example of common noise (see also Section~\ref{sec4_MFG_with_CN}).

The cost functions $f$ and $g$ appearing in \eqref{def:MFG-cost} takes the following form:
\begin{equation}
f(t, x, \nu,\alpha) = \half \alpha^2 - q \alpha(\bar \mu - x) + \frac{\eps}{2}(\bar \mu - x)^2,  \quad g(x, \nu) = \frac{c}{2}(\bar \mu - x)^2,
\end{equation}
which depend only on the mean $\bar \mu = \EE_{X\sim\mu}[X]$ of the first marginal of the state-action distribution $\nu$. 
It has been shown in~\cite{CaFoSu:15} that, in the MFG setting, the open-loop equilibrium is the same as the closed-loop Nash equilibrium, and it admits an explicit solution. Furthermore, it can be characterized using an MKV FBSDE system, that we omit for brevity; see~\cite{CaFoSu:15} for the details. If $\rho = 0$, then one can apply directly the method described above, with $y_0$ a function of $X_0$ and $z$ a function of $(t,X_t)$. When $\rho>0$, two changes need to be made: first, there is an extra process $Z^0$ to be learned, for which we use a neural network approximation as for $Z$; second, we expect the random variables $Z_t$ and $Z^0_t$ to depend not only on $X_t$ but also on the past of the common noise. In general, this would mean learning functions of the common noise's trajectory. However, in the present case, it is enough to rely on finite-dimensional information. Here, we add $\bar{\mu}_t$ as an input to the neural networks playing the roles of  $Z_t$ and $Z^0_t$, and this is sufficient to learn the optimal solution, see~\cite{CaFoSu:15}.

Figure~\ref{fig:ex-mfg-sysrisk-traj} displays three sample trajectories of $X$ and $Y$, obtained after training the neural networks for $Y_0, Z$ and $Z^0$, by simulating in a forward fashion the trajectories of $X$ and $Y$ using Monte Carlo samples and the same Euler-Maruyama scheme used in the numerical method. One can see that the approximation is better for $X$ than for $Y$, particularly towards the end of the time interval. This is probably because the BSDE is solved by guessing the initial point instead of starting from a terminal point, resulting in errors accumulated over time. Furthermore, \cite{carmona2019convergence2} shows that the results improve as the number of time steps, particles, and units in the neural network increase. In the numerical experiments presented here, we used the following parameters: $\sigma = 0.5, \rho = 0.5, q = 0.5, \epsilon = q^2 + 0.5 = 0.75, a = 1, c = 1.0$ and $T = 0.5$.

\begin{figure}[ht]
\centering
\subfloat[Trajectory of $X^i, i=1,2,3$]{
\includegraphics[width=.45\linewidth]{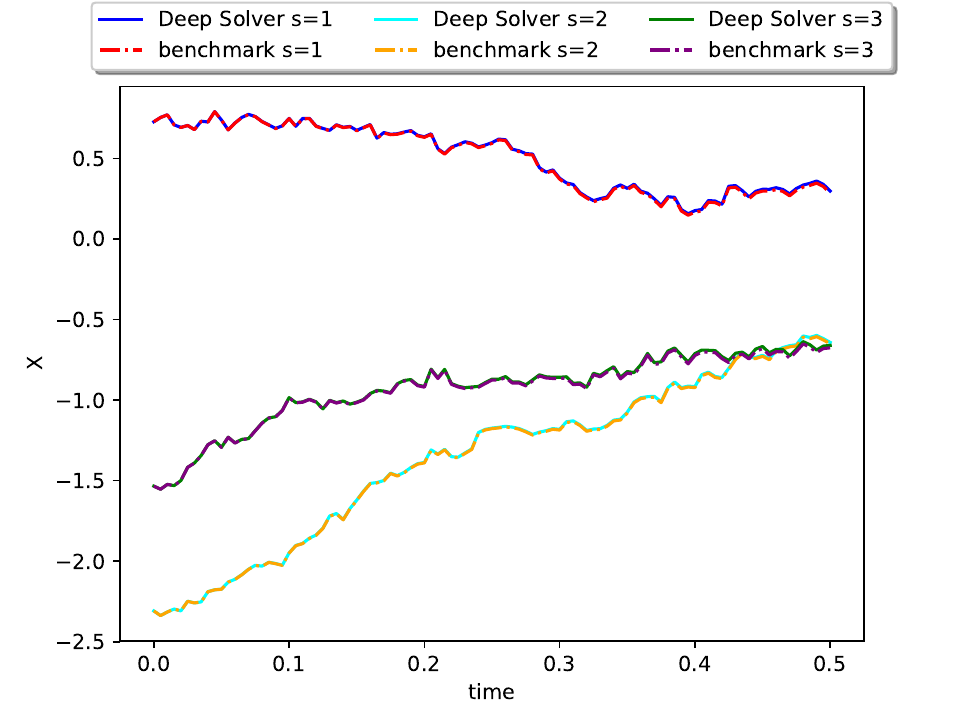}
}
\subfloat[Trajectory of $Y^i, i=1,2,3$] {
  \includegraphics[width=.45\linewidth]{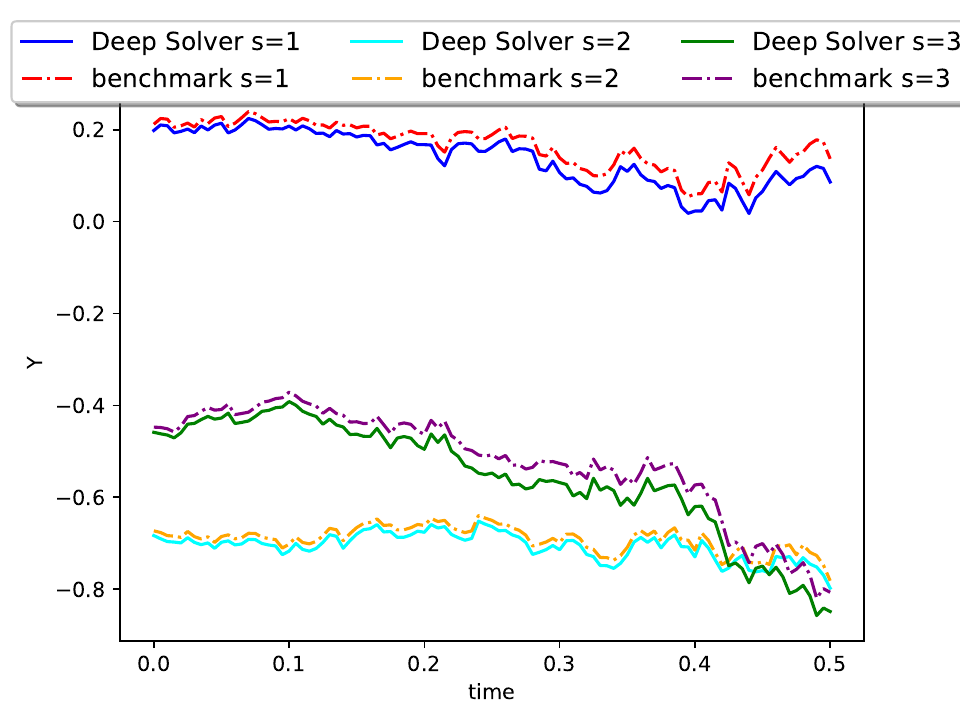}
  }
\caption{Systemic risk MFG example solved by the algorithm described in Section~\ref{sec4_MFGFBSDE}. Left: three sample trajectories of $X$ using the neural network approach ('Deep Solver' with full lines, in cyan, blue, and green) or using the analytical formula ('benchmark' with dashed lines, in orange, red and purple). Right: three sample trajectories of $Y$ (similar labels and colors). Note that the analytical formula satisfies the true terminal condition, whereas the solution computed by neural networks satisfies it only approximately since the trajectories are generated in a forward way starting from the learned initial condition. %
}
\label{fig:ex-mfg-sysrisk-traj}
\end{figure}

\subsubsection{PDE-based deep learning algorithms}

Besides direct parameterization of the controls and BSDE-based methods, it is also possible to adapt PDE-based methods to solve MFGs. In fact, such methods can be used in two different ways. First, one can use them to tackle the forward-backward PDE system characteriziting the distribution and the value function (see Section~\ref{sec:MFG-PDE-system}). One can also try to solve the Master equation (see Section~\ref{sec:background-master-eq}).

\paragraph{Deep learning for mean-field PDE systems}
\label{sec:MFPDE-deeplearning}

We now consider the PDE systems describing the equilibrium or social optimum in MFG or MFC, respectively. The Deep Galerkin Method (DGM) introduced in~\cite{SiSp:18} and reviewed in Section~\ref{sec:PDE} has been adapted to solve such PDE systems, see~\cite{al2022extensions,carmona2019convergence1,ruthotto2020machine,cao2020connecting,lin2020apac,AMSnotesLauriere}. We recall that the principle of the method is, for a single PDE, to replace the unknown function by a neural network and to optimize the parameters so as to minimize the residual of the PDE. 

For the sake of the presentation, we consider the MFG PDE system~\eqref{def:control-HJB-MFG}--\eqref{eq:MFG-general-KFP}. In line with the DGM method described in section~\ref{sec:PDE}, we proceed as follows. First, the MFG PDE system is rewritten as a minimization problem over the pair consisting of the density and the value function. The loss function is the sum of the two PDE residuals, as well as penalization terms for the initial and terminal conditions. Instead of the whole state space $\RR^d$, we focus on a compact subset $\tilde\dom \subset \RR^d$. If needed, extra penalization terms taking into account the boundary conditions can be added in the loss function. 
To be specific, we introduce the following loss function
\begin{equation}
\label{eq:loss-DGM-total}
	L(\mu,u) = L^\text{(KFP)}(\mu,u) + L^\text{(HJB)}(\mu,u),
\end{equation}
which is composed of one term for each PDE of the MFG system~\eqref{def:control-HJB-MFG}--\eqref{eq:MFG-general-KFP}. Each term is itself split into two terms: one for the residual inside the domain and one for the initial or terminal condition. The KFP loss function is
\begin{align}
\label{eq:loss-DGM-KFP}
	 L^\text{(KFP)}(\mu,u)
	 &= C^\text{(KFP)} \left\| 
        \displaystyle \partial_t \mu  - \sum_{i,j} \frac{\partial^2}{\partial_{x_i}\partial_{x_j}}\left( D_{i,j}\mu\right) + \diver\left( \mu b\right) \right\|^2_{L^2([0,T] \times \tilde\dom)} + C^\text{(KFP)}_0 \left\|  \mu(0)  - \mu_0 \right\|^2_{L^2(\tilde\dom)},
\end{align}
with $\nu_t = \mu_t \circ (I_d, \alpha(t,\cdot))^{-1}$ where $\ctrl(t,x) = \ctrl(t, x, \nu_t, \grad_x u(t,x), \Hess_x u(t,x))$, and $D$ and $b$ are defined as: 
\begin{align}
   D(t,x) = \frac{1}{2} \sigma(t, x, \nu_t, \alpha(t,x))\sigma(t,x, \nu_t, \alpha(t,x))\transpose, 
   \quad 
   b(t,x) = b(t,x,\nu_t,\ctrl(t,x)).
\end{align} 
The HJB loss function is
\begin{equation}
\label{eq:loss-DGM-HJB}
	 L^\text{(HJB)}(\mu,u)
	= C^\text{(HJB)} \left\| \partial_t u   + \min_{\alpha \in \mc{A}} H(\cdot,\cdot, \nu, \grad_x u, \Hess_x u, \alpha) \right\|^2_{L^2([0,T] \times \tilde\dom)}
	+ C^\text{(HJB)}_T \left\| u(T) - g(\cdot,  \mu(T) ) \right\|^2_{L^2(\tilde\dom)},
\end{equation}
with $H$ defined by~\eqref{def:control-H-MFG}. The weights $C^\text{(KFP)}, C^\text{(KFP)}_0, C^\text{(HJB)},$ and $C^\text{(HJB)}_T$ are positive constants that are used to tune the importance of each component relatively to the other components. If $(\mu,u)$ is a smooth enough solution to the PDE system~\eqref{def:control-HJB-MFG}--\eqref{eq:MFG-general-KFP}, then $L(\mu,u) = 0$. From here, the same strategy as in the DGM can be applied: one can look for an approximate solution using a class of parameterized functions for $\mu$ and $u$, replace the $L^2$ norms by integrals, and use samples to get Monte Carlo estimates; see~\cite{SiSp:18} and Section~\ref{sec:PDE} for more details.

\begin{remark} 
The same ideas can be applied to a variety of settings such as ergodic MFG~\cite{carmona2019convergence1}, non-separable Hamiltonians~\cite{ASSOULI2023113802} or finite-state MFGs~\cite{luo2023deep}. 
Each case may require some adjustments. For instance, in the case of ergodic MFGs, the initial and terminal conditions are replaced by normalization conditions; see~\cite{LaLi:2007}. Furthermore, if the PDE system was initially posed on a bounded domain and the solution had to satisfy boundary conditions, then these extra conditions could be dealt with by adding more penalty terms as in~\cite{carmona2019convergence1} or changing the architecture of the neural networks~\cite{cao2020connecting}. 
\end{remark}

\vskip 6pt

\paragraph*{Numerical illustration: a mean-field model of optimal execution.}

We now present an example based on a model of optimal execution. This model is similar to the one studied in Subsection \ref{sec:directMethod}. We consider a population of traders in which each trader wants to liquidate $Q_0$ shares of a given stock by a fixed time horizon $T$. At time $t\in[0,T]$, we denote by $S_t$ the price of the stock, by $Q_t$ the inventory ({\it i.e.}, number of shares) held by the representative trader, and by $X_t$ their wealth. These state variables are subject to the following dynamics 
\begin{equation*}
\begin{cases}
    dS_t = \gamma \bar{\mu}_t  dt + \sigma \ud W_t,\\
    dQ_t = \ctrl_t \ud t,\\
    dX_t = -\ctrl_t(S_t+\kappa \ctrl_t)\ud t.
\end{cases}
\end{equation*}
The evolution of the price $S$ is stochastic, representing that it can not be predicted with certainty. The randomness, scaled by $\sigma$, comes in through a standard Wiener process $W$. Furthermore, the drift of $S$ captures the permanent price impact $\gamma\bar{\mu}_t$ at time $t$. Here $\gamma>0$ is a multiplicative constant and $\bar{\mu}_t$ is the aggregate trading rate of all the traders.
The control $\ctrl_t$ at time $t$ corresponds to the individual rate of trading of the representative trader. Last, $\kappa>0$ is a constant that represents a quadratic transaction cost.

We assume that the representative agent tries
to maximize the following quantity, in which the first two terms reflect their payoff while the last two terms capture their risk aversion
$$
    \mathbb{E} \left[ X_T + Q_TS_T - A |Q_T|^2 - \phi \int_0^T |Q_t|^2 \ud t \right].
$$
 The constants $\phi>0$ and $A>0$ give weights to penalties for holding inventory through time and at the terminal time, respectively. 
 
 \begin{remark}
 Except for the fact that $\bar{\mu}_t$ is here endogenous, this is the model considered in \cite{MR3500455}, to which a deep learning method has been applied in~\cite{leal2020learning} to approximate the optimal control on real data.
In contrast with the model studied in Subsection \ref{sec:directMethod}, the model considered here is not linear-quadratic and the inventory is not directly subject to random shocks. We refer the interested reader to \cite{CarmonaWebster} and  \cite{CarmonaLeal} for more details and variants of these models.
\end{remark}

Although this problem is formulated with three state variables, we can actually reduce the complexity of the problem in the following way. When $(\bar{\mu}_t)_{0\le t\le T}$ is given, the optimal control of the representative agent can be found by solving an HJB equation. Following~ \cite{MR3500455}, the value function $V(t,x,s,q)$ can be decomposed as
 $V(t,x,s,q) = x + qs + v(t,q)$ for some function $v$ which is a solution to
$$
    -\gamma\bar{\mu} q = \partial_t v - \phi q^2 + \sup_\ctrl \{\ctrl \partial_q v - \kappa \ctrl^2\},
$$
with terminal condition $v(T,q) = - A q^2$. The maximizer in the supremum leads to the optimal control, which can be expressed as: $\ctrl^*_t(q) = \frac{\partial_q v(t,q)}{2\kappa}$. Based on this and going back to the consistency condition yields that, at equilibrium, the aggregate trading rate is
$$
    \bar{\mu}_t = \int \ctrl^*_t(q) \mu(t,dq) = \int \frac{\partial_q v(t,q)}{2\kappa} \mu(t, dq),
$$
 where $\mu(t,\cdot)$ is the distribution of inventories at time $t$ satisfying the KFP PDE:
$$
    \partial_t \mu + \partial_q\left( \mu  \frac{\partial_q v(t,q)}{2\kappa}\right) = 0, t \ge 0, \qquad \mu(0,\cdot) = \mu_ 0.
$$
As a consequence, the equilibrium solution of the MFG satisfies
\begin{equation}
\label{eq:CL-PDE-reduced}
\left\{
\begin{aligned}
    &\quad -\gamma\bar{\mu} q = \partial_t v - \phi q^2 + \frac{|\partial_q v(t,q)|^2}{4\kappa},
    \\
    &\quad \partial_t \mu + \partial_q\left( \mu  \frac{\partial_q v(t,q)}{2\kappa}\right) = 0,
    \\
    &\quad \bar{\mu}_t = \int \frac{\partial_q v(t,q)}{2\kappa} \mu(t, dq),
    \\
    &\quad \mu(0,\cdot) = \mu_ 0, v(T,q) = - A q^2. 
\end{aligned}
\right.
\end{equation}
The mean-field coupling between the two equations is non-local since it involves $\bar{\mu}_t$, and it combines the population distribution with the HJB solution. The value function, the associated control, and the mean of the distribution can be computed by solving a system of ODEs, which provides a benchmark to test numerical methods. We refer to \cite{cardaliaguet2017mean} for more details.

Moreover, \cite{al2022extensions} proposed a further change of variables to simplify the numerical computations and then used the DGM to approximate the solution of the transformed PDE system. Here, to simplify the presentation, we stick to the above PDE system~\eqref{eq:CL-PDE-reduced} and solve it directly using DGM. The initial and terminal conditions are imposed by penalization. For the non-local term $\bar{\mu}_t$, we use Monte Carlo samples to estimate the integral. In the implementation, we used the following values for the parameters: $T=1$, $\sigma=0.3$, $A = 1$, $\phi = 1$, $\kappa = 1$, $\gamma = 1$, and a Gaussian initial distribution with mean $4$ and variance $0.3$. To ensure that the neural network for the distribution always outputs positive values, we used on the last layer the exponential function as an activation function.   Figure~\ref{fig:ex-mfg-crowd-trade-distrib} shows the evolution of the distribution $m$. Figure~\ref{fig:ex-mfg-crowd-trade-value} shows the control obtained from the neural network approximating the HJB solution. The final distribution is concentrated near $0$, which is consistent with the intuition that the traders need to liquidate. Furthermore, the learned control coincides with the theoretical optimal control that can be computed by solving an ODE system~\cite{cardaliaguet2017mean}.

\begin{figure}[ht]
\subfloat[]{
  \includegraphics[width=0.45\linewidth]{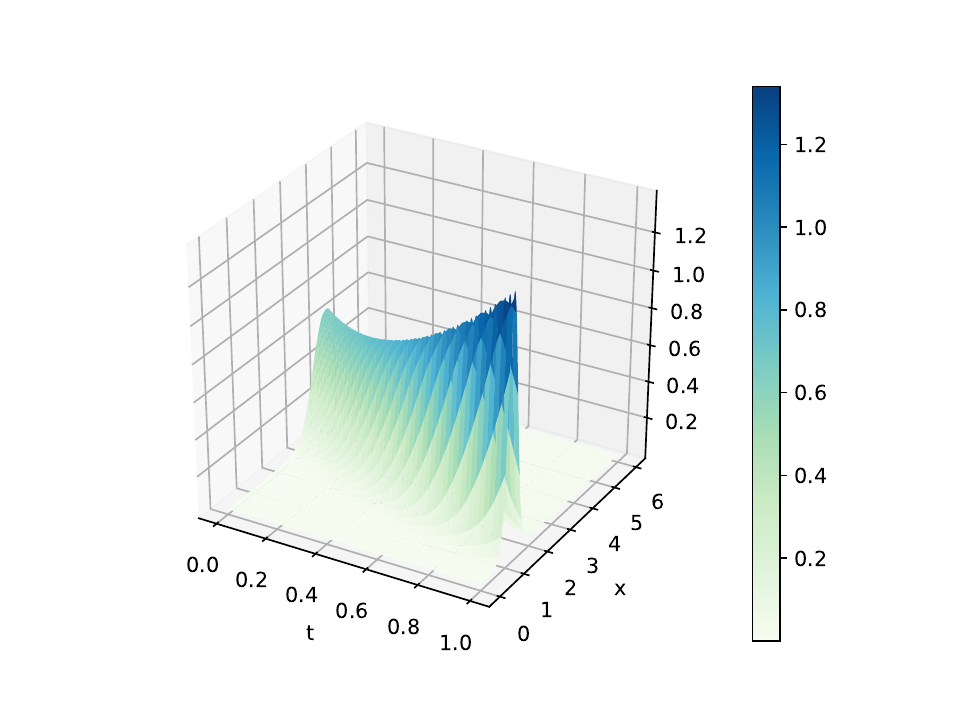}
  }
\subfloat[]{
\includegraphics[width=0.45\linewidth]{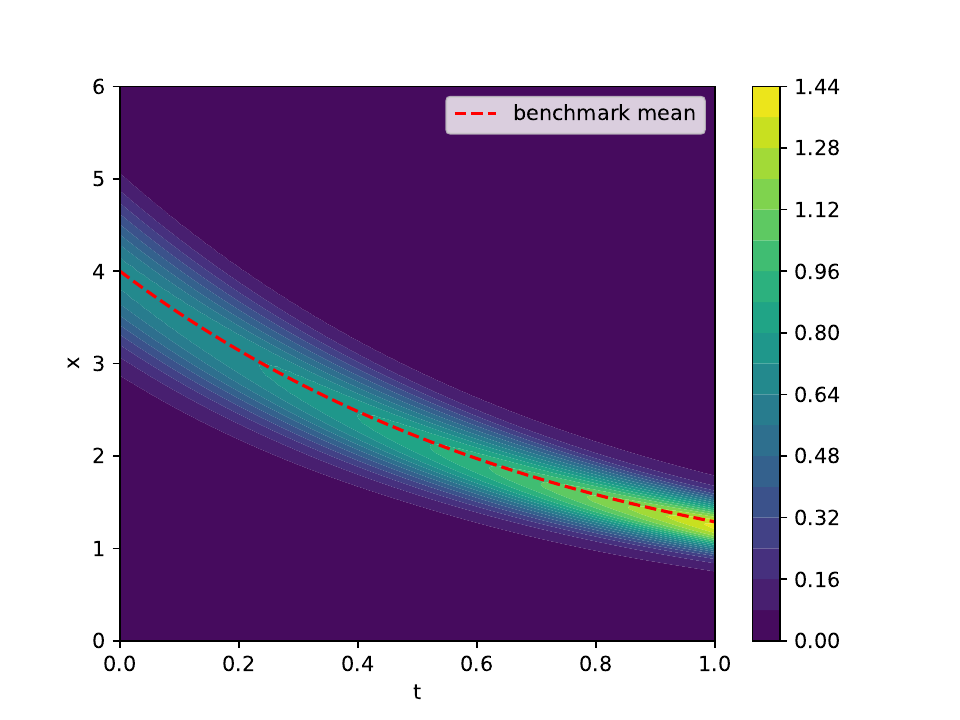}
}
\caption{Trade crowding MFG example in Section~\ref{sec:MFPDE-deeplearning} solved by DGM. Evolution of the distribution $m$. Left: surface with the horizontal axes representing time and space and the vertical axis representing the value of the density. Right: contour plot of the density with a dashed red line corresponding to the mean of the density computed by the semi-explicit formula. %
}
\label{fig:ex-mfg-crowd-trade-distrib}
\end{figure}

\begin{figure}[ht]
\subfloat[]{
\includegraphics[width=0.3\linewidth]{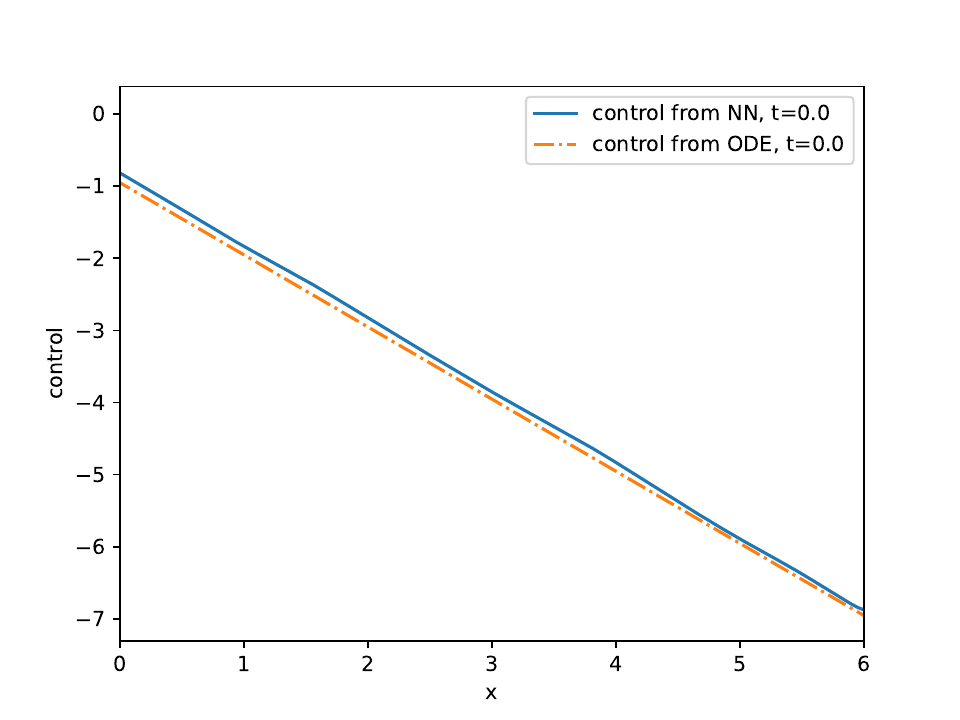}
}
\subfloat[]{
\includegraphics[width=0.3\linewidth]{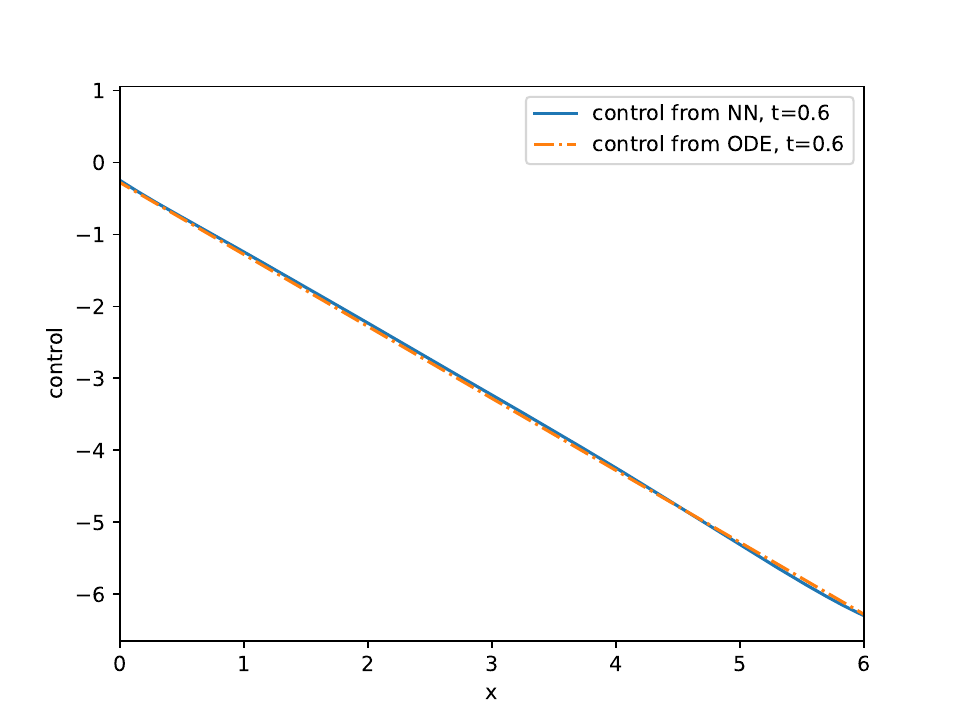}
}
\subfloat[]{
\includegraphics[width=0.3\linewidth]{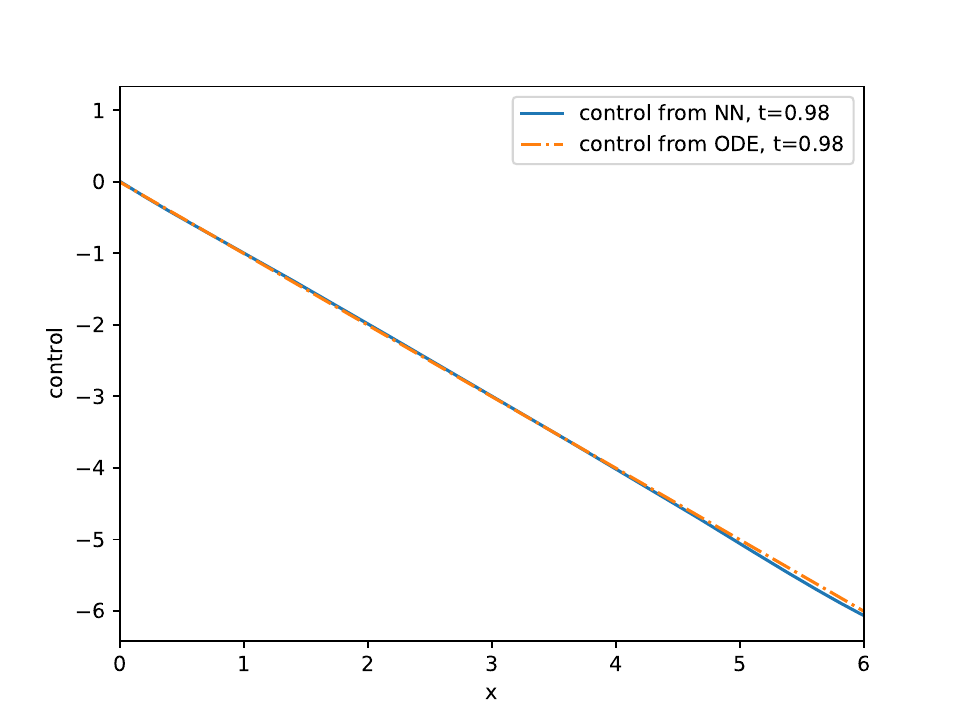}
}
\caption{Trade crowding MFG example in Section~\ref{sec:MFPDE-deeplearning} solved by DGM. Each plot corresponds to the control at a different time step: Optimal control $\ctrl^*$ (dashed line) and learned control (full line). %
}
\label{fig:ex-mfg-crowd-trade-value}
\end{figure}

\paragraph{Deep learning for mean-field master equation}
\label{sec:mastereq-deeplearning}

We now turn our attention to the question of solving the MFG master equation (Section~\ref{sec:background-master-eq}). Intuitively, the main motivation is to be able to approximate the value function of a representative player for any population distribution. This is in contrast with the methods presented above, which are based on controls that are fully decentralized in the sense that they are functions of the time and the state of the agent only. The fact that they do not depend on the population distribution is an advantage in that it simplifies the implementation, but it is also a limitation since the agent is not able to react to new distributions. For example, if the initial distribution is not known, the agent is not able to solve the forward equation and hence she is not able to anticipate the distribution at future times. The presence of common noise in the dynamics poses a similar challenge. For these reasons, being able to approximately solve the master equation is interesting for applications. When the state space is continuous, the distribution is an infinite dimensional object which is hard to approximate. For simplicity, we will thus focus here on a finite-state setting, in which case the distribution is simply a histogram. The convergence of finite-state MFGs to continuous-state MFGs has been studied for instance in~\cite{hadikhanloo2019finite}. Even though assuming the state space to be finite resolves the question of approximating the distribution, the master equation is posed on a high-dimensional space if the number of states is large. We will hence rely once again on neural networks to approximate the solution to this equation.

\smallskip
\noindent {\bf Master equation for finite state MFG.} We consider a finite-state MFG model based on the presentation of such models in~\cite[Section 7.2]{carmona2018probabilistic}. We consider a finite state space $\cE = \{e_1, \dots, e_d\}$ and an action space $\cA \subseteq \RR^\ctrldim$, which can be discrete or continuous. The states can be viewed as one-hot vectors, \textit{i.e.}, as the elements of a canonical basis of $\RR^d$. Then, the set of probability distributions on $\cE$ is the simplex $\{m \in \RR^d \,| \, \sum_{i=1}^d m_i = 1\}$, and we will sometimes write $m(x)$ instead of $m(\{x\})$. The running cost and the terminal cost are denoted by $f: \cE \times \cP(\cE) \times \cA \to \RR$ and  $g: \cE \times \cP(\cE) \to \RR$. The dynamics are given by a jump rate function denoted by $\lambda : \cE \times \cP(\cE) \times \cA \to \RR$. 
We denote by $\RR^\cE$ the set of functions from $\cE$ to $\RR$.

In this context, a finite state MFG equilibrium is a pair $(\hat{m}, \hat{\ctrl})$ with $\hat{m}: [0,T] \times \cE \to \RR$ and $\hat{\ctrl}: [0,T] \times \cE \to \cA$ such that
\begin{enumerate}
	\item $\hat{\ctrl}$ minimizes
\begin{align*}
	J^{MFG}_{\hat{m}}: \ctrl \mapsto  \EE \left[\int_0^T f(X_t^{\hat{m}, \ctrl}, \hat{m}(t,\cdot), \ctrl(t,X_t^{\hat{m}, \ctrl}) ) dt + g(X_T^{\hat{m}, \ctrl}, \hat{m}(T,\cdot)) \right],
\end{align*}
subject to: $X^{\hat{m}, \ctrl} = (X_t^{\hat{m}, \ctrl})_{t \ge 0}$ is a nonhomogeneous $\cE$-valued Markov chain with transition probabilities determined by the $Q$-matrix of rates $q^{\hat{m}, \ctrl}: [0,T] \times \cE \times \cE \to \RR$ given by
\begin{equation}
\label{master-eq-num:q-finite-MFG}
	q^{\hat{m}, \ctrl}(t, x,x') = \lambda(x, x', \hat{m}(t,\cdot), \ctrl(t,\cdot)), \qquad (t,x,x') \in [0,T] \times \cE \times \cE,
\end{equation}
and $X_0^{\hat{m}, \ctrl}$ has distribution with density $m_0$;
	\item For all $t \in [0,T]$, $\hat{m}(t,\cdot)$ is the law of $X_t^{\hat{m}, \hat{\ctrl}}$.
\end{enumerate}

The Hamiltonian of the problem is defined as
$$
	H(x, m, h) = \sup_{\ctrl \in \cA} -L(x, m, h, \ctrl)
$$ 
where $L: \cE \times \cP(\cE) \times \RR^\cE \times \cA \to \RR$ denotes the Lagrangian
$$
	L(x, m, h, \ctrl) = \sum_{x' \in \cE} \lambda(x, x', m, \ctrl) h(x') + f(x, m, \ctrl).
$$
Under suitable assumptions on the model, the supremum in the definition of $H$ admits a unique maximizer for every $(x, m, h) \in \cE \times \cP(\cE) \times \RR^\cE$, that we denote by
\begin{equation}
\label{master-eq-num:finiteMFG-ctrl-argmax}
	\ctrl^*(x, m, h) = \argmax_{\ctrl \in \cA} -L(x, m, h, \ctrl).
\end{equation}
The coefficients of the rates of the $Q$-matrix when using the optimal control are denoted by 
$$
	q^*(x,x',m,h) = \lambda\big(x, x', m, \ctrl^*(x,m, h)\big),
$$
where $q^*: \cE \times \cE \times \cP(\cE) \times \RR^\cE \to \RR$.

Similarly to the continuous setting (see Section~\ref{sec:mfg-theoretical-background}) the mean-field Nash equilibrium can be characterized using a forward-backward system of deterministic or stochastic equations. Using the deterministic approach, the optimal conditions take the form of an ODE system (instead of a PDE system as in the continuous space case). The system is composed of a forward ODE for the mean-field $m: [0,T] \times \cE \to \RR$ and a backward ODE for the value function $u: [0,T] \times \cE \to \RR$. Under suitable assumptions (see, \textit{e.g.}, \cite[section 7.2]{carmona2018probabilistic}), there is a unique MFG equilibrium $(\hat{m}, \hat{\ctrl})$, which is characterized by:
$$
	\hat{\ctrl}(t,x) = \ctrl^*(x, \hat m(t,\cdot), \hat u(t,\cdot)), 
$$  
where $\ctrl^*$ is defined by~\eqref{master-eq-num:finiteMFG-ctrl-argmax} and $(u,m)$ solves the forward-backward system
\begin{equation}
    \label{master-eq-num:ODE-system-finiteMFG}
    \begin{dcases}
        \displaystyle 0 
	   = - \partial_t \hat u(t, x) + H(x, \hat m(t,\cdot), \hat u(t,\cdot)),
	    \quad (t,x) \in [0,T) \times \cE,,
        \\
        0 
	    = \partial_t \hat m(t, x) - \sum_{x' \in \cE} \hat m(t, x') q^*(x', x, \hat m(t,\cdot), \hat u(t,\cdot)), 
	    \quad (t,x) \in (0,T] \times \cE,
	    \\
	    \hat u(T,x) = g(x, \hat m(T,\cdot)), \qquad \hat m(0,x) = m_0(x), 
	    \quad x \in \cE.
    \end{dcases}
\end{equation}

This ODE system can be solved using for example techniques discussed in previous sections for forward-backward PDE or SDE systems. However, this assumes that the initial distribution $m_0$ is known and when it is unknown, new techniques are required. We thus consider the master equation.

As in the continuous space case described in Section~\ref{sec:background-master-eq}, the solution to the master equation makes the dependence of $\hat u$ and $\hat m$ completely explicit. In the present discrete space setting, the master equation can be written as follows (see, \textit{e.g.}, \cite[section 7.2]{carmona2018probabilistic})
 \begin{equation}
 \label{master-eq-num:master-finiteMFG}
 	-\partial_t \cU(t,x,m) 
	+ H(x,m,\cU(t, \cdot, m))
	- \sum_{x' \in \cE} h^*(m ,\cU(t, \cdot, m))(x') \frac{\partial \cU(t, x ,m)}{\partial m(x')} = 0, 
 \end{equation}
for $(t,x,m) \in [0,T] \times \cE \times \cP(\cE)$, with the terminal condition $\cU(T,x,m) = g(x,m)$, for $(x,m) \in \cE \times \cP(\cE)$. The function $h^*: [0,T] \times \cP(\cE) \times \RR^\cE \times \cE \to \RR$ is defined as
$$
	h^*(m, u)(x') = \sum_{x \in \cE} \lambda(x, x', m, \ctrl^*(x, m, u)) m(x).
$$
Besides a simple representation of probability distributions, the fact that the state space is finite has another advantage: we do not need to involve the notions of derivative with respect to a measure discussed in Section~\ref{sec:background-master-eq}. Instead, we can rely on standard partial derivatives with respect to the finite-dimensional inputs of $\cU$. As a matter of fact, in the above equation, $\displaystyle \frac{\partial \cU(t, x ,m)}{\partial m(x')}$ denotes the standard partial derivative of $\RR^d \ni m \mapsto \cU(t,x,m)$ with respect to the coordinate corresponding to $x'$ (recall that $m$ is viewed as a vector of dimension $d$). The analog of~\eqref{eq:masterfield-to-u} in the continuous case is 
\begin{equation}
\label{eq:master-eq-to-hjb-sol-finitestate}
    \cU(t,x,\hat m(t)) = \hat u(t,x),
\end{equation}
where $\hat m = (\hat m(t))_t$ is the mean-field equilibrium distribution flow. Notice that both $\hat m$ and $\hat u$ implicitly depend on the initial distribution $m_0$, but $\cU$ does not.

The master equation~\eqref{master-eq-num:master-finiteMFG} is posed on a possibly high dimensional space since the number $d$ of states can be large. To numerically solve this equation, we can thus rely on deep learning methods for high-dimensional PDEs, such as the DGM introduced in~\cite{SiSp:18} and already discussed above in Sections~\ref{sec:PDE} and~\ref{sec:MFPDE-deeplearning}. 
For the sake of completeness, let us mention that this technique boils down to approximating $\cU$ by a neural network, say $\cU_\theta$ with parameters $\theta$, and using SGD to adjust the parameters $\theta$ such that the residual of~\eqref{master-eq-num:master-finiteMFG} is minimized and the terminal condition is satisfied. SGD as described in Algorithm~\ref{algo:SGD-generic} in Appendix~\ref{sec:SGD-var} is used, where a sample is $\xi = (t,x,m) \in [0,T] \times \cE \times \cP(\cE)$ and the loss function is
\begin{equation}
\label{master-eq-num:master-finiteMFG-residual}
	\mathfrak{L}(\cU_\theta, \xi) = \left|\partial_t \cU_\theta(t,x,m) - H(x,m,\cU_\theta(t,x,m))
	+ \sum_{x' \in \cE} h^*(m ,\cU_\theta(t,x,m))(x') \frac{\partial \cU_\theta(t,x,m)}{\partial m(x')}\right|^2.
\end{equation}

\paragraph*{Numerical illustration: A Cybersecurity model.} 
Here we present an example of the application of the above method. We consider the cybersecurity introduced in~\cite{MR3575619}; see also~\cite[Section 7.2.3]{carmona2018probabilistic}. Each player owns a computer and her goal is to avoid being infected by a malware. The state space is denoted by $\cE = \{DI, DS, UI, US\}$, which represents the four possible states in which a computer can be depending on its protection level -- defended (D) or undefended (U) -- and on its infection status  -- infected (I) or susceptible (S) of infection. The player can choose to switch its protection level between D and U. The change is not instantaneous so the player can only influence the transition rate. We represent by ``$1$'' the fact that the player has the intention to change its level of protection (be it from D to U or from U to D). On the other hand, ``$0$'' corresponds to the situation where the player does not try to change her protection level. So the set of possible actions is $\cA = \{0,1\}$. When the action is equal to $1$, the change of level of protection takes place at a rate  $\rho >0$. A computer in states DS or US might get infected either directly by a hacker or by getting the virus from an infected computer. We denote by $v_H q_{inf}^D$ (resp. $v_H q_{inf}^U$) the rate of infection from a hacker if the computer is defended (resp. undefended). We denote by $\beta_{UU}\mu(\{UI\})$ (resp. $\beta_{UD}\mu(\{UI\})$) the rate of infection from an undefended infected computer if the computer under consideration is undefended (resp. defended). Likewise, we denote by $\beta_{DU}\mu(\{DI\})$ (resp. $\beta_{DD}\mu(\{DI\})$) the rate of infection from a defended infected computer if the computer under consideration is undefended (resp. defended). Note that these rates involve the distribution since the probability of getting infected should increase with the number of infected computers in the rest of the population. Last, an infected computer can recover and switch to the susceptible state at rate $q_{rec}^D$ or $q_{rec}^U$ depending on whether it is defended or not.
These transition rates can be summarized in a matrix form: for $m \in \cP(\cE), a \in \cA$, 
$$
	\lambda(\cdot, \cdot, m, a) 
	= \left( \lambda(x, x', m, a) \right)_{x,x' \in \cE}
	= \begin{pmatrix}
	\dots 	& 		P^{m,a}_{DS \rightarrow DI}	&	 \rho a 	&	0
	\\
	q_{rec}^D 	& 	\dots 		&	 0	&	\rho a
	\\
	\rho a 	& 	0 		&	 \dots	&	P^{m,a}_{US \rightarrow UI}
	\\
	0	&	\rho a	&	q_{rec}^U	& \dots
	\end{pmatrix},
$$
where 
\begin{align*}
	&P^{m,a}_{DS \rightarrow DI} = v_H q_{inf}^D + \beta_{DD} m(\{DI\})  + \beta_{UD} m(\{UI\}) ,
	\\
	&P^{m,a}_{US \rightarrow UI} = v_H q_{inf}^U + \beta_{UU} m(\{UI\}) + \beta_{DU} m(\{DI\}).
\end{align*}
The dots ($\dots$) on each row stand for the value such that the sum of the coefficients on this row equals $0$. 

We assume that each player wants to avoid seeing her computer being infected, but protecting a computer costs some resources. So the running cost is of the form
$$
	f(t, x, \nu, \ctrl) = -\left[ k_D \indic_{\{DI, DS\}}(x) + k_I \indic_{\{DI, UI\}}(x)\right],
$$
where $k_D>0$ is a protection cost to be paid whenever the computer is defended, and $k_I>0$ is a penalty incurred if the computer is infected. We consider $g \equiv 0$ (no terminal cost).

By using the DGM, we train a neural network $\cU_\theta$ to approximate the solution $\cU$ to the master equation. Equation~\eqref{eq:master-eq-to-hjb-sol-finitestate} provides us with a way to check how accurate this approximation is:  we can fix an initial distribution, solve the forward-backward ODE system (which is easy given the initial condition), and then compare the value of the neural network $\cU_\theta$ evaluated along the equilibrium flow of distributions with the solution to the backward ODE for the value function. To be specific, for every $m_0$, we first compute the equilibrium value function $\hat u^{m_0}$ and the equilibrium flow of distributions $\hat m^{m_0})$. We then evaluate $\cU_\theta(t,x,\hat m^{m_0}(t,\cdot))$ for all $t \in [0,T]$ and check how close it is to $\hat u^{m_0}(t,x)$ for each of the four possible states $x$.  Figures~\ref{AMS-num-fig:finiteMFG-cyber-Master-m0-1}--\ref{AMS-num-fig:finiteMFG-cyber-Master-m0-3} show that the two curves (for each state) coincide for at least three different initial conditions. This means that, using the DGM, we managed to train a neural network that accurately represents the value function of a representative player for various distributions at once. In the numerical experiments, we used the following values for the parameters
\begin{align*}
&\beta_{UU} = 0.3,
\beta_{UD} = 0.4,
\beta_{DU} = 0.3,
\beta_{DD} = 0.4,
\qquad v_H = 0.2,
\lambda = 0.5,
\\
&q_{rec}^D = 0.1, 
q_{rec}^U = 0.65, 
q_{inf}^D = 0.4, 
q_{inf}^U = 0.3, 
\qquad k_D = 0.3, k_I = 0.5.
\end{align*}

\begin{figure}[ht]
\subfloat[]{		\includegraphics[width=.45\columnwidth]{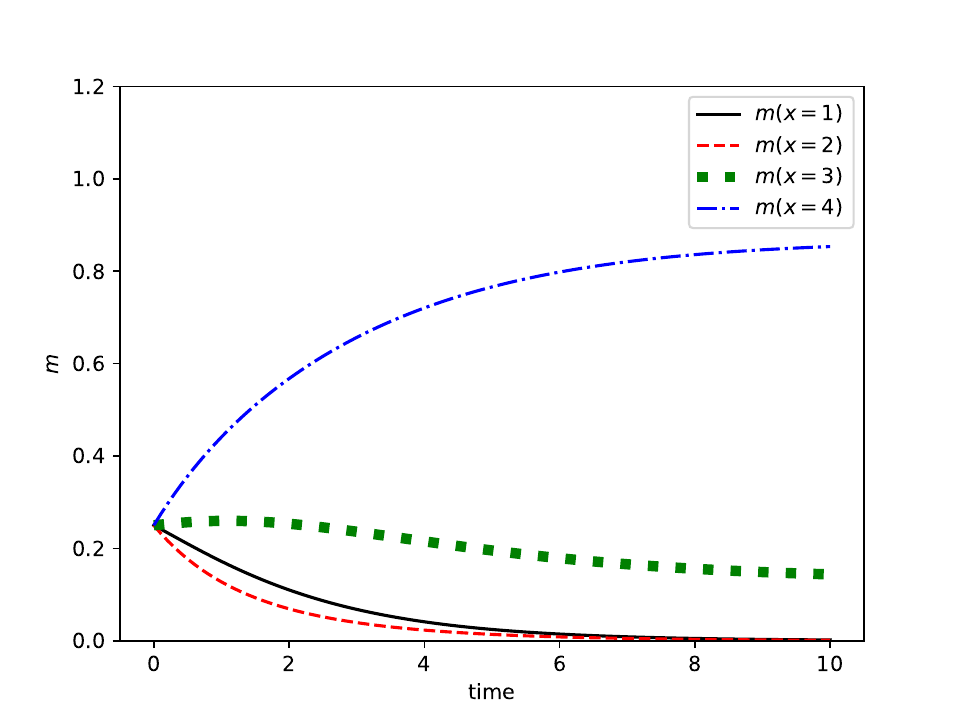}
}
\subfloat[]{
		\includegraphics[width=.45\columnwidth]{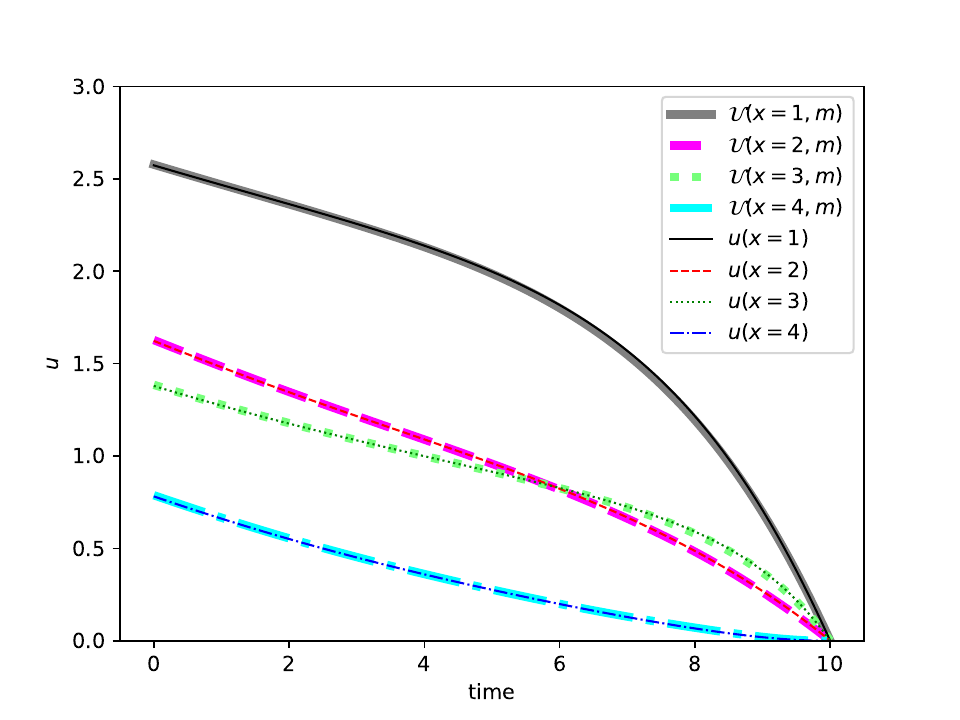}
}
	\caption{MFG Cybersecurity example in Section~\ref{sec:mastereq-deeplearning}. Test case 1: Evolution of the distribution $m^{m_0}$ (left) and the value function $u^{m_0}$ and $\cU(\cdot, \cdot, m^{m_0}(\cdot))$ (right) for $m_0 = (1/4, 1/4, 1/4, 1/4)$. First published in~\cite{AMSnotesLauriere} by the American Mathematical Society.} 
 	\label{AMS-num-fig:finiteMFG-cyber-Master-m0-1}
\end{figure}

\begin{figure}[ht]
\subfloat[]{
		\includegraphics[width=.45\columnwidth]{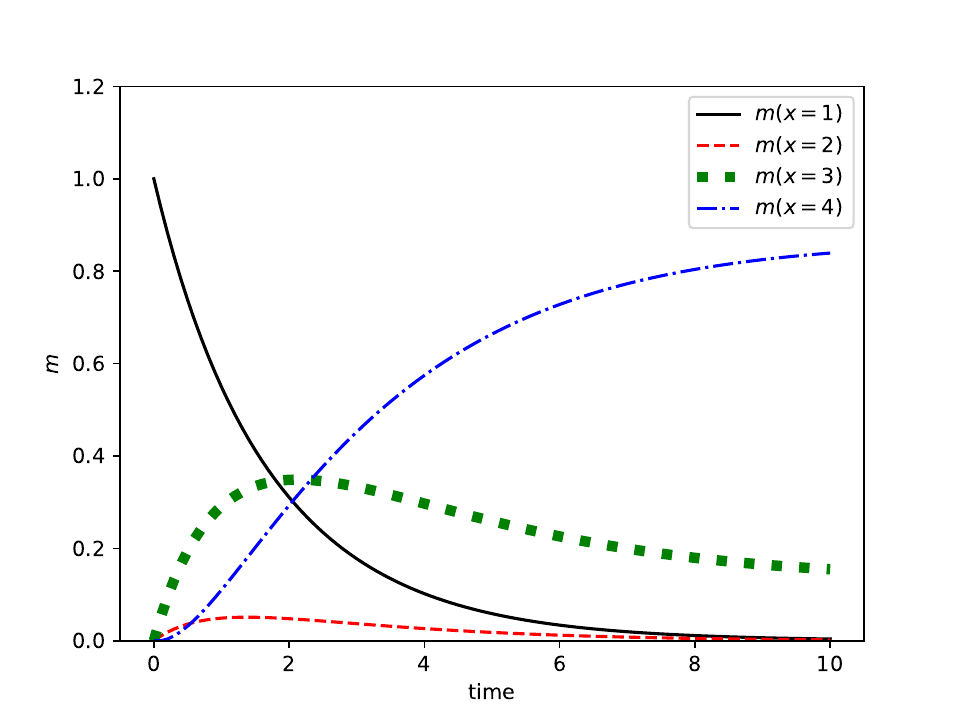}
}
\subfloat[]{
\includegraphics[width=.45\columnwidth]{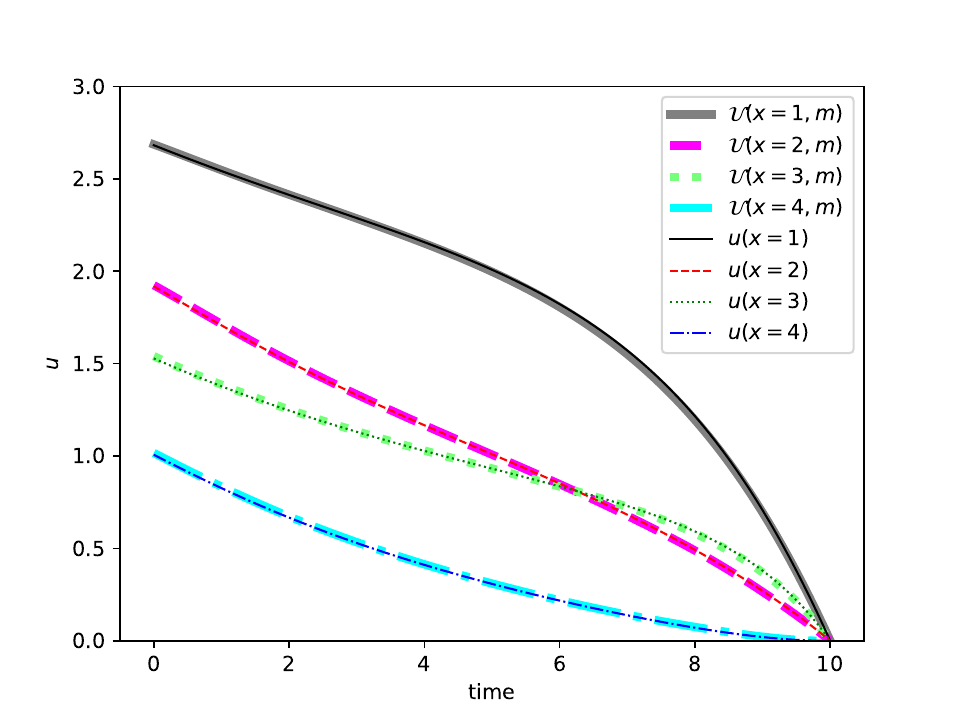}
}
	\caption{MFG Cybersecurity example in Section~\ref{sec:mastereq-deeplearning}. Test case 2: Evolution of the distribution $m^{m_0}$ (left) and the value function $u^{m_0}$ and $\cU(\cdot, \cdot, m^{m_0}(\cdot))$ (right) for  $m_0 = (1, 0, 0, 0)$. First published in~\cite{AMSnotesLauriere} by the American Mathematical Society.} 
 	\label{AMS-num-fig:finiteMFG-cyber-Master-m0-2}
\end{figure}

\begin{figure}[ht]
\subfloat[]{
\includegraphics[width=.45\columnwidth]{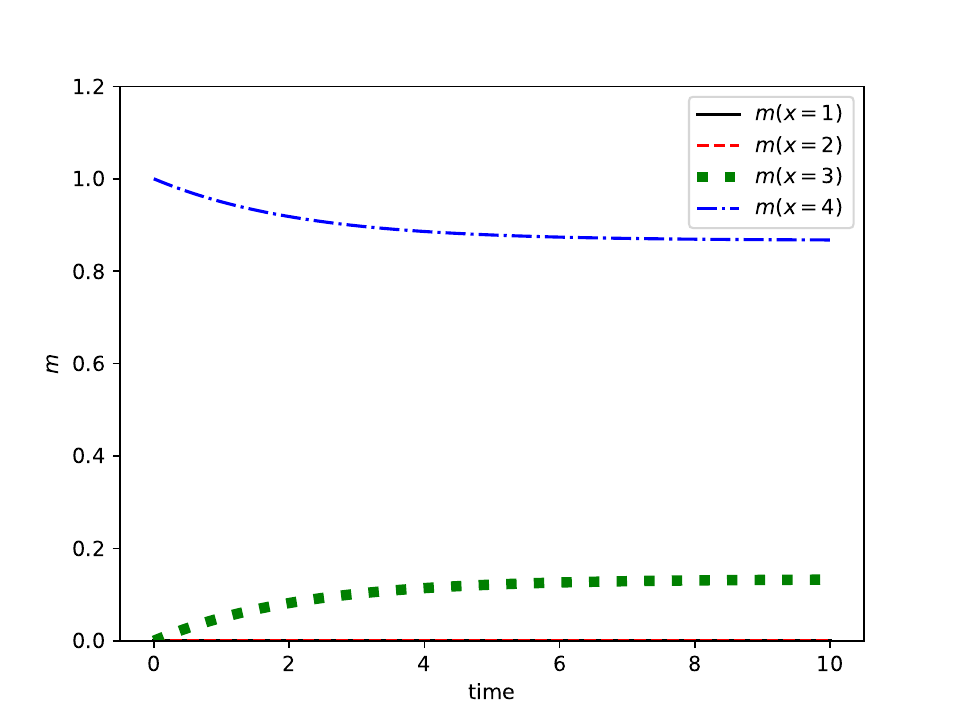}
}
\subfloat[]{
\includegraphics[width=.45\columnwidth]{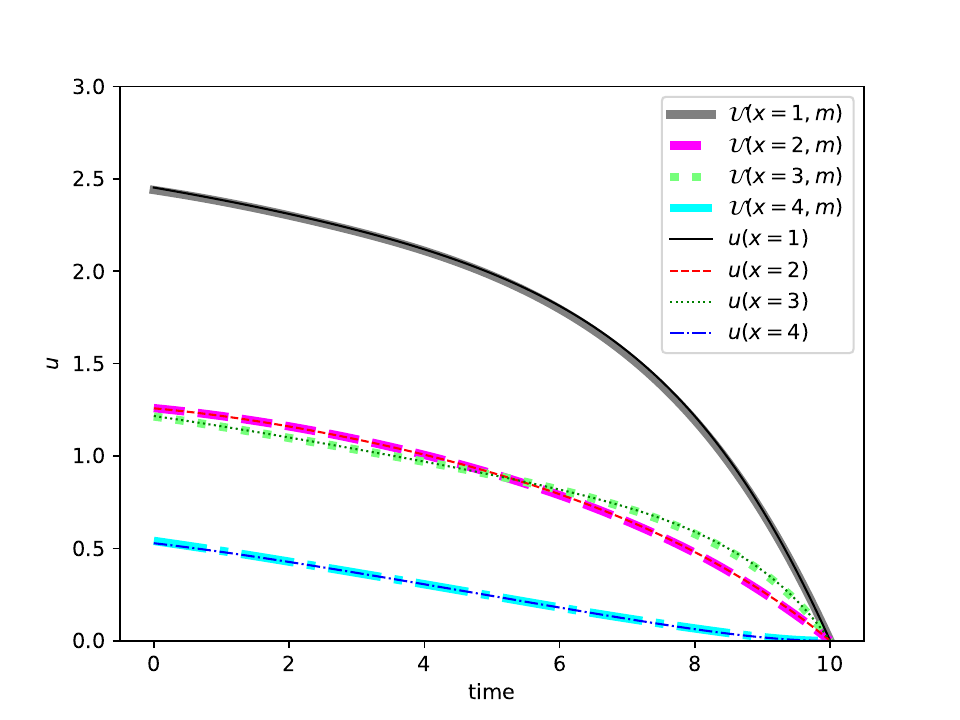}
}
	\caption{MFG Cybersecurity example in Section~\ref{sec:mastereq-deeplearning}. Test case 3: Evolution of the distribution $m^{m_0}$ (left) and the value function $u^{m_0}$ and $\cU(\cdot, \cdot, m^{m_0}(\cdot))$ (right) for  $m_0 = (0, 0, 0, 1)$. First published in~\cite{AMSnotesLauriere} by the American Mathematical Society.} 
 	\label{AMS-num-fig:finiteMFG-cyber-Master-m0-3}
\end{figure}

\section{Reinforcement Learning}\label{sec:RL}

All the previous methods rely, in one way or another, on the fact that the cost functions $f$ and $g$ as well as the drift $b$ and the volatility $\sigma$ (\emph{cf.} \eqref{def:control-Xt}--\eqref{def:control-cost}) are known. However, in many applications, coming up with a realistic and accurate model is a daunting task. It is sometimes impossible to guess the form of the dynamics, or the way the costs are incurred. This has motivated the development of so-called model-free methods. The reinforcement learning (RL) theory provides a framework for studying such problems.  Intuitively, an agent evolving in an environment can take actions and observe the consequences of her actions: the state of the environment (or her own state) changes, and a cost is incurred to the agent. The agent does not know how the new state and the cost are computed. The goal for the agent is then to learn an optimal behavior by trial and error.

Numerous algorithms have been developed under the topic of RL; see, {\it e.g.}, the surveys and books \cite{kaelbling1996reinforcement,Busoniu08,li2017deep,sutton2018reinforcement,hambly2021recent}. Most of them focus on RL itself, with state-of-the-art methods in single-agent or multi-agent problems and some provide theoretical guarantees of numerical performances. We aim to review its connections to stochastic control and games, as well as the mean-field setting. We shall start by discussing how problems in Section~\ref{sec:SCP} are formulated as single-agent RL\footnote{The terminology RL is from the perspective of artificial intelligence/computer science. In the operation research community, it is often called approximate dynamic programming (ADP) \cite{powell2007approximate}.}. Although we here focus on the traditional presentation of RL in discrete time, let us mention that a continuous-time stochastic optimal control viewpoint on RL has also been studied, see {\it e.g.}, \cite{munos2006policy,zhou2021actor} for policy gradient algorithms,  \cite{wang2020continuousmeanvar,wang2020reinforcement} for a mean-variance portfolio problem and for generic continuous time and space problems, and \cite{zhou2023policy,zhou2024solving} for convergence analysis of policy gradient methods. 
Furthermore, \cite{jia2021policy,jia2021policy2} studied policy-based and value function-based RL methods.
It has also been extended in several directions, such as variance reduction techniques~\cite{kobeissi2022variance}, risk-aware problems~\cite{jaimungal2022robustrl} and mean-field games~\cite{guo2022mfgentropyregu} and \cite{firoozi2022exploratory}.

Since the environment is unknown to the agent, a recurring question that needs to be addressed is the tradeoff between exploration and exploitation, whether in the single-agent, multiple-agent, or mean-field settings. Exploitation involves choosing actions that the agent believes will yield the highest immediate rewards based on her current knowledge, or equivalently; the agent follows the best-known strategies to maximize its short-term gains. Exploration, on the other hand, means taking actions that the agent has less information about, even if it might not result in the highest immediate rewards. The purpose of exploration is to gather more data and learn about the environment to improve the agent's long-term performance. Pure exploitation can lead to suboptimal decisions if the agent's initial knowledge is incomplete or incorrect. She may miss out on potentially better actions. Pure exploration, while informative, can be inefficient and may lead to delayed or reduced rewards, as the agent keeps trying unproven actions.

Effective RL algorithms aim to strike a balance between these opposing objectives. Alongside classical approaches like $\epsilon$-greedy and Upper-Confidence-Bound for exploration, recent developments have expanded our understanding of managing this trade-off to optimize an agent's learning and decision-making process in dynamic and uncertain environments. Notable contributions include \cite{wang2020reinforcement}, which elucidates the trade-off through entropy regularization from a continuous-time stochastic control perspective and offers theoretical support for Gaussian exploration, particularly in the context of a linear-quadratic regulator. \cite{guo2022mfgentropyregu} explores entropy regularization and devises a policy-gradient algorithm for exploration within the mean-field setting. Furthermore, \cite{delarue2021exploration} demonstrates that common noise can serve as an exploration mechanism for learning the solution of a mean-field game, through a linear-quadratic model.

\subsection{Reinforcement learning for stochastic control problems}
\label{sec:RLforSCP}
Recall the stochastic control problems studied in Section~\ref{sec:SCP}, see~\eqref{def:control-cost}--\eqref{def:control-Xt}. 
In this section, we will assume that the agent can not directly access $b, \sigma, f$ and $g$, but can observe the ``next step'' information given the current state and control. We consider the time discretized problem
\begin{align}
 &\checkX_{t_{n+1}} = \checkX_{t_n} + b(t_n, \checkX_{t_n}, \alpha_{t_n}) \Delta t + \sigma(t_n, \checkX_{t_n}, \alpha_{t_n}) \Delta \check W_{t_n}, \label{eq:rl-Xt-discrete} \\
&\min_{(\alpha_{t_n})_{n=0,\dots,N_T-1}} \EE\left[\sum_{n=0}^{N_T-1} f(t_n, \checkX_{t_n}, \alpha_{t_n}) \Delta t + g(\checkX_T)\right] \label{eq:rl-Xt-cost},
\end{align}
where
\begin{equation}
    0 = t_0 < t_1 < \ldots < t_{N_T} = T, \text{ with } t_n - t_{n-1} = \Delta t = T/N_T,
\end{equation}
is the temporal discretization on $[0,T]$ as before. By doing so, the system is Markovian, and can be viewed as a Markov decision process (MDP).

\subsubsection{Markov decision processes}\label{sec:MDP}

Problem \eqref{eq:rl-Xt-discrete}--\eqref{eq:rl-Xt-cost} can be recast as an MDP, which is a tuple $(\mc{X},\mc{A},p, f, g, N_T)$, where
\begin{itemize}
    \item $\mc{X}$ is the set of states called  the {state space};
    \item $\mc{A}$ is the set of actions called the {action space};
    \item $N_T<+\infty$ is the time horizon;
    \item $p: \mc{X} \times \mc{A} \times \{0,\Delta t, 2\Delta t, \dots,T\} \to \mathcal{P}(\mc{X})$ is the transition kernel, and $p(\cdot \vert x, a, t_n)$ is a probability density function; 
    \item $f:\{0,\Delta t, 2\Delta t, \dots,T\} \times \mc{X} \times \mc{A} \to \mathbb{R}$ is the one-step cost function, and $f(t_n, x, a)$ is the immediate cost at time $t_n$ at state $x$ due to action $a$;
    \item $g: \mc{X} \to \mathbb{R}$ is the terminal cost function, and $g(x)$ is the terminal cost at the final time $N_T$.
\end{itemize} 
 A large part of the RL literature focuses on the infinite horizon setting with discounted costs. Furthermore, the state space and action space are usually discrete (and they are often in fact finite), in which case  $p(x' \vert x, a, t_n) = \PP(X_{t_{n+1}} = x' \vert X_{t_n} = x, \alpha_{t_n} = a)$ is the probability to go to state $x'$ at time $t_{n+1}$ if at time $t_n$ the state $x$ and the action is $a$. 
 However, for the sake of consistency with the previous sections and the literature on optimal control, we stick to the finite horizon and continuous space setting in this section.

In model-free RL, the agent typically uses multiple episodes to learn the control that optimizes \eqref{eq:rl-Xt-discrete} with a simulator. In one episode of learning, the agent-environment interaction is as follows: Starting with $X_0 \in \mc{X}$, the agent chooses $\alpha_0 \in \mc{A}$, pays a one-step cost $f(0, X_0, \alpha_0)\Delta t$ and finds herself in a new state $X_{t_1}$; the process continues, forming a sequence:
\begin{equation}\label{eq:rl-sequence}
    X_0, \; \alpha_0, \; f(0, X_0, \alpha_0)\Delta t, \; X_{t_1}, \; \alpha_{t_1}, \; f(t_1, X_{t_1}, \; \alpha_{t_1})\Delta t, \; \ldots,\; X_T, \; g(X_T).    
\end{equation}

Under the Euler scheme \eqref{eq:rl-Xt-discrete}, given the state-action pair $(\check X_{t_n}, \alpha_{t_n}) = (x, a)$ at time $t_n$, $X_{t_{n+1}}$ follows a normal distribution $\mc{N}(x + b(t_n, x, a)\Delta t, \sigma^2(t_n, x, a)\Delta t)$.

In RL, there are four main components: policy, reward signal, value function, and optionally, a model of the environment. The MDP provides a mathematical framework to describe the agent-environment interface. A \emph{policy} $\pi: \{0,\Delta t, 2\Delta t, \dots,T\} \times \mc{X} \to \mc{P}(\mc{A})$ is a mapping from the state space to the probability space of the action space, and $\pi_t(a\vert x)$ describes the probability of choosing action $a$ at state $x$, which is in general random. 

The \emph{value function} associated to a specific policy $\pi$ is denoted by $V^\pi_t$ and defined as the expected cost when starting from $x$ at time $t$ and following $\pi$ thereafter, {\it i.e.},
\begin{equation}
    V_{t_n}^\pi(x) = \EE^\pi\left[\sum_{j=n}^{N_T-1} f(t_j, X_{t_j}, \alpha_{t_j})\Delta t + g(X_T) \vert X_{t_n} = x\right],
\end{equation}
where the superscript $\pi$ over the expectation means that, at each time step, the action is sampled according to $\pi$. 
Similarly, the \emph{action-value function} $Q^\pi_t$ associated to $\pi$ is defined as the expected cost when starting from $x$ at time $t$, taking the action $a$ and then following $\pi$, {\it i.e.},
\begin{equation}
    Q^\pi_{t_n}(x, a) = \EE^\pi\left[\sum_{j=n}^{N_T-1} f(t_j, X_{t_j}, \alpha_{t_j})\Delta t + g(X_T) \vert X_{t_n} = x, \; \alpha_{t_n} = a\right].
\end{equation}
Both functions satisfy the dynamic programming equations (also called Bellman equations),
\begin{align}
    &V^\pi_{t_n}(x) = \int_{a \in \mc{A}}\pi_{t_n}(a\vert x) \int_{x' \in \mc{X}} p(x' \vert x, a, t_n)[f(t_n, x, a)\Delta t + V^\pi_{t_{n+1}}(x')] \ud x' \ud a, \\
    & Q_{t_n}^\pi(x, a) = \int_{x' \in \mc{X}} p(x' \vert x, a, t_n) [f(t_n, x, a) \Delta t + V_{t_{n+1}}^\pi(x')] \ud x',
\end{align}
with terminal conditions $V_T^\pi(x) = Q_T^\pi(x, a) = g(x)$, where we have simplified the subscript $t_{N_T} = T$. The goal of RL is to identify the optimal $\pi^\ast = (\pi^\ast_t)_{t}$ that minimizes $V_t^\pi(x)$ for every $t$ and $x \in \mc{X}$. 

To this end, one also works with the \emph{optimal value function} $V^\ast_t(x) = \inf_\pi V_t^\pi(x)$ and the \emph{optimal action-value function} defined as $Q^\ast_t(x, a ) = \inf_\pi Q^\pi_t(x, a)$, which satisfy and the optimal Bellman equation reads,
\begin{align}
    &V^\ast_{t_n}(x) = \inf_{a \in \mc{A}} \int_{x' \in \mc{X}} p(x' \vert x, a, t_n)[f(t_n, x, a) \Delta t + V^\ast_{t_{n+1}}(x')] \ud x',
    \\
    &Q^\ast_{t_n}(x, a) = f(t_n, x, a) \Delta t +  \int_{x' \in \mc{X}} p(x' \vert x, a, t_n) \inf_{a' \in \mc{A}} Q^\ast_{t_{n+1}}(x', a') \ud x',
\end{align}
with terminal conditions $V_T^\ast(x) = Q_T^\ast(x, a) = g(x)$.

Model-free RL aims at computing $\pi^*$ without using the knowledge of the transition probability kernel $p$, and by instead relying on samples of transitions $X_{t_{n+1}} \sim p(x' \vert x, a, t_n)$. 
There are primarily two categories of learning methods: value-based methods and policy gradient methods. In the continuous time framework, the connection between policy evaluation in value-based methods and policy gradient methods has been developed in \cite{jia2021policy,jia2021policy2}.

\paragraph{Value-based methods}
For value-based methods, the workflow can be summarized as follows: starting with an arbitrary policy $\pi$, evaluate its value, improve the policy, and repeat until convergence:
\begin{equation}
    \pi_0 \to V^{\pi_0} \to \pi_1 \to V^{\pi_1} \to \ldots \pi_\ast \to V^{\ast}.
\end{equation}
The symbol $\pi_{\mt{k}} \to V^{\pi_{\mt{k}}}$ denotes a \emph{policy evaluation}, and the symbol $\pi_{\mt{k}} \to V^{\pi_{\mt{k}+1}}$ denotes a \emph{policy improvement}. 
Evaluating a given policy $\pi_{\mt{k}} \to V^{\pi_{\mt{k}}}$ exactly is not possible since we assume that $p(\cdot\vert x, a, t_n)$ is unknown. The Temporal-Difference (TD) learning remedies this issue by updating $V^\pi_{t_n}(x)$ with one sample drawn according to $X_{t_{n+1}} \sim p(x' \vert x, a, t_n)$,
$$
V^\pi_{t_n}(X_{t_n}) \leftarrow V^\pi_{t_n}(X_{t_n}) + \beta[f(t_n, X_{t_n}, \alpha_{t_n})\Delta t + V_{t_{n+1}}^\pi(X_{t_{n+1}}) - V_{t_n}^\pi(X_{t_n})],
$$
where $\beta>0$ is a learning rate. 
This is the simplest TD method, usually denoted by TD(0). To unify TD methods and MC methods, one can view the latter as updating $V^\pi_{t_n}$ using the entire sequence of observed cost from time $t_n$ until the end of the episode $T$. The $n$-step TD method lies in between and consists of simulating $n$ Monte Carlo samples to update $V^\pi$.

TD learning can also be applied to action-value functions, for example by using the update rule:
$$
    Q^{\pi}_{t_n}(X_{t_n}, \alpha_{t_n}) \leftarrow Q^{\pi}_{t_n}(X_{t_n}, \alpha_{t_n}) + \beta[f(t_n, X_{t_n}, \alpha_{t_n})\Delta t + Q^{\pi}_{t_{n+1}}(X_{t_{n+1}}, \alpha_{t_{n+1}}) - Q^{\pi}_{t_n}(X_{t_n}, \alpha_{t_n})], 
$$
where $X_{t_{n+1}}, \alpha_{t_{n+1}}$ are random samples from \eqref{eq:rl-Xt-discrete} and from $Q$ plus some randomization using for instance the $\eps$-greedy policy, which picks the currently optimal action with probability $1-\epsilon$ and, with probability $\epsilon$, picks an action uniformly at random. This approach is called SARSA.
Then the optimal action-value function $Q^\ast$ can be learned as follows: choose $\alpha_{t_n}$ according to $Q$ plus $\epsilon$-greedy for some exploration, then update $Q$ using SARSA. This method falls into the category of \emph{on-policy} algorithms since it evaluates or improves the policy that is used to make decisions. In fact, it uses an $\epsilon$-greedy policy to balance between learning an optimal behavior and behaving non-optimally for exploration, so it learns the value function for a sub-optimal policy that still explores. \emph{Off-policy} methods, on the contrary, use different policies for evaluation and data generation. \emph{Q-learning} may be the earliest well-known off-policy algorithm, which directly approximates $Q^\ast$ using the update rule:
$$
Q^{\pi}_{t_n}(X_{t_n}, \alpha_{t_n}) \leftarrow Q^{\pi}_{t_n}(X_{t_n}, \alpha_{t_n}) + \beta[f(t_n, X_{t_n}, \alpha_{t_n})\Delta t + \max_a Q^{\pi}_{t_{n+1}}(X_{t_{n+1}}, a) - Q^{\pi}_{t_n}(X_{t_n}, \alpha_{t_n})].
$$

\paragraph{Policy gradient methods}\label{sec:MDP-PGM}

This section describes some methods that aim at directly learning an optimal policy without deducing it from the value function. They use a parameterized class of policies. We denote by $\pi_t(a\vert x; \theta)$ the probability of taking action $a$ at state $x$ with parameter $\theta$. In practice, this can be a linear function $\theta\transpose\mathrm{f}(x, a)$ where $\mathrm{f}(x, a)$ is called feature vector, or a neural network taking $x$ as input and outputting a probability distribution over actions. Policy gradient methods update the policy parameter $\theta$ based on the gradient of some performance measure $L(\theta)$, with updates of the form
$$
    \theta \leftarrow \theta - \beta \widehat{\nabla J(\theta)},
$$
where $\widehat{\nabla L(\theta)}$ denotes an estimation of $\nabla L(\theta)$ based on Monte Carlo samples. A natural choice of $L(\theta)$ is the value function $V^{\pi_\theta}$ we aim to minimize. According to the policy gradient theorem,
$$
    \nabla V^{\pi_\theta}_{t_n}(x) = \EE_\pi \left[\int_{\mc{A}} Q^\pi_{t_n}(x, a)\nabla \pi_{t_n}(a\vert x; \theta) \ud a\right].
$$
Multiplying the first term by $\pi_{t_n}(a \vert x; \theta)$, dividing the second term by $\pi_{t_n}(a \vert x; \theta)$, replacing $a$ by a sample $\alpha_{t_n}$, and using $\EE_\pi[G_{t_n}\vert x, a] = Q^\pi_{t_n}(x, a)$ leads to the REINFORCE algorithm~\cite{williams1992simple},
$$
\theta \leftarrow \theta - \beta G_{t_n} \frac{\nabla_\theta \pi_{t_n}(\alpha_{t_{n}}\vert X_{t_n}; \theta)}{\pi_{t_n}(\alpha_{t_{n}}\vert X_{t_n}; \theta)},
$$
where $G_{t_n} = \sum_{n' = n+1}^{N_T-1} f(t_{n'}, \check X_{t_{n'}}, \alpha_{t_{n'}})\Delta t + g(\check X_T)$ denotes the cumulated cost from time $t_{n+1}$ to $T$.

\begin{remark}[Theoretical analysis]
Convergence of policy gradient has been studied in various settings. As for the model-based framework, LQ problems have attracted a particular interest since the optimal control can be written as a linear function of the state. Global convergence in the infinite horizon setting has been proved by Fazel et al. in~\cite[Theorems 7 and 9]{fazel2018global}. This result has been extended in various directions, such as the finite horizon setting in~\cite{hambly2021policy}, the neural setting \cite{wang2019neural}, problems with entropy regularization \cite{cen2022fast}, and MFC~\cite{carmona2019linear}, to cite just a few examples.
\end{remark}

With an additional parameterized value function $V_{t_n}(x; \theta')$, this leads to the actor-critic algorithm (see, {\it e.g.}, \cite[Section~13.5]{sutton2018reinforcement} or \cite{degris2012off}), 
\begin{align}
& \delta_{t_n} = f(t_n, X_{t_n}, \alpha_{t_n})\Delta t + V_{t_{n+1}}(X_{t_{n+1}}; \theta') - V_{t_n}(X_{t_n}; \theta'),\\
&\theta' \leftarrow \theta' - \beta' \delta_{t_n}\nabla_{\theta'} V_{t_n}(X_{t_n}; \theta'),\\
&\theta \leftarrow \theta - \beta \delta_{t_n} \nabla_\theta \ln \pi_{t_n}(\alpha_{t_n}\vert X_{t_n}; \theta).
\end{align}

Both REINFORCE and actor-critic methods mentioned above stochastically select an action $a$ when in state $x$ according to the parameter $\theta$, {\textit i.e.}, $a \sim \pi_{t_n}(\cdot \vert x, \theta)$. For some problems, it is more appropriate to look for a deterministic policy $\alpha_{t_n}(x; \theta) \in \mc{A}$. To ensure exploration, one can use an off-policy approach: a stochastic policy $\tilde\pi_{t_n}(a\vert x)$ is used to choose the action, and a deterministic policy $\alpha_{t_n}(x; \theta)$ is learned to approximate the optimal one.
An example of such methods is the deterministic policy gradient (DPG) \cite{pmlr-v32-silver14}, which is an off-policy actor-critic algorithm that learns a deterministic target policy $\alpha_{t_n}(x; \theta)$ from an exploratory behavior policy $\tilde\pi_{t_n}(a\vert x)$. In particular, a differentiable critic $Q(x, a; \theta')$ is used to approximate $Q^{\alpha(\cdot;\theta)}(x, a)$ and is updated via Q-learning: at each step, we  sample $\alpha_{t_n}$ from $\tilde\pi_{t_n}(a\vert x)$ and
\begin{align}
    & \delta_{t_n} = f(t_n, X_{t_n}, \alpha_{t_n})\Delta t + Q_{t_{n+1}}(X_{t_{n+1}}, \alpha_{t_{n+1}}(X_{t_{n+1}}; \theta); \theta') - Q_{t_n}(X_{t_n}, \alpha_{t_n}; \theta'),\\
&\theta' \leftarrow \theta' - \beta' \delta_{t_n}\nabla_{\theta'} Q_{t_n}(X_{t_n}, \alpha_{t_n}; \theta'),\\
&\theta \leftarrow \theta - \beta \delta_{t_n} \nabla_\theta \alpha_{t_n}(X_{t_n}; \theta)\nabla_a Q(X_{t_n}, \alpha_{t_n}; \theta')\vert_{a = \alpha_{t_n}(X_{t_n}; \theta)}.
\end{align}

When using neural networks to approximate $Q^{\alpha(\cdot;\theta)}$ and the deterministic policy $\alpha_{t_n}(x; \theta)$, one can use the Deep DPG (DDPG) algorithm \cite{lillicrap2015continuous}, which is based on the same intuition as DPG. For the sake of robustness it uses the ``replay-buffer'' idea borrowed from the Deep Q Network (DQN) algorithm, see~\cite{mnih2015human}: the network parameters are learned in mini-batches rather than online by using a replay buffer, so that correlation between samples are kept minimal. Another pair of networks $Q'(x, a; \hat \theta')$ and $\alpha'_{t_n}(x; \hat \theta)$ are copied from $Q(x, a; \theta')$ and $\alpha_{t_n}(x; \theta)$ for calculating the target value in order to improve the stability. At each step, an action $\alpha_{t_n}$ is sampled from $\alpha_{t_n}(X_{t_n}; \theta) + \mc{N}_{t_n}$ where $\mc{N}_{t}$ is a noise process for exploration; then the cost $f(t_n, X_{t_n}, \alpha_{t_n})\Delta t$ and the new state $X_{t_{n+1}}$ are observed and saved to the buffer. A mini-batch of $N$ transitions $(X_{t_n}, \alpha_{t_n}, f, X_{t_{n+1}})$ are sampled from the buffer, acting as supervised learning data for the critic $Q(x, a; \theta')$. The loss to be minimized is the mean-squared error of $Q_{t_n}(X_{t_n}, \alpha_{t_n}; \theta')$ and $f(t_n, X_{t_n}, \alpha_{t_n})\Delta t +  Q'_{t_{n+1}}(X_{t_{n+1}}, \alpha'_{t_{n+1}}(X_{t_{n+1}}; \hat \theta); \hat\theta')$. The actor network and both copies are updated via
\begin{align}
 & \theta \leftarrow \theta - \beta \frac{1}{N} \sum_{i}\nabla_a Q_{t_n}(x,a;\theta')\vert_{x = X_{t_n}^i, a = \alpha_{t_n}(X_{t_n}^i; \theta)} \nabla_\theta \alpha_{t_n}(X_{t_n}; \theta), \\
 & \hat \theta' \leftarrow \tau \theta' + (1-\tau) \hat \theta', \quad \hat \theta \leftarrow \tau \theta + (1-\tau) \hat \theta,
\end{align}
where the superscript $i$ indicates the $i^{th}$ sample from the mini-batch, and $\tau \ll 1$ is used to slowly track the learnt counterparts $\theta$ and $\theta'$.

\subsubsection{Mean-field MDP and reinforcement learning for mean-field control problems}
\label{sec5_RLforMFC}

We now consider the MFC setting discussed in Section~\ref{sec:directMethod} as an extension of standard OC and we present an RL framework for this setting. MFC can be viewed as an optimal control problem in which a ``state'' is a population configuration. However an ``action'' is not a finite-dimensional object but rather a function providing a control for every individual state. Intuitively, in discrete time, this yields an MDP of the form $(\mc{P}(\mc{X}), \mc{F}_{\mc{A}}, \bar{p}, \bar {f}, \bar{g}, N_T)$, where \begin{itemize}
    \item The state space is the set $\mc{P}(\mc{X})$ of probability measures on $\mc{X}$; 
    \item The action space $\mc{F}_{\mc{A}}$ is a suitable subset of $\mc{A}^{\mc{X}}$, the set of functions from $\mc{X}$ to $\mc{A}$;
    \item The transition kernel is given by  
$$
    \bar{p}:   \{t_0,t_1,\dots,T\} \times \mc{P}(\mc{X}) \times \mc{F}_\mc{A} \to \mc{P}(\mc{P}(\mc{X})), 
    \quad \bar{p}(\cdot|t, \mu, \bar a) = \delta_{\int p(\cdot|t, x, \mu, \bar{a}(x)) \mu(x) dx} \,,
$$
    meaning that with probability one, the new mean field state is given by one transition of the population distribution. Here $\mu$ represents a population distribution, $\bar a$ is an action at the population level, and $\bar{p}(\cdot|t, \mu, \bar a)$ is the distribution of the next population distribution, which is a Dirac mass at the next population distribution since there is no common noise in the present model;
    
    \item The running and terminal cost functions are given by
$$
    \bar{f}: \{t_0,t_1,\dots,T\} \times \mc{P}(\mc{X}) \times \mc{F}_{\mc{A}} \to\RR, \qquad \bar{f}(t, \mu, \bar{a}) = \int_{x} f(t, x, \mu, \bar{a}(x)) \mu(x) \ud x \,,
$$
and
$$
    \bar{g}: \mc{P}(\mc{X}) \to\RR, \qquad \bar{g}(\mu) = \int_{x} g(x, \mu) \mu(x) \ud x.
$$

\end{itemize}
Such MDPs have been referred to as mean field MDPs (MFMDP for short) in the literature~\cite{gast2012mean,carmona2019modelfree,gu2019dynamicmfc,gu2021mean,motte2019mean,bayraktar2022finite}. These MDPs can be rigorously studied using the tools developed for instance by Bertsekas and Shreve in~\cite{bertsekasshreve1996stochastic}. 
Since this problem fits in the framework of MDPs, one can directly apply RL methods in principle. For instance, the Q-function of the MDP naturally satisfies a dynamic programming principle; see~\cite{carmona2019modelfree,gu2019dynamicmfc,gu2021mean,motte2019mean}. Note that, if there is no common noise (as in the setting presented above), the evolution of the population distribution is purely deterministic. 

To implement RL methods for MFC, the main difficulties are related to handling the distribution and the class of controls. In particular, we note that
\begin{itemize}
    \item If $\mc{X}$ is finite, then the state of the MDP, namely $\mu$, is a finite-dimensional vector; if $\mc{A}$ is also finite, then $\mc{F}_{\mc{A}}$ can simply be taken as $\mc{A}^{\mc{X}}$, which is a finite set as well;
    \item If $\mc{X}$ is not finite, then $\mu$ is infinite-dimensional and likewise for the elements of $\mc{A}^{\mc{X}}$.
\end{itemize}

One simple approach is to discretize $\mc{P}(\mc{X})$ and $\mc{A}^{\mc{X}}$, and then use standard RL techniques for finite state, finite action MDPs, such as the ones described in Section~\ref{sec:MDP}. For instance tabular Q-learning has been used {\it e.g.} in~\cite{carmona2019modelfree,gu2021mean} in the first case above by identifying $\mc{P}(\mc{X})$ with the simplex $\Delta_{\mc{X}}$ in dimension $|\mc{X}|$ and by approximating the latter with an $\epsilon$-net. However, this approach does not scale well when the number of states is large or when $\mc{X}$ is continuous. In this case, one can use RL methods for continuous state space, such as deep RL methods, see for instance~\cite{carmona2019modelfree}.

\begin{remark}[Theoretical analysis]
    The convergence of Q-learning for MFMDP has been analyzed in~\cite{carmona2019modelfree} and~\cite{gu2021mean} using tabular or kernel-based methods respectively. The convergence of a policy gradient method for LQ MFC has been proved in~\cite{carmona2019linear} based the ideas of~\cite{fazel2018global}. 
\end{remark}

For the sake of illustration, we provide an example in a setting where $\mc{X}$ is finite. Let $d = |\mc{X}|$ be the number of states. As mentioned above, we view $\mc{P}(\mc{X})$ as the $d$-dimensional simplex $\Delta_{\mc{X}}$. In this case, the MFMDP is an MDP over a finite-dimensional continuous state space. To avoid discretizing the space, deep RL methods rely on neural networks to efficiently approximate the value function or the policy.

\paragraph*{Numerical illustration: A Cybersecurity model revisited.} 
We consider the cybersecurity model introduced in~\cite{MR3575619} (see also~\cite[Section 7.2.3]{carmona2018probabilistic}) and that we already discussed in Section~\ref{sec:mastereq-deeplearning}. We revisit this problem from the point of view of MFC, meaning that the players cooperate to jointly minimize the social cost. 

To be able to tackle this problem using RL, we discrete time using a mesh $\{t_n = n \Delta t, n = 0,1,2,\dots, N_T\}$ where $\Delta t = T / N_T>0$. The total cost for the whole population is
\begin{align*}
	J(\ctrl) 
	=  \sum_{n=0}^{N_T-1} \bar f(\mu_{t_n}, \ctrl(t_n, \cdot)) \Delta t,
\end{align*}
under the constraint that the evolution of distribution is given by
\begin{equation}
    \label{eq:cyber-mfc-dyn-mu}
	\mu_{t_{n+1}} 
	= \bar p(\mu_{t_n}, \ctrl(t_n,\cdot))
	= (\mu_{t_n})\transpose (I + P^{\ctrl(t_n,\cdot), \mu_{t_n}} \Delta t), \qquad n = 0, 1, \dots, N_T-1,
\end{equation}
with a given initial condition $\mu_0$. The population-wise cost function 
$\bar f: \mc{P}(\mc{X}) \times \mc{A}^{\mc{X}} \to \RR$ is defined based on the individual cost function $f$ by
$$
	\bar f(m, \ctrl) = \sum_{x\in\mc{X}} f(x, m, \ctrl(x) ) m(x), \qquad (m, \ctrl) \in \mc{P}(\mc{X}) \times \mc{A}^\mc{X}, 
$$
and $P^{\ctrl, m}$ denotes the matrix whose coefficients are given by
$$
	P^{\ctrl, m}(x', x) = \lambda(x', x, m, \ctrl(x')), \qquad (x',x,m,\ctrl) \in \mc{X} \times \mc{X} \times \mc{P}(\mc{X}) \times \mc{A}^\mc{X}.
$$
From this formulation, we see that the problem fits in the framework of MFMDPs, or MDPs with finite horizon and continuous space, the state being the distribution.

In~\cite{AMSnotesLauriere}, the solution is learned using tabular Q-learning after discretizing the simplex: replacing $\mc{P}(\mc{X})$ by an $\epsilon$-net with a finite number of distributions allows one to replace the MFMDP by a finite-state MFMDP on which tabular RL methods can be applied. This approach is convenient in that tabular methods typically have fewer hyperparameters and furthermore convergence results are easier to obtain. However, the main drawback is that such methods do not scale well to very large state space. In our case, discretizing the simplex requires a large number of points when the number of states increases.

Alternatively, the value function can be approximated directly on the simplex $\mc{P}(\mc{S})$, without any discretization. For example, we can replace the Q-function by a neural network and employ deep RL techniques to train the parameters. Here we follow the approach proposed in~\cite{carmona2019modelfree} and we focus on deterministic controls. The control and the value function are approximated by neural networks and trained using the DDPG method~\cite{lillicrap2015continuous}, which has been reviewed in Section~\ref{sec:MDP-PGM}. Since this method allows the control to take continuous values, we replace $A=\{0,1\}$ by $A =[0,1]$ (without changing the transition rate matrix), which amounts to letting the player choose the intensity with which she seeks to change her computer's level of protection. 

We aim at learning the solution for various distributions. To train the neural networks, we sample at each iteration a random initial distribution $\mu_0$ and generate a trajectory in the simplex by following the dynamics~\eqref{eq:cyber-mfc-dyn-mu}. Figure~\ref{fig:cyber-ddpg-cn-testing} displays the evolution of the population when using the learned control starting from five initial distributions of the testing set and one initial distribution of the training set. The testing set of initial distributions is: $\{(0.25,0.25,0.25,0.25),$ $(1, 0, 0, 0),$ $(0, 0, 0, 1),$ $(0.3, 0.1, 0.3, 0.1),$ $(0.5, 0.2, 0.2, 0.1)\}$. We see that, in this setting, the distribution always evolves towards a configuration in which there is no defended agents, and the proportion of undefended infected and undefended susceptible are roughly $0.43$ and $0.57$, respectively.

\begin{figure}[!htb]
\centering
\includegraphics[width=0.4\textwidth]%
{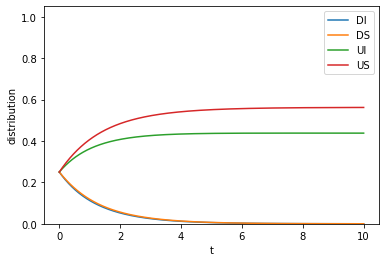}%
\includegraphics[width=0.4\textwidth]%
{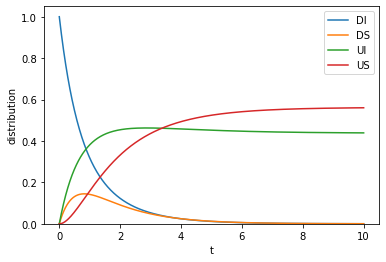}\\
\includegraphics[width=0.4\textwidth]%
{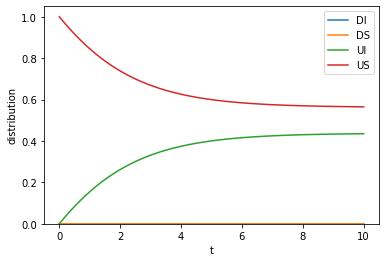}%
\includegraphics[width=0.4\textwidth]%
{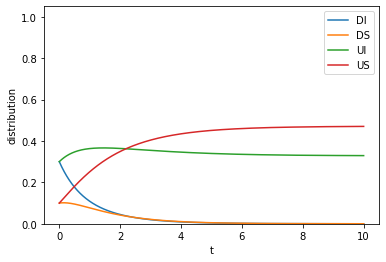}\\
\includegraphics[width=0.4\textwidth]%
{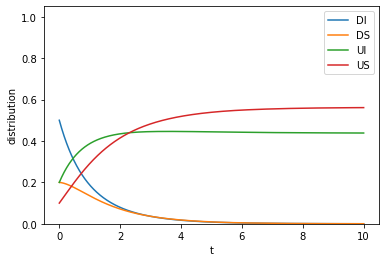}
  \caption{
      Cybersecurity MFC model solved with DDPG in Section~\ref{sec5_RLforMFC}: Evolution of the population distribution for five initial distributions. 
  }
  \label{fig:cyber-ddpg-cn-testing}
\end{figure}

\subsection{Reinforcement learning for stochastic differential games}
\subsubsection{Multi-agent reinforcement learning (MARL)}
Multi-agent reinforcement learning (MARL) studies reinforcement learning methods for multiple learners. The main difficulty is that, when several agents learn while interacting, from the point of view of each agent, the environment is non-stationary. Another issue is the question of scalability, which arises when the number of learners is very large. However, for a small number of agents, MARL has led to recent breakthrough results; see {\it e.g.} in autonomous driving ~\cite{shalev2016safe}, the game of Go~\cite{silver2016mastering}, or video games such as Star Craft~\cite{vinyals2019grandmaster}. 

Several viewpoints can be adopted. Relying on dynamical systems theory, one approach is to consider that each agent uses a learning algorithm, and to study the resulting behavior of the group of agents viewed as a system evolving in discrete or continuous time. Another approach, based on game theory and closer to the topics discussed in Section~\ref{sec:SDG}, is to look for notions of solutions such as Nash equilibria and to design algorithms that let the agents learn such solutions. A typical example is Nash Q-learning, in which every player runs their own version of Q-learning simultaneously with the other players. Each player tries to compute its optimal Q-function, but the optimal policy of player $i$ depends on the policies implemented by the other players. To be specific, consider an $N$-player game as in Section~\ref{sec:Nplayer} but now in discrete time. Note that the problem faced by player $i$ is not an MDP with state $X^i$ because the cost and dynamics of player $i$ depend on the other players. Assume the players use a strategy profile $\bm\pi = (\pi^1,\dots,\pi^N)$. Then the Q-function of player $i$ is: for $\bm{x} = (x^1,\dots,x^N)$ and $\bm{a} = (a^1,\dots,a^N)$,
\begin{equation}
    Q^{i,\bm\pi}_{t_n}(\bm{x}, \bm{a}) = \EE^{\bm\pi}\left[\sum_{j=n}^{N_T-1} f^i(t_j, \bm{X}_{t_j}, \bm{\alpha}_{t_j})\Delta t + g^i(\bm{X}_T) \Big\vert \bm{X}_{t_n} = \bm{x}, \; \bm{\alpha}_{t_n} = \bm{a}\right].
\end{equation}
Hu and Wellman proposed in~\cite{hu2003Nash} a version of Q-learning for (infinite horizon discounted) $N$-player games, called Nash Q-learning, and identified conditions under which this algorithm converges to a Nash equilibrium. The method can be adapted with deep neural networks, as done for instance in~\cite{casgrain2019deepnashq}.  We refer the interested reader to, {\it e.g.}, \cite{Busoniu08,tuyls2012multiagent,bloembergen2015evolutionary,lanctot2017unified,yang2020overview,zhang2021multi,gronauer2021multi} for more details on MARL. Recently, \cite{gu2021mean,gu2021mean2} also studied mean-field control RL in a decentralized way using cooperative MARL.

\subsubsection{Reinforcement learning for mean-field games}
\label{sec5_RLforMFG}

We now turn our attention to RL methods for MFG. As pointed out in Section~\ref{sec:MFG}, finding a mean-field Nash equilibrium boils down to (1) finding a control that is optimal for a representative infinitesimal player facing the equilibrium distribution flow, and (2) computing the induced distribution flow, which should match the equilibrium one. These two elements can be tackled alternatively, as described in Section~\ref{sec:Nplayer} in the N-player case and in Section~\ref{sec4_MFG_with_CN} in the mean-field case. The first part is a standard optimal control problem, which can thus be tackled using standard RL techniques, see Section~\ref{sec:MDP}. In this setting, we assume that the agent who is learning can repeat experiments of the following form: given the current state, the agent chooses an action (or a sequence of actions), and the environment returns the new state as well as the reward (or a sequence of states and rewards). In the representative player's MDP, the distribution enters as a parameter that influences the reward and dynamics, but is fixed when the player learns an optimal policy. During such experiments, we generally assume that the population distribution is fixed, and it is updated after a number of iterations; see {\it e.g.}~\cite{guo2019learning,EliePerolatLauriereGeistPietquin-2019_AFP-MFG}. Alternatively, we can assume that it is updated at every iteration but at a slow rate; see {\it e.g.}~\cite{subramanianpolicy,angiuli2020unified,xie2021learning}. Most of the literature thus far focuses on tabular methods. A few works have used deep RL methods to compute the best response. For example, DDPG has been used in~\cite{EliePerolatLauriereGeistPietquin-2019_AFP-MFG}, soft actor-critic (SAC) has been used for a flocking model in~\cite{perrin2021mfgflockrl}, while deep Q-learning or some variants of it has been used in~\cite{pmlr-cui21-approximately,perrin2021generalizationmfg,lauriere2022scalable}. Recently, several works have studied the advantages and the limitations brought by the regularization of the policy through penalization terms in the cost function~\cite{anahtarci2020mfgqregu,pmlr-cui21-approximately,guo2022mfgentropyregu}. We refer to~\cite{lauriere2022learningmfgsurvey} for a survey of learning algorithms and reinforcement learning methods to approximate MFG solutions.

\paragraph*{Numerical illustration: an example with explicit solution.}
For the sake of illustration, we consider an MFG model which admits an explicit solution in the continuous time ergodic setting. The model has been introduced and solved in~\cite{almulla2017two}. The MFG is defined as follows. The state space is the one-dimensional unit torus, {\it i.e.}, $\mathbb{T} = [0,1]$ with periodic boundary conditions. The action space is $\mathbb{R}$ (or in practice any bounded interval containing $[-2\pi,2\pi]$, which is the range of the equilibrium control). The drift function is 
$$
  b(x,m,a) = a.
$$
The running cost is
$$
	f(x,m,a) = \tilde f(x) + \frac{1}{2}|a|^2 + \log(m),
$$
where the first term is the following cost, which encodes spatial preferences for some regions of the domain
$$
	\tilde f(x)= -2 \pi^2 \sin(2 \pi x) + 2 \pi^2 \cos(2 \pi x)^2 - 2 \sin(2 \pi x).
$$
In the ergodic MFG setting, the objective of an infinitesimal representative player is to minimize
$$
    \lim_{T \to +\infty}\frac{1}{T} \EE\left[\int_0^T f(X_t, \mu_t(X_t), \alpha_t(X_t)) \ud t \right],
$$
where $X$ is controlled by $\alpha$. Here $\mu_t$ is assumed to have a density for every $t \ge 0$, and we identify it with its density. So $\mu_t(X_t)$ denotes the value of the density of $\mu_t$ at $X_t$.
The equilibrium control and the equilibrium mean-field distribution are respectively given by
$$
  a^*:x\mapsto 2 \pi \cos(2\pi x)\;\quad \mbox{and} \;\quad 
 \mu^*:x\mapsto  \frac{e^{2 \sin(2\pi x)}}{\int_{\mathbb{T}} e^{2 \sin(2\pi y)} \ud y}\;.
$$
We use fictitious play~\cite{cardaliaguet2015learning} combined with a deep RL algorithm to learn the best response at each iteration for the solution. The problem is in continuous state and action spaces and admits a deterministic equilibrium control. Hence, following~\cite{EliePerolatLauriereGeistPietquin-2019_AFP-MFG}, at each iteration, we solve the representative player's MDP using DDPG~\cite{lillicrap2015continuous} reviewed in Section~\ref{sec:MDP-PGM}.

The plots in Figure~\ref{fig:analytical-mfg-ddpg} are borrowed from~\cite{EliePerolatLauriereGeistPietquin-2019_AFP-MFG}. The left plot displays the $L^2$  distance between the true equilibrium control and the control learnt by the algorithm. The right plot shows the stationary distribution learnt by the algorithm, which is to be compared with the distribution described in~\cite{almulla2017two} for the ergodic problem. Although the two problems are slightly different (one being in the infinite horizon discounted setting and the other one in the ergodic setting),  we can see that the distribution has the same shape, for suitable choises of parameters.    
We refer to \cite{EliePerolatLauriereGeistPietquin-2019_AFP-MFG} for more details on the implementation and the choice of parameters and hyperparameters for the results shown in Figure~\ref{fig:analytical-mfg-ddpg}.
\begin{figure}[!htb]
\centering
\includegraphics[width=0.45\textwidth]%
{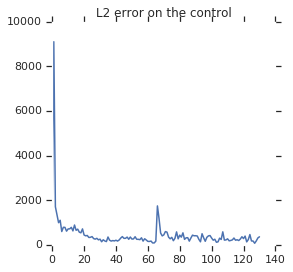}%
\includegraphics[width=0.45\textwidth]%
{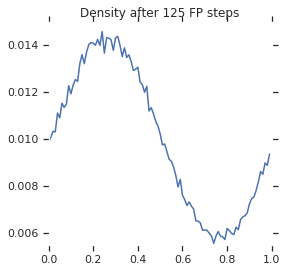}
  \caption{
      MFG described in Section~\ref{sec5_RLforMFG}, solved with fictitious play and DDPG. Left: $L^2$ error on the analytical control; right: stationary distribution. Results obtained by 125 iterations of fictitious play. 
  }
  \label{fig:analytical-mfg-ddpg}
\end{figure}

\section{Conclusion and Future Directions}\label{sec:conclusion}

This paper reviews recent developments in machine learning methods for stochastic optimal control and games, with a special focus on emerging applications of deep learning and reinforcement learning to these problems. Despite the rapidly growing number of recent works, many questions remain to be investigated further. We hope this survey will generate more interest in this topic and attract more researchers to work on it. Besides the material already reviewed in this survey, we outline a few research directions below.

First, the main goal of this survey was to provide an overview of existing methods, with a harmonized presentation of the types of problems studied in this field. However, the literature still lacks unified and thorough comparison of all existing methods on common benchmark problems. Indeed, most existing works use different assumptions for the their theoretical analysis and, on the numerical side, focus on illustrating the performance of one method on examples which are not necessarily the same as examples considered in other works. To understand better in which case each method is the most suitable, it would be important to provide a detailed comparison of the assumptions required for the analysis and to perform rigorous numerical comparisons on common problems.

Most of the methods presented here lack full analysis on the theoretical side. The mathematical foundations of deep learning are attracting growing interest, and recent results could help analyze the methods described in this paper. The main motivation underlying the use of deep networks is their ability to cope with the curse of dimensionality. However, rigorously phrasing and proving such a statement has only been done in particular cases. Analyzing the generalization capability of neural networks is typically done by splitting the analysis into several types of errors, such as approximation, estimation, and optimization errors. Bounds on the approximation and estimation errors can generally be obtained based on the regularity of the function to be approximated, which can be difficult in the context of differential games. Furthermore, bounding the optimization error is even more challenging since it involves not only the definition of the game but also the optimization algorithm. Due to these difficulties, estimating these errors remains an open question for most methods discussed in this survey.

From a practical viewpoint, an important question related to neural network-based methods is the choice of hyperparameters. The most obvious one is the architecture of the neural network. In many cases, a feedforward fully connected architecture provides good performances ({\it e.g.}, for deep BSDE, DBDP, Sig-DFP). However, in other cases ({\it e.g.}, DGM, RNN for problems with delay, as discussed in this survey), ad hoc architectures seem necessary to reach the best results. In any case, architectures undoubtedly play a crucial role in the performance of every deep learning method, and a careful design is, in general, what leads to state-of-the-art results in high dimensions. However, most deep learning methods for differential games presented in this survey have been limited to proof of concepts. As such, exploring more sophisticated architectures is a natural progression towards achieving better numerical performance. Once the neural network architecture is fixed, the next important step is to determine the hyperparameters of the optimization method, such as the initialization of network parameters, the learning rate, and the mini-batch size, which are crucial factors for ensuring fast convergence. However, finding a systematic rule for choosing these hyperparameters {\it a priori} remains a challenge. A common approach is to try several values and measure the empirical convergence speed on problems for which the solution is known, using either an analytical formula or another numerical method. This task is complex due to the interdependent influences of hyperparameters. For problems without benchmarks, finding suitable hyperparameters is even more challenging. To the best of our knowledge, the literature does not yet provide a comprehensive understanding of how to choose hyperparameters and measure algorithm performance without benchmark solutions. Although we did not discuss this aspect in the present survey, finding efficient heuristics for hyperparameter tuning is certainly an interesting direction. Another related question is how to assess the convergence of algorithms that compute Nash equilibria in games since the objective is to find a fixed point, not an optimizer.

A direction that has received little attention thus far regarding specific problems related to MFGs is numerical methods that can work even when a common noise affects the entire population. The difficulties that arise numerically are connected to the difficulty of solving such MFGs from a theoretical viewpoint. We have presented the Sig-DFP method to tackle MFGs with common noise, focusing on mean-field interactions through moments. Common noise appears in applications, for instance, in the form of aggregate shocks in macroeconomics. Therefore, it is worth developing further machine learning algorithms to deal with MFGs with common noise and general interactions. Currently, we lack efficient ways to parameterize, represent, and discretize probability measures defined on a continuous-state space.

Another aspect related to concrete applications of the methods presented in this survey pertains to the resources needed to train deep neural networks. For model-based methods and even more for model-free reinforcement learning methods, sophisticated models typically need a vast number of training episodes, leading to two challenges. First, as the model complexity grows, the massive computational cost required to learn the solution becomes prohibitive. Second, for real-world applications, Monte Carlo simulations will be replaced by real data, but we generally have much fewer data points than the number of samples used by most deep learning methods described in this survey. Therefore, it will be very interesting to design more sample-efficient methods and establish sharp estimates of their sample complexity.

Last but not least, to the best of our knowledge, the methods presented in this survey have only been applied to relatively simple models for academic research purposes. However, a significant motivation for the development of machine learning methods is to enable us to efficiently solve more realistic optimal control problems and games. We hope that this survey can contribute to fostering interactions between theoretical research and applied research communities, leading to concrete applications in real-world problems.

\section*{Acknowledgement}
R.H.~was partially supported by the NSF grant DMS-1953035, the Faculty Career Development Award, the Research Assistance Program Award, the Early Career Faculty Acceleration funding, and the Regents' Junior Faculty Fellowship at University of California, Santa Barbara. Some parts of the review paper have been used for teaching special topic graduate classes at the University of California, Santa Barbara, and R.H. appreciates all the feedback from the audience of these classes. R.H. and M.L. are grateful to all their co-authors of the papers mentioned in this review.

\appendix

\section{Deep Learning Tools}\label{sec:DLtools}

In this section, we briefly review neural networks and stochastic gradient descent, which are two of the main tools of modern machine learning. We refer to \cite{higham2019deep} for a more comprehensive mathematical introduction to deep learning for applied mathematicians and to~\cite{goodfellow2016deep} for more background on deep learning.

\subsection{Neural network architectures}
\label{subsec:NN-architectures}

We start by introducing the feedforward fully connected architecture, before discussing recurrent neural networks and long short-term memory networks.

\subsubsection{Feedforward fully connected neural networks} Feedforward neural networks (FNNs) are the most common type of neural networks. We denote by
\begin{align*}
	\mathbf{L}^\rho_{d_1, d_2} = 
	&\Big\{ \phi: \RR^{d_1} \to \RR^{d_2} \,\Big|\,  \exists (w, \beta) \in \RR^{d_2 \times d_1} \times \RR^{d_2}, \forall  i \in \{1,\dots,d_2\}, 
	\;\phi(x)_i = \rho\Big(\beta_i + \sum_{j=1}^{d_1} w_{i,j} x_j\Big) \Big\} 
\end{align*}
the set of layer functions with input dimension $d_1$, output dimension $d_2$, and activation function $\rho: \RR \to \RR$. Typical choices for $\rho$ are ReLU (positive part), identity, sigmoid, or hyperbolic tangent,
\begin{equation}\label{eq:activation}
    \rho_{\text{ReLU}}(x) = \max\{x,0\}, \quad \rho_{\text{Id}}(x) = x, \quad \rho_{\text{s}}(x) = \frac{1}{1+e^{-x}}, \quad \rho_{\tanh}(x) = \tanh(x).
\end{equation}
Building on this notation and denoting by $\circ$ the composition of functions, we define
\begin{align*}
	\bN^{\rho,\tilde\rho}_{d_0, \dots, d_{\ell+1}} 
	= 
	&\Big\{ \phi_\ell \circ \phi_{\ell-1} \circ \dots \circ \phi_0 \,\Big|\,  (\phi_i)_{i=0, \dots, \ell-1} \in \bigtimes_{i=0}^{i=\ell-1} \mathbf{L}^\rho_{d_i, d_{i+1}},   \phi_\ell \in \mathbf{L}^{\tilde\rho}_{d_{\ell}, d_{\ell+1}}  \Big\} \, 
	\notag
\end{align*}
 as the set of regression neural networks with $\ell$ hidden layers and one output layer, the activation function of the output layer being $\tilde\rho$. The number $\ell$ of hidden layers, the numbers $d_0$, $d_1$, $\cdots$ , $d_{\ell+1}$ of units per layer, and the activation functions, are the components of what is called the architecture of the network. Once it is fixed, the actual network function $\varphi\in \bN^{\rho,\tilde\rho}_{d_0, \dots, d_{\ell+1}} $ is determined by the remaining parameters
 $$
 \theta=(\beta^{(0)}, w^{(0)},\beta^{(1)}, w^{(1)},\cdots\cdots,\beta^{(\ell-1)}, w^{(\ell-1)},\beta^{(\ell)}, w^{(\ell)}),
 $$
defining the functions $\phi_0$, $\phi_1$, $\cdots$ , $\phi_{\ell-1}$ and $\phi_\ell$ respectively. Let us denote by $\Theta$ the set of values for such parameters. For each $\theta\in\Theta$, the function computed by the network will be denoted by $\varphi^\theta \in \bN^{\rho, \tilde \rho}_{d_0, \dots, d_{\ell+1}}$ when we want to stress the dependence on the parameters.

To alleviate the presentation, we will follow the convention of using vector and matrix notations, and here activation functions are implicitly applied coordinate-wisely. Then
$$
    \varphi^{\theta}(x) = \tilde\rho\left(\beta^{(\ell)} + w^{(\ell)} \rho\left( \beta^{(\ell-1)} + w^{(\ell-1)} \rho\left(  \dots \beta^{(0)} + w^{(0)} x\right)\right) \right).
$$

\subsubsection{Recurrent neural networks}\label{sec:RNNdetails} 

Although FNNS are universal approximators, they are not very suitable to handle path-dependent properties of the state process, which are important for instance when the stochastic control problem or the game has delay features. The idea of recurrent neural networks (RNNs)~\cite{rumelhart1986learning} is to make use of sequential information, and thus provide a natural framework for overcoming these issues. In fact, RNNs have already shown great success in, {\it e.g.}, natural language processing and handwriting recognition \cite{graves2013generating,graves2013speech,graves2009offline}. Many variants exist and below we shall focus on one such variant, but the generic architecture can be described as follows: the neural network takes two inputs, $x$ and $h$, and produces two outputs, $y$ and $h'$, as follows:
\begin{align*}
    h' &= \rho\left(\beta^{(1)} + w^{(1,1)} h + w^{(1,2)} x \right),
    \\
    y &= \tilde\rho\left(\beta^{(2)} + w^{(2)} h' \right),
\end{align*}
where $\rho,\tilde \rho$ are two activation functions, and the parameters of the neural network are vectors $\beta^{(1)},\beta^{(2)}$ of suitable sizes, and matrices $w^{(1,1)}, w^{(1,2)}, w^{(2)}$ of suitable sizes. 

Given a sequence of data points $(x_k)_{k \ge 0}$, which can represent the discrete-time trajectory of state process for instance, and an initial input $h_0$, an RNN can be used recursively to produce the sequence $(y_k, h_k)_{k \ge 1}$ defined by
\begin{align*}
   h_{k} &= \rho\left(\beta^{(1)} + w^{(1,1)} h_{k-1} + w^{(1,2)} x_{k-1} \right),
    \\
    y_{k} &= \tilde \rho\left(\beta^{(2)} + w^{(2)} h_{k} \right).
\end{align*}
Here, $h_{k+1}$ encodes information that is transmitted from iteration $k$ to iteration $k+1$. This information is produced using previous information $h_{k}$ and the current data point $x_{k}$. It is used to compute the output $y_{k}$ associated with the current input $x_{k}$, and the future information $h_{k+2}$.

Based on the idea of using an architecture in a recurrent way, many generalizations of the above simple neural network have been proposed. We next present one of them.

\subsubsection{Long short-term memory}\label{sec:LSTMdetails} 
One of the most common types of RNN is the long short-term memory (LSTM) neural network~\cite{hochreiter1997long}. The advantage of an LSTM is the ability to deal with the vanishing gradient problem and data with lags of unknown duration. 
An LSTM is composed of a series of units, each of which corresponds to a timestamp, and each unit consists of a cell $\mathfrak{c}$ and three gates: input gate $\mathfrak{i}$, output gate $\mathfrak{o}$, and forget gate $\mathfrak{f}$. Among these components, the cell keeps track of the information received so far, the input gate captures to which extent new input information flows into the cell, the forget gate captures to which extent the existing information remains in the cell, and the output gate controls to which extent the information in the cell will be used to compute the output of the unit. Given a data sequence $(x_k)_{k \geq 0}$ and an initial input $h_0$, the information flows are
\begin{equation}
\label{eq:lstm_gate}
\begin{aligned}
    &\text{forget gate: } \mathfrak{f}_k = \rho_{\text{s}}(W_f x_k + U_f h_{k-1} + b_f),\\
    &\text{input gate: } \mathfrak{i}_k =  \rho_{\text{s}}(W_i x_k + U_i h_{k-1} + b_i), \\
    &\text{ontput gate: } \mathfrak{o}_k =  \rho_{\text{s}}(W_o x_k + U_o h_{k-1} + b_o), \\
    &\text{cell: } \mathfrak{c}_k = \mathfrak{f}_k \odot \mathfrak{c}_{k-1} + \mathfrak{i}_k \odot \rho_{\tanh}(W_c x_k + U_c h_{k-1} + b_c), \\
    &\text{output of the } k^{th} \text{ unit: } h_k = \mathfrak{o}_k \odot  \rho_{\tanh}(\mathfrak{c}_k),
\end{aligned}
\end{equation}
where the operator $\odot$ denotes the Hadamard product, $W_f$, $W_i$, $W_o$, $W_c$,  $U_f$, $U_i$, $U_o$, $U_c$, $b_f$, $b_i$, $b_o$ and $b_c$ are neural network parameters of compatible sizes, and $\rho_{\text{s}}$ and $\rho_{\tanh}$ are activation functions given in \eqref{eq:activation}.

\subsubsection{Expressive power of neural networks}
The first theories about neural networks date back to 1989 by Cybenko \cite{cybenko1989approximation} and by Hornik, Stinchcombe and White \cite{hornik1989multilayer}, concerning the approximation capabilities of feedforward networks within a given function space of interest. Hornik \cite{hornik1991approximation} then extended the results to approximating the function's derivatives, and Leshno, Lin, Pinkus and Schocken \cite{leshno1993multilayer} proved results under arbitrary nonpolynomial activation functions. These results are referred to as universal approximation theorems. See also~\cite{Pinkus:99}.

In the past decade, the mathematical theory has been greatly developed, which is complemented by unprecedented advances in highly parallelizable graphics processing units (GPUs), the introduction of new network architectures, and the development of GPU-enabled algorithms. For instance, in terms of approximation theories, other types of neural networks have been investigated, including RNN \cite{schafer2006recurrent}, convolutional neural networks (CNNs) \cite{zhou2020universality} and graph neural networks (GNNs) \cite{keriven2021universality}. Concerning the expressive power of neural networks,  \cite{hanin2019universal} analyzed it from the depth point view, while \cite{lu2017expressive,hanin2017approximating,park2020minimum} considered a width perspective. Several works tackle the question of how neural networks can tackle the curse of dimensionality; see, \textit{e.g.}, \cite{hanin2019deep,grohs2021lower,jentzen2021proof,hutzenthaler2020proof}. For further discussion on the mathematical theory of deep learning, we refer to, \textit{e.g.}, \cite{berner_grohs_kutyniok_petersen_2022}.

\subsection{Stochastic gradient descent and its variants}
\label{sec:SGD-var} 

The process of adjusting the parameters of a parameterized system, such as a neural network, in order to optimize a loss function is called \emph{training}. Stochastic gradient descent (SGD) is one of the most popular methods to train neural network parameters, for example for the aforementioned FNNs, RNNs and LSTMs.

Consider a generic optimization problem: minimize over $\varphi$,
$$
    J(\varphi) = \mathbb{E}_{\xi \sim \nu}[\mathfrak{L}(\varphi, \xi)],
$$
where $\xi$ follows a distribution $\nu$ and $\mathfrak{L}$ is a loss function. Using a neural network with a given architecture as an approximator for $\varphi$, the goal becomes to minimize over $\theta$,
$$
    J(\theta) = \mathbb{E}_{\xi \sim \nu}[\mathfrak{L}(\varphi_\theta, \xi)].
$$
Even if $\mathfrak{L}$ is known, the loss cannot be computed exactly when $\nu$ is unknown. If one does not have access to $\nu$ but only to samples drawn from $\nu$, one can use SGD, described in Algorithm~\ref{algo:SGD-generic}, which relies on an empirical risk minimization problem,
$$
    J^{S, N}(\theta) = \frac{1}{N} \sum_{i=1}^N \mathfrak{L}(\varphi_\theta, \xi^i),
$$
where $N$ is the number of training samples of $\xi$ and we denote the sample set by $S = (\xi^1,\dots,\xi^N)$.

\begin{algorithm}[htp]
\begin{algorithmic}
\STATE {\bfseries Input: }{An initial parameter $\theta_0\in\Theta$.  A mini-batch size $N_{\text{Batch}}$. A number of iterations $M$. A sequence of learning rates $(\beta_m)_{m = 0, \dots, M-1}$.}
\STATE {\bfseries Output: }{Approximation of $\theta^*$}
  \FOR{$m = 0, 1, 2, \dots, M-1$}
    \STATE Sample a minibatch of $N_\text{Batch}$ samples $S = (\xi^{i})_{i=1,\dots,N_\text{Batch}}$ where $\xi^i$ are i.i.d. drawn from $\nu$\;
    \STATE Compute the gradient $\nabla J^{S, N_\text{Batch}}(\theta_{m})$ \;
    \STATE Update $\theta_{m+1} =  \theta_{m} -\beta_m \nabla J^{S, N_\text{Batch}}(\theta_{m})$ \;
  \ENDFOR
  \STATE {Return $\theta_{M}$}
\end{algorithmic}
\caption{Stochastic Gradient Descent (SGD)\label{algo:SGD-generic}}
\end{algorithm}

SGD is generally used with a moderately large size of mini-batch, which reduces the computational cost of each iteration and can furthermore help escape local minima. In practice, the choice of the learning rate can be crucial to ensure convergence. A popular way to adjust the learning rate is the Adam method~\cite{kingma2014adam}  which is summarized in Algorithm~\ref{algo:ADAM-generic} and can be viewed as an adaptive momentum accelerated SGD.  The computation of the gradient $\nabla J^{S, N}(\theta)$  with respect to $\theta$ can be done automatically by libraries such as {\tt TensorFlow} or {\tt PyTorch}, which perform this computation efficiently by using backpropagation.

\begin{algorithm}[htp]
\begin{algorithmic}
\STATE {\bfseries Input: }{Stepsize $\alpha$. Exponential decay rates for the moment estimates $\beta_1,\beta_2 \in [0,1)$. Initial parameter $\theta_0$. Small parameter for numerical stability $\epsilon$. }
\STATE {\bfseries Output: }{Approximation of $\theta^*$}
  \STATE Initialize first moment vector $\bar M_0$ and second moment vector $\bar V_0$ \;
  \FOR{$m = 0, 1, 2, \dots, M-1$}
    \STATE Sample a minibatch of $N_\text{Batch}$ samples $S = ((\xi^{i})_{i=1,\dots,N_{\text{Batch}}}$ where $x^i$ are i.i.d. drawn from $\nu$\;
    \STATE Compute the gradient $g_m = \nabla J^{S}(\theta_{m})$ \;
    \STATE Update biased first moment estimate: $\bar M_m = \beta_1 \bar M_{m-1} + (1-\beta_1) g_m$ \;
    \STATE Update biased second moment estimate: $\bar V_m = \beta_2 \bar V_{m-1} + (1-\beta_2) g_m^2$ \;
    \STATE Compute biased-corrected first moment estimate: $\hat M_m = \bar M_m / (1 - \beta_1^m)$ \;
    \STATE Compute biased-corrected second moment estimate: $\hat V_m = \bar V_m / (1 - \beta_2^m)$ \;
    \STATE Set $\theta_{m+1} =  \theta_{m} -\alpha \hat M_m / (\sqrt{\hat V_m} + \epsilon)$ \;
  \ENDFOR
  \STATE {Return $\theta_{M}$}
\end{algorithmic}
\caption{ADAM: Adaptive Moment Estimation \label{algo:ADAM-generic}}
\end{algorithm}

\section{Preliminaries on SDDE}\label{supp:SDDE}

In this section, we provide some background on stochastic differential delay equations. Let $\mathscr{C}_{-\delta} = C([-\delta, 0], \RR^d)$ be the Banach space of all $\RR^d$-valued continuous functions defined on $[-\delta,0]$ endowed with the supremum norm
\begin{equation}
   \|y\|_{\mathscr{C}_{-\delta}} = \sup_{-\delta \leq s \leq 0}|y_s|, \;\forall y \in \mathscr{C}_{-\delta}.
\end{equation}
This space is represents portions of trajectories. 
The drift $b$ and volatility $\sigma$ coefficients, and the running $f$ and terminal costs $g$ are deterministic functionals
\begin{equation}
    (b,\sigma, f) : [0,T] \times \mathscr{C}_{-\delta} \times \mc{A} \to (\RR^d, \RR^{d \times m}, \RR); \quad g: \mathscr{C}_{-\delta} \to \RR.
\end{equation}
Denote by $L^2(\Omega, \mathscr{C}_{-\delta})$ the space of all $\MCF$-measurable stochastic processes, \emph{i.e.},
\begin{equation}
    \Omega \ni \omega \to X(\omega) \in \mathscr{C}_{-\delta} \text{ is in } L^2(\Omega, \mathscr{C}_{-\delta}), \text{ iff. } \int_\Omega \|X(\omega)\|_{\mathscr{C}_{-\delta}}^2 \ud \PP(\omega) < \infty.
\end{equation}
The space $L^2(\Omega, \mathscr{C}_{-\delta})$ endowed with the semi-norm $\|X\|_{L^2(\Omega, \mathscr{C}_{-\delta})} = [\int_\Omega \|X(\omega)\|_{\mathscr{C}_{-\delta}}^2 \ud \PP(\omega)]^{1/2}$ is a complete space. We assume that the initial path $\underline \varphi \in L^2(\Omega, \mathscr{C}_{-\delta})$ and is independent of the Brownian motion $W$,  and the existence of a solution $X$ to the SDDE \eqref{def:delay-Xt} is considered in $L^2(\Omega, C([-\delta, T], \RR^d))$. Let $(\mc{F}_t)_{t\geq 0}$ be the filtration supporting $W$ and $\varphi$, and let $C([0,T], L^2(\Omega, \mathscr{C}_{-\delta}))$ be the space of all $L^2$-continuous $\mathscr{C}_{-\delta}$-valued $\MCF_t$-adapted processes $P: [0,T] \ni t \to \underline P_t \in L^2(\Omega, \mathscr{C}_{-\delta})$ endowed with the semi-norm
\begin{equation}
    \|P\|_{C([0,T], L^2(\Omega, \mathscr{C}_{-\delta}))} = \sup_{0 \leq t \leq T} \| \underline P_t \|_{L^2(\Omega, \mathscr{C}_{-\delta})}.
\end{equation}
The trajectory $\underline X_t$ solving the SDDE \eqref{def:delay-Xt} is considered in $C([0,T], L^2(\Omega, \mathscr{C}_{-\delta}))$.

Usually, one requires uniform Lipschitz conditions in the second variable of $b$ and $\sigma$ to ensure the existence and uniqueness of strong solutions to the SDDE \eqref{def:delay-Xt}, that is,
\begin{equation}
    \|(b, \sigma)(t, y_1, \ctrl) - (b, \sigma)(t, y_2, \ctrl) \|_{L^2} \leq L \|y_1 - y_2\|_{L^2(\Omega, \mathscr{C}_{-\delta})}, \; \forall t \in [0, T] \text{ and } y_1, y_2 \in L^2(\Omega, \mathscr{C}_{-\delta}). 
\end{equation}
See detailed analysis in Mohammed's monographs \cite{mohammed1984stochastic,mohammed1998stochastic}. Assumptions on $f$ and $g$ will ensure the expected cost \eqref{def:delay-cost} is finite.

\section{Pseudo-codes of Algorithms}\label{supp:pseudocode}

\begin{algorithm}[h t p]
		\caption{Deep Fictitious Play for Finding Open-loop Nash Equilibrium \label{algo:DFP-open}}
	\begin{algorithmic}[1]
		\REQUIRE $N$ = \# of players, $N_T$ = \# of subintervals on $[0,T]$, $M$ = \# of training paths, $M'$ = \# of out-of-sample paths for final evaluation, $\bm{\alpha}^0 = \{\alpha_{t)n}^{i,0} \in \mc{A} \subset \RR^k, i \in \mc{I}\}_{n=0}^{N_T-1}$ = initial belief, $\bm{X}_0 = \{x_0^i \in \RR^d, i \in \mc{I}\}$ = initial states
		  \STATE Create $N$ separated deep neural networks as described in Eq.~\eqref{def:NN-OL} 
		  \STATE  Generate $M$ discrete sample path of BM: $\bm{W} = \{W_{t_n}^{i} \in \RR^m, i \in \mc{I}\cup \{0\} \}_{n=1}^{N_T} $
		  \STATE $\mt{k} \gets 0$
		  \REPEAT
		  \FORALLP{$i \in \mc{I}$} 
		  \STATE $\mt{k} \gets \mt{k}+1$
		  \STATE (Continue to) Train $i^{th}$ NN with data $\{\bm{X}_0, \bm\alpha^{-i,\mt{k}-1} = \{\alpha_{t_n}^{j,\mt{k}-1}, j \in \mc{I}\setminus \{i\}\}_{n=0}^{N_T-1}, \bm W\}$
		  \STATE Obtain the approximated optimal strategy  $\alpha^{i,\mt{k}}$ and cost $J^i(\alpha^{i,\mt{k}}; \bm\alpha^{-i,\mt{k}-1})$ 
		  \ENDFOR
		  \STATE Collect optimal policies at stage $\mt{k}$: $\bm{\alpha}^\mt{k} \gets (\alpha^{1,\mt{k}}, \ldots, \alpha^{N,\mt{k}})$
		  \STATE Compute relative change of cost $\displaystyle err^\mt{k} = \max_{i \in \mc{I}}\left\{\frac{\abs{J^i(\alpha^{i,\mt{k}}; \bm\alpha^{-i,\mt{k}-1}) - J^i(\alpha^{i,\mt{k}-1}; \bm{\alpha}^{-i,\mt{k}-2})}}{ J^i(\alpha^{i,\mt{k}-1}; \bm{\alpha}^{-i,\mt{k}-2})}\right\}$

		  \UNTIL$err^\mt{k}$ go below a threshold
		  
		\STATE Generate $M'$ out-of-sample paths of BM for final evaluation
		\STATE $\mt{k}' \gets 0$
		\REPEAT 
		\STATE $\mt{k}' \gets \mt{k}'+1$
		\STATE Evaluate $i^{th}$ NN with  \{$\bm X_0$, $\bm{\alpha}^{-i,\mt{k}'-1}$, out-of-sample paths\}, $\forall i \in \mc{I}$
		\STATE Obtain $\alpha^{i,\mt{k}'}$ and $J^{i,\mt{k}'} = J^i(\alpha^{i,\mt{k}'};\bm{\alpha}^{-i, \mt{k}'-1})$ $\forall i \in \mc{I}$
		\UNTIL $J^{i,\mt{k}'}$ converges in $\mt{k}'$, $\forall i \in \mc{I}$
		\RETURN The optimal policy $\alpha^{i,\mt{k}'}$, and the final cost for each player $J^{i,\mt{k}'}$
	\end{algorithmic}
\end{algorithm}

\begin{algorithm}[h t p]
\caption{Deep Fictitious Play for Finding Markovian Nash Equilibrium \label{def_algorithm1}}
    \begin{algorithmic}[1]
	\REQUIRE $N$ = \# of players, $N_T$ = \# of subintervals on $[0,T]$, $\mt{K}$ = \# of total stages in fictitious play, $N_{\text{sample}}$ = \# of sample paths generated for each player at each stage of fictitious play, $N_{\text{SGD\_per\_stage}}$ = \# of SGD steps for each player at each stage, $N_{\text{batch}}$ = batch size per SGD update, $\balpha^0\colon$ the initial policies that are smooth enough
	    \STATE Initialize $N$ deep neural networks to represent $u^{i,0}, i \in \mc{I}$
		\FOR{$\mt{k} \gets 1$ to $\mt{K}$}
		\FORALLP{$i \in \mc{I}$}
		\STATE  Generate $N_\text{sample}$ sample paths $\{\check \bX_{t_n}^{i}\}_{n=0}^{N_T}$ according to \eqref{eq:disc_X_path} and the realized optimal policies $\balpha^{-i, \mt{k}-1}(t_n, \check {\bm X}_{t_n}^{i})$
		\FOR{$\ell \gets 1$ to $N_{\text{SGD}\_\text{per}\_\text{stage}}$}
		    \STATE Update the parameters of the $i^{th}$ neural network one step with $N_{\text{batch}}$ paths using the SGD algorithm (or its variant), based on the loss function \eqref{eq:disc_objective}
		  \ENDFOR
		  \STATE Obtain the approximate optimal policy  $\alpha^{i,\mt{k}}$  according to \eqref{def_alphaast}
		  \ENDFOR
		  \STATE Collect the optimal policies at stage $\mt{k}$: $\bm{\alpha}^\mt{k} \gets (\alpha^{1,\mt{k}}, \ldots, \alpha^{N,\mt{k}})$
		  \ENDFOR
		\RETURN The optimal policy $\balpha^{\mt{K}}$
	\end{algorithmic}
\end{algorithm}

\begin{algorithm}[h t p]
   \caption{The Sig-DFP Algorithm}
   \label{alg:sig-dfp}
\begin{algorithmic}
   \STATE {\bfseries Input:} $b, \sigma, \sigma_0, f, g, \iota$ and $X_0(\omega_i), \{W_{t_n}(\omega_i)\}_{n=0}^{N_T}, \{W^0_{t_n}(\omega_i)\}_{n=0}^{N_T}$ for $i=1,2,\dots, N$; $\mt{K}$: rounds for FP;    $B$: minibatch size; $N_{\text{batch}}$: number of minibatches. 
   \STATE Compute the signatures of $\hat{W}^0_{t_n}(\omega_i)$ for $i=1, \dots, N$, $n=1,\dots, N_T$;
   \STATE Initialize $\hat{\nu}^{(0)}$, $\theta$;
   \FOR{$\mt{k}=1$ {\bfseries to} $\mt{K}$}
   \FOR{$r=1$ {\bfseries to} $N_{\text{batch}}$}
   \STATE Simulate the $r^{th}$ minibatch of $X^{(\mt{k})}(\omega_i)$ using $\hat{\nu}^{(\mt{k}-1)}$ and compute $J_B(\theta, \hat{\nu}^{(\mt{k}-1)})$;
   \STATE Minimize $J_B(\theta, \hat{\nu}^{(\mt{k}-1)})$ over $\theta$, then update $\alpha(\cdot; \theta)$;
   \ENDFOR
   \STATE Simulate $X^{(\mt{k})}(\omega_i)$ with the optimized $\alpha(\cdot; \theta^\ast)$, 
   for $i=1, \dots, N$;
   \STATE Regress $\iota(X^{(\mt{k})}_0(\omega_i), \alpha^{(\mt{k})}_0(\omega_i))$, $\iota(X^{(\mt{k})}_{T/2}(\omega_i), \alpha^{(\mt{k})}_{T/2}(\omega_i))$, $\iota(X^{(\mt{k})}_{T}(\omega_i), \alpha^{(\mt{k})}_T(\omega_i))$ on $S^M(\hat{W}^0_0{\omega_i})$, $S^M(\hat{W}^0_{t_{N_T/2}}(\omega_i))$, $S^M(\hat{W}^0_{T}(\omega_i))$ to get $\tilde l^{(\mt{k})}$; 
   \STATE Update $\hat{l}^{(\mt{k})} = \frac{\mt{k}-1}{\mt{k}}\hat{l}^{(\mt{k}-1)} + \frac{1}{\mt{k}}\tilde l^{(\mt{k})}$;
   \STATE Compute $\hat{\nu}^{(\mt{k})}$ by $\hat{\nu}^{(\mt{k})}_{t_n}(\omega^i) = \langle \hat{l}^{(\mt{k})}, S^M(\hat{W}^0_{t_n}(\omega_i)) \rangle$, for $i=1,2,\dots, N, n=1,\dots, N_T$;
   \ENDFOR
   \STATE {\bfseries Output:} the optimized $\alpha_\varphi^\ast$ and  $\hat {l}^{(\mt{K})}$.
\end{algorithmic}
\end{algorithm}

\section{List of Acronyms}

\begin{longtable}{l|l}
ADAM    & Adaptive Moments Gradient Descent\\
BSDE    & Backward Stochastic Differential Equation \\
2BSDE   & Second Order Backward Stochastic Differential Equation \\
CDC     & Centers for Disease Control and Prevention \\
DNN     & Deep Neural Networks \\
DBDP    & Deep Backward Dynamic Programming \\
Deep BSDE & Deep Backward Stochastic Differential Equation \\
DDPG & Deep Deterministic Policy Gradient \\
DGM     & Deep Galerkin Method \\
DFP     & Deep Fictitious Play \\
DPG     & Deterministic Policy Gradient \\
DPP     & Dynamic Programming Principle \\
FBSDE   & Forward Backward Stochastic Differential Equation \\
FNN     & Feedforward Neural Network \\
FP      & Fokker-Planck \\
GRU     & Gated Recurrent Unit\\
HJB     & Hamilton-Jacobi-Bellman \\
KFP     & Kolmogorov-Fokker-Planck \\
LQ      & Linear-Quadratic \\
LSTM    & Long Short-term Memory \\
MARL    & Multi-agent Reinforcement Learning \\
MC      & Monte Carlo \\
MDP     & Markov Decision Process \\
MFC     & Mean-Field Control \\
MFG     & Mean-Field Game \\
MFMDP   & Mean-Field Markov Decision Process \\
MKV     & McKean-Vlasov \\
MKV FBSDE & McKean-Vlasov Forward Backward Stochastic Differential Equation \\
NE      & Nash Equilibrium \\
NN & Neural Network\\
ODE     & Ordinary Differential Equation \\
PDE     & Partial Differential Equations \\
PINN    & Physics Informed Neural Network \\
ReLu    & Rectified Linear Unit \\
RL      & Reinforcement Learning \\
RNN     & Recurrent Neural Network \\
SC      & Stochastic Control \\
SDE     & Stochastic Differential Equation \\
SDDE    &  Stochastic Differential Delay Equation \\
SDFP    & Scaled Deep Fictitious Play \\
SGD     & Stochastic Gradient Descent \\
Sig-DFP & Signatured Deep Fictitious Play \\
TD      & Temporal-Difference
\end{longtable}

\section{List of Frequently Used Notations}
\begin{longtable}{l|l}
    \hline
    \hline
    \multicolumn{2}{l}{Deep learning related notations}
    \\
    \hline
    $\theta$ & neural network parameters 
    \\
    $\varphi^\theta$ & neural network with parameters $\theta$ \\
    $J(\theta)$ & optimization function for training the parameters $\theta$ \\
    \hline
    \hline
    \multicolumn{2}{l}{Optimal control related notations}
    \\
    \hline
    $x$ & state variable
    \\
    $(X_t)_{t \in [0,T]}$ & state process
    \\
    $(Y_t, Z_t)_{t \in [0,T]}$ & backward and adjoint processes 
    \\
    $(\ctrl_t)_{t \in [0,T]}$, $(\beta_t)_{t \in [0,T]}$ & control process
    \\
    $(W_t)_{t \in [0,T]}$, $(W^0_t)_{t \in [0,T]}$ & Wiener process
    \\
    $d$ & dimension of the state
    \\
    $k$ & dimension of the action
    \\
    $m$ & dimension of the noise
    \\
    $b(t, x, \ctrl)$ & drift coefficient of the state process
    \\
    $\sigma(t, x, \ctrl)$ & diffusion coefficient of the state process
    \\
    $f(t, x, \alpha)$ & instantaneous cost function
    \\
    $g(x)$ & terminal cost function
    \\
    $J(\alpha)$ & total cost function associated with the control $\alpha$
    \\
    $\check J(\theta)$ & total cost function associated with the discretized stochastic control problem 
    \\
    & with neural network parameters $\theta$
    \\
    $(\overline X_t)_{t \in [0,T]}$ & empirical average of $N$-players' states
    \\     
    $\underline X_t= \bigl(\underline X_t(s)\bigr)_{s \in [-\delta, 0]}$ & trajectory of $X_t$ from time $t-\delta$ to $t$ %
    \\
    $T$ & time horizon
    \\
     $\Hess_x u(t, x)$ & Hessian matrix of $u$ with respect to $x$ \\
    $H$ & Hamiltonian function
    \\
    $u(t, x)$ & value function
    \\
    \hline
    \hline
    \multicolumn{2}{l}{Stochastic games related notations}
    \\ 
    \hline
    $X^i_t$ & state process for player $i$\\
    $\alpha^i_t$ & action function for player $i$\\
    $\bm{\alpha} = [\alpha^1, \alpha^2, \ldots, \alpha^N]$ & a collection of all players' strategy profiles \\
    $\bm{\alpha}^{-i} = [\alpha^1, \ldots, \alpha^{i-1}, \alpha^{i+1}, \ldots, \alpha^{N}]$ & the strategy profiles excluding player $i$'s \\
    $\bm{W} = [W^0, W^1, \ldots, W^N]$ & $(N+1)$-vector of $m$-dimensional independent Brownian motions \\
    $\mathbb{F} = \{\MCF_t, 0\leq t \leq T\}$ & the augmented filtration generated by $\bm{W}$ \\
    $J^i(\boldsymbol{\alpha})$ & cost function for player $i$ \\
    $\mathrm{k}$ & stage index in deep fictitious play \\
    $m$ & mean process (or conditional mean if there is common noise)
    \\
    $\mu$ & state distribution process
    \\
    $\mu^N$ & empirical state distribution process
    \\
    $\nu$ & state-action distribution process
    \\
    $\nu^N$ & empirical state-action distribution process

\end{longtable}

\bibliographystyle{myplain}

\bibliography{Reference_surveyMLSCG}

\end{document}